\documentclass[final,1p,times]{elsarticle}

\newcommand{\f}{\frac}

\newcommand{\tbf}{\textbf}

\newcommand{\mbf}{\mathbf}
\newcommand{\mal}{\mathcal}
\newcommand{\mo}{\mathcal{O}}

\newtheorem{theorem}{Theorem}[section]
\newtheorem{remark}{Remark}[section]

\newtheorem{lemma}{Lemma}[section]
\newtheorem{corollary}{Corollary}[section]

\usepackage{amssymb}
\usepackage{latexsym}
\usepackage{hyperref}
\usepackage{color}
\usepackage{amsmath}
\usepackage{multirow}

\journal{       }
\begin{document}

\begin{frontmatter}


\numberwithin{equation}{section}

\title{High order numerical methods based on quadratic spline collocation method and averaged L1 scheme for the variable-order time fractional mobile/immobile diffusion equation}


\author[1]{Xiao Ye}
\author[1]{Jun Liu\corref{Liu}}\ead{liujun@upc.edu.cn}
\author[2]{Bingyin Zhang}
\author[2]{Hongfei Fu}
\author[1]{Yue Liu}
\address[1]{College of Science, China University of Petroleum (East China), Qingdao, Shandong 266580, China}
\address[2]{School of Mathematical Sciences, Ocean University of China, Qingdao, Shandong 266100, China}
\cortext[Liu]{Corresponding author.}

\begin{abstract}
	In this paper, we consider the variable-order time fractional mobile/immobile diffusion (TF-MID) equation in two-dimensional spatial domain, where the fractional order $\alpha(t)$ satisfies $0<\alpha_{*}\leq \alpha(t)\leq \alpha^{*}<1$. We combine the quadratic spline collocation (QSC) method and the $L1^+$ formula to propose a QSC-$L1^+$ scheme. It can be proved that, the QSC-$L1^+$ scheme is unconditionally stable and convergent with $\mo(\tau^{\min{\{3-\alpha^*-\alpha(0),2\}}} + \Delta x^{2}+\Delta y^{2})$, where $\tau$, $\Delta x$ and $\Delta y$ are the temporal and spatial step sizes, respectively. With some proper assumptions on $\alpha(t)$, the QSC-$L1^+$ scheme has second temporal convergence order even on the uniform mesh, without any restrictions on the solution of the equation. We further construct a novel alternating direction implicit (ADI) framework to develop an ADI-QSC-$L1^+$ scheme, which has the same unconditionally stability and convergence orders. In addition, a fast implementation for the ADI-QSC-$L1^+$ scheme based on the exponential-sum-approximation (ESA) technique is proposed. Moreover, we also introduce the optimal QSC method to improve the spatial convergence to fourth-order. Numerical experiments are attached to support the theoretical analysis, and to demonstrate the effectiveness of the proposed schemes.
\end{abstract}


\end{frontmatter}


\section{Introduction}
Over the past several decades, the fractional partial differential equations (FPDEs) have attracted more and more attention as a tool for
modeling various physical phenomena with memory or hereditary properties, such as damping laws,
diffusion processes and viscoelastic behavior and so on \cite{Wuli-2,Wuli-4,Wuli-1}.
Recent studies showed that the variable order fractional PDEs are more powerful tools for modeling many multiphysics phenomena, where
the properties of the materials or systems evolve with time \cite{Wang-1,Sun-2}.
As an important class of variable-order FPDEs, the variable-order time fractional mobile/immobile diffusion (TF-MID) equation
describes the transport characteristics of particle in the fluid, and provides a more realistic model for
solute diffusion transport in heterogeneous porous media \cite{M/IM-1,Wang-2}.

In recent years, various numerical methods for FPDEs have been proposed, such as
finite difference methods \cite{FD-4,FD-2,FD-1,L1-1}, finite element methods \cite{FE-1,FE-3},
finite volume methods \cite{FV-1,FV-2,FV-3}, spectral methods \cite{S-3,ADI-Li-1}, quadratic spline collocation (QSC) method \cite{QSC-7,C-1,C-2} and so on. The QSC method gives an approximation to the solution of the original differential equations in the quadratic spline space. Since the QSC method employs smoother basis functions, and needs less degrees of freedom than some classical methods for the same number of grid points, then the QSC method results in algebraic systems of relatively smaller scale. The QSC method, as well as its optimal version, have been widely applied for various kinds of integer-order PDEs \cite{KE,PSoper2,PSoper1}.

For time FPDEs,  piecewise interpolation based numerical method is one of the main strategies for discretization. Sun and Wu \cite{L1-1} derived
the $L1$ scheme for the fractional diffusion-wave equation.
Lin and Xu \cite{L1-Xu-1} constructed a stable $L1$ scheme for time fractional diffusion equation. Li et al. \cite{L1-3} used the linearized $L1$-Galerkin finite element
method to solve the multidimensional nonlinear time-fractional Schr\"{o}dinger equation. Alikhanov \cite{L21-1} proposed the $L_{2-1\sigma}$ formula for Caputo fractional
derivative to achieve second convergence order. 
Lv and Xu \cite{L2-1} proposed the $L2$ scheme based on parabolic interpolation, to achieve high-order accuracy.
Quan and Wang \cite{L2-2} established the energy stability of high-order $L2$-type schemes for time fractional phase-field equations.
Shen et al. \cite{L1+1} developed $L1^+$ scheme for constant order Caputo fractional derivative on suitably graded meshes to achieve second-order convergence in time.
Ji et al. \cite{L1+2} employed $L1^+$ scheme to solve the time-fractional molecular beam epitaxial
models with constant order Caputo fractional derivative. The $L1^+$ scheme can achieve second-order convergence for functions with enough regularity just based on the piecewise linear interpolation.

Variable-order fractional differential operators, like their constant-order counterparts, are nonlocal and weakly singular.
But the construction of numerical discretization is more complicated, and the numerical analysis is more difficult.
Zeng et al. \cite{Zeng-1} proposed spectral collocation methods for variable-order fractional advection-diffusion
equation. 
Zheng and Wang \cite{VO-Zheng-1} proposed an $L1$ scheme for a hidden-memory variable-order space-time fractional
diffusion equation.
Du et al. \cite{Sun-1} developed a temporal second-order finite difference scheme for the variable-order
time-fractional wave equation. We proposed and analyzed a first-order numerical method based on the classical $L1$ scheme, for variable-order TF-MID equation with variable diffusive coefficients \cite{Vo-Liu-Fu}.

For multi-dimensional problems, the alternating direction implicit (ADI) method is an efficient solution strategy, and it can divide the multi-dimensional problem into a series of independent one-dimensional problems. ADI methods have also been widely used for variants of FPDEs.
Ran and Zhang \cite{ADI-Ran-1} proposed compact ADI difference schemes to solve a class of spatial fractional nonlinear
damped wave equations in two space dimensions. 
Qiu et al. \cite{ADI-Qiu-1} presented the ADI Galerkin finite element method to solve the distributed-order
time-fractional mobile/immobile equation in two dimensions. We proposed the QSC method in the ADI framework for two-dimensional space fractional diffusion equation \cite{C-1}.

The historical dependence of the time fractional operators results in the computational complexity $\mathcal{O}(n^2)$, with $n$ the number of the time levels, which is much more expensive than that of the integer-order operator. In order to reduce computational cost, many methods have been proposed to accelerate the evaluation of the fractional derivatives. Jiang et al. \cite{Fast-SOE-1} employed sum-of-exponentials (SOE) technique
with the $L1$ scheme for constant-order time fractional derivatives, which reduced the computational cost to $\mathcal{O}(n\log^{2}n)$.
Liao et al. \cite{L1+2} applied the SOE technique to speed up the evaluation of the $L1^+$ formula for constant-order fractional derivatives. For variable-order FPDEs, Zhang et al. \cite{Fast-ESA-1} approached the singular kernel in Caputo fractional derivatives by the exponential-sum-approximation (ESA) technique, which reduced the computational cost to
$\mathcal{O}(n\log^{2}n)$. Based on the ESA technique, they developed a fast temporal second-order scheme with $L$2-1$_{\sigma}$ formula
in \cite{Fast-ESA-2}. 

In this paper, we first combine the $L1^+$ formula in time discretization with the QSC method in space discretization to propose the QSC-$L1^+$ scheme, for solving the variable-order TF-MID equation in two-dimensional space domain. Such a scheme is a one-step method, and easy to implement. We will prove that the scheme is unconditionally stable and convergent with the order $\mo(\tau^{\min{\{3-\alpha^*-\alpha(0),2\}}} + \Delta x^{2}+\Delta y^{2}))$. Then, we design a novel ADI framework to produce an ADI-QSC-$L1^+$ scheme, where the error caused by alternating direction is much smaller than that caused by $L1^+$ formula. Numerical tests show that the ADI-QSC-$L1^+$ scheme preserves almost the same observation error as the QSC-$L1^+$ scheme.
Furthermore, the fast computation based on the ESA technique with properly chosen parameters for $L1^+$ formula of the variable-order differential operator is constructed, which leads to
the ADI-QSC-F$L1^+$ scheme, and it can reduce the computational cost and the memory requirement effectively. In addition, we employ the optimal QSC method in space by introducing proper perturbations to get the optimal ADI-QSC-F$L1^+$ scheme with fourth-order convergence in space, which results in the numerical solution with a desired accuracy with much less spatial meshes.

The outline of this paper is as follows. In Section \ref{sect:QSCL1+}, we propose the QSC-$L1^+$ scheme for the variable-order TF-MID equation.
Then the unconditional stability and convergence are proved in Section \ref{sect:stability}. In Section \ref{sect:ADIQSCL1+}, we introduce the ADI framework to develop the ADI-QSC-$L1^+$ scheme, and anaylze the unconditional stability and convergence. In Section \ref{sect:speed}, we respectively consider the fast implementation based on the ESA technique along the time direction and the optimal QSC method in the space domain, to reduce the computational cost.
Numerical experiments are presented in Section \ref{sect:numerical}, to verify the theoretical results of
the proposed schemes. Finally, conclusions are given in Section \ref{sect:conclusions}.

Throughout this paper, we use $C_i$ and $Q_i$ to denote positive constants which are independent of the temporal and spatial step sizes.

\section{Variable-order TF-MID equation and the QSC-$L1^+$ scheme}
\label{sect:QSCL1+}
In this section, we consider the numerical solution of
the following two-dimensional variable-order TF-MID equation, which is used to model the anomalously diffusive transport \cite{Wang-2}.
\begin{equation}\label{equation1}
   u_{t}(x,y, t)+{ }_{0}^{C} D_{t}^{\alpha(t)} u(x,y, t)
   = \kappa \mathcal{L} u(x,y, t)+f(x,y, t), \quad (x,y, t) \in \Omega \times(0, T],
\end{equation}
subjecting to the initial condition
\begin{equation}\label{equation2}
    u(x,y,0)=u^{0}(x,y), \quad (x,y) \in \bar{\Omega}=\Omega\cup\partial \Omega,
\end{equation}
and the boundary condition
\begin{equation}\label{equation3}
    u(x,y,t)= 0, \quad (x,y,t) \in \partial \Omega \times(0, T],  
\end{equation}
where constant $\kappa >0$ is the diffusion coefficient, $f$ is a given source function, $u^{0}$ is the initial function. $\Omega=(x_L,x_R)\times(y_L,y_R)$ is a rectangular domain,
and $\partial \Omega$ is the boundary.
$\mathcal{L}$ is the spatial elliptic operator with
\begin{equation*}
 \begin{aligned}
  \mathcal{L} u=\frac{\partial^{2}u}{\partial x^{2}}+\frac{\partial^{2}u}{\partial y^{2}}.
 \end{aligned}
\end{equation*}
$[0,T]$ is the time interval,
and $\alpha(t) \in C[0,T]$ is the variable time fractional order which satisfies the following conditions
\begin{equation}  \label{alpt1}
 \begin{aligned} 
  0<\alpha_{*}\leq \alpha(t)\leq \alpha^{*}<1,\ t\in[0,T],\ \ \lim_{t \rightarrow 0+}(\alpha(t)-\alpha(0)) \ln t \text { exists }.
 \end{aligned}
\end{equation}
The variable-order Caputo fractional
differential operator ${}_{0}^{C} D_{t}^{\alpha(t)}$ is usually used to describe the subdiffusive transport of a large amount of solute particles , which is defined as
\begin{equation*}
 \begin{aligned}
  { }^{C}_0 D_{t}^{\alpha(t)} u(x,y,t):=\int_{0}^{t} \omega_{1-\alpha(t)}(t-s)\partial_{s}u(x,y,s)ds,
 \end{aligned}
\end{equation*}
where
\begin{equation} \label{kernel}
 \begin{aligned}
  \omega_{1-\beta}(t):=\frac{t^{-\beta}}{\Gamma(1-\beta)}.
 \end{aligned}
\end{equation}
The term $u_{t}(x,y, t)$ in (\ref{equation1}) describes the the Fickian diffusive transport of the remaining portion of the total solute mass.
The function $u(x,y,t)$ is to be determined. 

For illustrating the singularity of the solution at the initial time, one can define the
weighted Banach space involving time $C_{\mu}^{m}((0, T] ; \mathcal{X})$ with $m \geq 2,0 \leq \mu<1$ \cite{Zheng-19},
\begin{equation*}
	\begin{aligned}
	 &C_{\mu}^{m}((0, T] ; \mathcal{X}) 
	 :=\left\{v \in C^{1}([0, T] ; \mathcal{X}):\|v\|_{C_{\mu}^{m}((0, T] ; \mathcal{X})}<\infty\right\}, \\
	 &\|v\|_{C_{\mu}^{m}((0, T] ; \mathcal{X})} 
	 :=\|v\|_{C^{1}([0, T] ; \mathcal{X})}+\sum_{l=2}^{m} \sup _{t \in(0, T]} t^{l-1-\mu}\left\|\frac{\partial^{l} v}{\partial t^{l}}\right\|_{\mathcal{X}}.
	\end{aligned}
\end{equation*}
In addition, the eigenfunctions $\left\{\varphi_{i}\right\}_{i=1}^{\infty}$ of the Sturm-Liouville problem
\begin{equation*}
	\begin{aligned}
	\mathcal{L} \varphi_{i}(x, y)=\lambda_{i} \varphi_{i}(x, y),\ (x, y) \in \Omega ; \quad \varphi_{i}(x, y)=0,\ (x, y) \in \partial \Omega
	\end{aligned}
\end{equation*}
form an orthonormal basis in $L^{2}(\Omega)$. The eigenvalues $\left\{\lambda_{i}\right\}_{i=1}^{\infty}$ are positive and nondecreasing
that tend to $\infty$ with $i$. By the theory of sectorial operators, we can define the fractional Sobolev space
\begin{equation*}
	\begin{aligned}
		\breve{H}^{\gamma}(\Omega):=\left\{v \in L^{2}(\Omega):|v|_{\breve{H}^{\gamma}}^{2}
		:=\sum_{i=1}^{\infty} \lambda_{i}^{\gamma}\left(v, \varphi_{i}\right)^{2}<\infty\right\},
	\end{aligned}
\end{equation*}
with the norm $\|v\|_{\breve{H}^{\gamma}}:=\left(\|v\|_{L^{2}}^{2}+|v|_{\breve{H}^{\gamma}}^{2}\right)^{1/2}$. 
Furthermore, $\breve{H}^{\gamma}(\Omega)$ is a subspace of the fractional Sobolev space ${H}^{\gamma}(\Omega)$ that can be characterized by \cite{H-gamma}
\begin{equation*}
	\begin{aligned}
		\breve{H}^{\gamma}(\Omega)
		=\left\{v \in H^{\gamma}(\Omega)
		: \mathcal{L}^{s} v(x, y)=0,(x, y) \in \partial \Omega, s<\gamma / 2\right\},
	\end{aligned}
\end{equation*}
and the norms $|v|_{\breve{H}^{\gamma}}$ and $|v|_{{H}^{\gamma}}$ are equivalent in $\breve{H}^{\gamma}$.

With conditions (\ref{alpt1}) and suitable assumptions on the data, the following important
lemma ensures the regularity and well-posedness of the solution of model (\ref{equation1})--(\ref{equation3}).
\begin{lemma}[\cite{M/IM-1, Zheng-1}] \label{ZhengZe}
	If condition \eqref{alpt1} holds and $u^{0} \in \check{H}^{\gamma+2}$, $f \in H^{d}\left([0, T] ; \breve{H}^{\gamma}(\Omega)\right)$
for $\gamma>1 / 2$ and $d>1 / 2$, then the variable-order TF-MID model \eqref{equation1}--\eqref{equation3} have a unique solution
$u \in {C}^{1}\left([0, T] ; \breve{H}^{\gamma}(\Omega)\right)$ and
\begin{equation} \label{ZhengZe1}
	\begin{aligned}
	&\|u\|_{{C}\left([0, T] ; \breve{H}^{s}(\Omega)\right)} 
	\leq 
	Q\left(\left\|u^{0}\right\|_{\breve{H}^{s}(\Omega)}+\|f\|_{H^{d}\left([0, T] ; \breve{H}^{\max \{s-2,0\}}(\Omega)\right)}\right), \\
	&\|u\|_{{C}^{1}\left([0, T] ; \breve{H}^{s}(\Omega)\right)} 
	\leq 
	Q\left(\left\|u^{0}\right\|_{\breve{H}^{s+2}(\Omega)}+\|f\|_{H^{d}\left([0, T] ; \breve{H}^{s}(\Omega)\right)}\right)
	\end{aligned}
\end{equation}
for $0 \leq s \leq \gamma$. Here $ Q=Q\left(\alpha^{*},\|\alpha\|_{{C}[0, T]}, T, d\right)$.

Moreover, suppose that $u_{0} \in \check{H}^{\gamma+6}$, $f \in H^{1}\left([0, T] ; \check{H}^{s+4}\right) \cap H^{2}\left([0, T] ; \check{H}^{s+2}\right) \cap H^{3}\left([0, T] ; \check{H}^{s}\right)$ for $s\geq 0$, $\alpha \in C^{2}[0, T]$, and \eqref{alpt1} holds.
If  $\alpha(0)>0$, 
we have  $u \in C^{3}\left((0, T] ; \check{H}^{\gamma}(0, L)\right) \cap C_{1-\alpha(0)}^{3}\left((0, T];\check{H}^{\gamma}(0, L)\right)$  and
\begin{equation} \label{ZhengZe2}
	\begin{aligned}
		\|u\|_{C_{1-\alpha(0)}^{3}\left((0, T] ; \check{H}^{\gamma}(0, L)\right)} 
		\leq 
		C_{0} \left( 
		\left\|u_{0}\right\|_{\check{H}^{\gamma+6}(0, L)}
		+\|f\|_{H^{1}\left(\check{H}^{s+4}\right)}
		+\|f\|_{H^{2}\left(\check{H}^{s+2}\right)}
		+\|f\|_{H^{3}\left(\check{H}^{s}\right)} 
		\right).
	\end{aligned}
\end{equation}

\end{lemma}

Next, we will consider $L1^{+}$ discretization in time and QSC discretization in space.
\subsection{The $L1^+$ formula for time discretization}
Given a positive integer $N$, we define a uniform partition
on the time interval $[0, T ]$ by $t_n = n\tau$, for $n = 0, 1, . . . , N$ with $\tau = \frac{T}{N}$.
For function $v( t)$, piecewise linear interpolant on the temporal
mesh is denoted by $\Pi v( t)$,
and we define $\theta v\left( t\right)=v\left(t\right)-\Pi v\left(t\right)$. By the Taylor's expansion with integral remainder, we can obtain
\begin{equation*}
	\begin{aligned}
		\theta v\left(t\right)
		=\int_{t_{n-1}}^{t} (t-s) \partial_{s}^{2} v\left(s\right) \mathrm{d} s
		+\frac{1}{\tau}\int_{t_{n-1}}^{t_n} \left(t-t_{n-1}\right)\left(t_{n}-s\right) \partial_{s}^{2} v\left(s\right) \mathrm{d} s,
		\quad t_{n-1} \leq t \leq t_{n},\ 1 \leq n \leq N.
	\end{aligned}
\end{equation*}
Based on Lemma \ref{ZhengZe}, we have $\left\|\frac{\partial^{2} v}{\partial t^{2}}\right\|_{\mathcal{X}} \leq C_{1}t^{-\alpha(0)}$.
Therefore, it can be verified that
\begin{equation} \label{theta-v}
	\begin{aligned}
		\left|\theta v\left(t\right)\right| 
		\leq 
		C_{2} \tau\left(t_{n}^{1-\alpha(0)}-t_{n-1}^{1-\alpha(0)}\right), \quad t_{n-1} \leq t \leq t_{n},\ 1 \leq n \leq N.
	\end{aligned}
\end{equation}

We denote by $v^{n} $ the approximation of $v( t)$ at time instant $t=t_{n}$.
Let $\mathfrak{T} =\left\{v^{n},\ n=0,1,\cdots,N\right\}$ be a temporal grid function space, we define
\begin{equation*}
 \begin{aligned}
  \boldsymbol{\delta}_{t} v^{n-\frac{1}{2}} = \frac{v^n-v^{n-1}}{\tau}
  \quad \text{and} \quad
  v^{n-\frac{1}{2}}=\frac{v^n+v^{n-1}}{2}.
 \end{aligned}
\end{equation*}
Then, taking the mean
value of the $L$1 discretization of ${}_{0}^{C}D_{t}^{\alpha(t)} v( t)$ over $[t_{n-1}, t_{n}]$, that is
\begin{equation}\label{L1+-1}
	\begin{aligned} 
		\frac{1}{\tau}\int_{t_{n-1}}^{t_{n}} {}_{0}^{C}D_{t}^{\alpha(t)} v( t)dt 
		=\frac{1}{\tau}\int_{t_{n-1}}^{t_{n}} {}_{0}^{C}D_{t}^{\tilde{\alpha}_{n}}v( t)dt +r_{1,n} ,
	\end{aligned}
\end{equation}
where $\tilde{\alpha}_{n} := \alpha_{n-\frac{1}{2}}$, for $n=1,2,\cdots,N$. Based on trapezoidal formula, we can verify that $r_{1,n}=\mo\left(\tau^{2}\right)$. 
For completeness, we give a detailed proof in \ref{r1n}. 
Moreover, we have
\begin{equation} \label{L1+operaer}
 \begin{aligned}
 &\frac{1}{\tau}\int_{t_{n-1}}^{t_{n}} {}_{0}^{C}D_{t}^{\tilde{\alpha}_{n}}v( t)dt \\
 &=\frac{1}{\tau}\int_{t_{n-1}}^{t_{n}}\int_{0}^{t}\omega _{1-\tilde{\alpha}_{n}}(t-s)\partial_{s} \Pi v( s)dsdt +r_{2,n}\\
 &=\frac{1}{\tau}\int_{t_{n-1}}^{t_{n}}\sum_{k=1}^{n-1}\int_{t_{k-1}}^{t_{k}}
     \frac{(t-s)^{-\tilde{\alpha}_{n}}}{\Gamma(1-\tilde{\alpha}_{n})}\cdot \frac{v^{k} -v^{k-1} }{\tau}dsdt
  +\frac{1}{\tau}\int_{t_{n-1}}^{t_{n}}\int_{t_{n-1}}^{t} \frac{(t-s)^{-\tilde{\alpha}_{n}}}{\Gamma (1-\tilde{\alpha}_{n})}\cdot \frac{v^{n} -v^{n-1} }{\tau}dsdt +r_{2,n}.
 \end{aligned}
\end{equation}
The local truncation error $r_{2,n}=\mo\left(\tau^{2}t_{n}^{-\tilde{\alpha}_{n}-\alpha(0)}\right)$,
based on the truncation error analysis in \cite{L1+2,L1+1,WangL1+}. The detailed proof will be given in \ref{r2n}.
By some fundamental calculations, the integration in (\ref{L1+-1}) can be expressed as
\begin{equation}\label{L1+operaer-2}
 \begin{aligned}
  \frac{1}{\tau}\int_{t_{n-1}}^{t_{n}} {}_{0}^{C}D_{t}^{\alpha(t)} v( t)dt
   &=\sum_{k=1}^{n} a_{n-k+1}^{(n)}\left(v^{k} -v^{k-1} \right)+r_{1,n}+r_{2,n}\\
  &:=\bar{\boldsymbol{\delta}_{t}}^{\tilde{\alpha}_{n}} v^{n-\frac{1}{2}} +r_{1,n}+r_{2,n},\quad n = 1, 2, \cdots ,N,
 \end{aligned}
\end{equation}
where
\begin{equation} \label{L1+-cof}
 \begin{aligned}
  a_{n-k+1}^{(n)}=\frac{1}{\tau^{2}} \int_{t_{n-1}}^{t_{n}} \int_{t_{k-1}}^{\min \left\{t, t_{k}\right\}} \omega_{1-\tilde{\alpha}_{n}}(t-s)dsdt,\quad k=1,2,\cdots,n.
 \end{aligned}
\end{equation}

The discretization (\ref{L1+operaer-2}) for the variable-order Caputo time fractional derivative ${}_{0}^{C}D_{t}^{\alpha(t)} v( t)$,
with the coefficients (\ref{L1+-cof}), is called $L1^+$ formula.
The coefficients $a_{n-k+1}^{(n)}$ for $k=1,2,\cdots,n$, satisfy the following lemma, which will be used in the numerical analysis below.
\begin{lemma}[\cite{L1+2,L1+1}]\label{lemma:L1+CO}
At time instant $t=t_{n}$, the coefficients $\left\{a_{n-k+1}^{(n)},k=1,2,\cdots,n\right\}$ of L1$^+$ scheme satisfy\\
 \begin{equation*}
  \begin{aligned}
  a_{2}^{(n)}>a_{3}^{(n)}>\cdots>a_{n}^{(n)}>0.
  \end{aligned}
 \end{equation*}
\end{lemma}

\subsection{The QSC method for space discretization}
Let $M_x$ and $M_y$ be two positive integers. We define the uniform spatial partitions of $[x_L,x_R]$ and $[y_L,y_R]$ as
\begin{equation*}
 \begin{aligned}
  \triangle_x:=\{x_L=x_{0}<x_{1}<\ldots<x_{M_x}=x_{R}\},\quad \triangle_y:=\{y_L=y_{0}<y_{1}<\ldots<y_{M_y}=y_{R}\},
 \end{aligned}
\end{equation*}
respectively, with corresponding mesh sizes $\Delta x=\frac{x_R-x_L}{M_x}$ and $\Delta y=\frac{y_R-y_L}{M_y}$.
Furthermore, let $\triangle :=\triangle_x \times \triangle_y$ be the mesh partition of $\bar{\Omega}$.

Define the quadratic splines space along each spatial direction as
\begin{equation*}
 \begin{aligned}
  \mathcal{V}_{x}&:=\left\{v\in C^{1}\left(x_{L},x_{R}\right),\ \left.v\right|_{\left[x_{i-1},x_{i}\right]}\in\mbf{P}_{2}\big(\triangle_{x}\big),i=1,2,\ldots,M_{x}\right\}, \\
  \mathcal{V}_{y}&:=\left\{v \in C^{1}\left(y_{L}, y_{R}\right),\ \left.v\right|_{\left[y_{j-1},y_{j}\right]}\in\mbf{P}_{2}\big(\triangle_{y}\big),j=1,2,\ldots, M_{y}\right\},
 \end{aligned}
\end{equation*}
where $\mbf{P}_{2}(\cdot)$ represent the set of quadratic polynomials in a single variable. Besides, let
\begin{equation*}
 \begin{aligned}
  \mathcal{V}_{x}^{0}:=\left\{v \in \mathcal{V}_{x},\ v\left(x_{L}\right)=v\left(x_{R}\right)=0\right\},
  \quad \mathcal{V}_{y}^{0}:=\left\{v \in \mathcal{V}_{y},\ v\left(y_{L}\right)=v\left(y_{R}\right)=0\right\},
 \end{aligned}
\end{equation*}
and denote by $\mathcal{V}^{0}:=\mathcal{V}_{x}^{0} \otimes \mathcal{V}_{y}^{0}$ the space of piecewise biquadratic polynomials with respect to the
spatial partition $\triangle$, which satisfy the homogeneous Dirichlet boundary conditions $(\ref{equation3})$.

Now we consider the basis functions for the space $\mathcal{V}^{0}$. First, let
\begin{equation*}
 \begin{aligned}
  \phi(x)=\frac{1}{2}
   \left\{
    \begin{aligned}
  &x^{2},                   &0\leq x \leq 1, \\
  &-2(x-1)^{2}+2(x-1)+1,    &1\leq x \leq 2, \\
  &(3-x)^{2},               &2\leq x \leq 3, \\
  &0,                       &\text{elsewhere}.
    \end{aligned}
  \right.
 \end{aligned}
\end{equation*}
We define the quadratic B-splines
\begin{equation} \label{basisf}
 \begin{aligned}
  \phi_{j}(x)=\phi\left(\frac{x-x_{L}}{\Delta x}-j+2\right), \quad j=0,1, \cdots, M_{x}+1,
 \end{aligned}
\end{equation}
and choose $\left\{\phi_{j}(x), \ j=0, \cdots ,M_{x}+1 \right\}$ as the basis function of $\mathcal{V}_{x}^{0}$.
Similarly, we can define the basis functions $\left\{\phi_{j}(y), \ j=0, \cdots ,M_{y}+1 \right\}$ for
$\mathcal{V}_{y}^{0}$ with replacing the variable $x$ by the variable $y$. Then, the basis functions of $\mathcal{V}^0$ can be defined as
the tensor product of the basis functions for the spaces $\mathcal{V}_{x}^{0}$ and $\mathcal{V}_{y}^{0}$.
Thus, the quadratic spline solution $u_{h}^{n} \in \mathcal{V}^{0}$ of the model (\ref{equation1}) can be represented as
\begin{equation} \label{qsc}
 \begin{aligned}
  u_{h}^{n}(x, y)=\sum_{i=0}^{M_{x}+1} \sum_{j=0}^{M_{y}+1} c_{ij}^{n} \phi_{i}(x) \phi_{j}(y),\quad n=1,\cdots ,N,
 \end{aligned}
\end{equation}
where the coefficients $c_{ij}^{n}$ are degrees of freedom (DOFs). In order to determine the DOFs,
we define the midpoints of $\triangle$ as
\begin{equation*}
 \begin{aligned}
  \xi =\left\{\left(\xi _{i}^{x}, \xi _{j}^{y}\right),\ i=1,2, \ldots, M_{x} ; \  j=1,2, \cdots, M_{y}\right\},
 \end{aligned}
\end{equation*}
where $\left\{\xi_{i}^{x}=\frac{1}{2}\left(x_{i-1}+x_{i}\right), \  i=1,2, \cdots, M_{x}\right\}$ and
$\left\{\xi_{j}^{y}=\frac{1}{2}(y_{j-1}+y_{j}), \  j=1,2, \cdots, M_{y}\right\}$.
Denote $\partial{\xi}= \left\{\left(\xi _{i}^{x}, \xi _{j}^{y}\right),\ i\in \{0, M_{x}+1\}\right.$ or $\left.j\in\{0,M_{y}+1\}\right\}$,
where $\xi_{0}^{x}=x_{0}$, $\xi_{M_{x}+1}^{x}=x_{M_{x}}$
and
$\xi_{0}^{y}=y_{0}$, $\xi_{M_{y}+1}^{y}=y_{M_{y}}$
are  boundary points of $\triangle_x$ and $\triangle_y$, respectively. Thus, we choose $\bar{\xi}=\xi \cup \partial\xi$ as the collocation points.
For convenience, we also define index sets $\Lambda=\left\{ (i,j),\ (\xi _{i}^{x}, \xi _{j}^{y}) \in \xi \right\}$,
$\partial\Lambda=\left\{ (i,j),\ (\xi _{i}^{x}, \xi _{j}^{y}) \in \partial \xi \right\}$ and $\bar{\Lambda}=\Lambda \cup \partial\Lambda$.

 Based on the definition $(\ref{basisf})$,
we can get the following lemma by some fundamental calculations.

\begin{lemma}[\cite{Basis function}]\label{lemma:basisf}
$\textit{(}I\textit{)}.$ For the basis function $\phi_{0}(x)$, we have
\begin{equation*}
 \begin{aligned}
  \phi_{0}\left(\xi_{i}^{x}\right)=\frac{1}{8}
   \left\{
    \begin{aligned}
    4,&\quad i=0,\\
    1,&\quad i=1,\\
    0,&\quad \text{else,}
   \end{aligned}
  \right.
  \quad \phi_{0}^{\prime\prime}\left(\xi_{i}^{x}\right)=\frac{1}{{\Delta x}^2}
   \left\{
    \begin{aligned}
     1, & \quad i=1, \\
     0, & \quad i=2,3,\ldots,M_{x}+1.
    \end{aligned}
  \right.\\
   \end{aligned}
\end{equation*}
$\textit{(}II\textit{)}.$ For the basis functions $\phi_{j}(x)$ with $j = 1,\ldots,M_{x}$, we have
\begin{equation*}
 \begin{aligned}
  \phi_{j}\left(\xi_{0}^{x}\right)=\frac{1}{8}
  \left\{
   \begin{aligned}
    4, &\quad j=1, \\
    0, &\quad \text {else,}
  \end{aligned}\quad \right.
 \end{aligned}
\quad
 \begin{aligned}
  \phi_{j}\left(\xi_{M_{x}+1}^{x}\right)=\frac{1}{8}
  \left\{
   \begin{aligned}
    4, &\quad j=M_{x}, \\
    0, &\quad \text {else,}
  \end{aligned}\quad \right.
 \end{aligned}
\end{equation*}
and for $i = 1, \ldots , M_{x}$,
\begin{equation*}
 \begin{aligned}
 \phi_{j}\left(\xi_{i}^{x}\right)=\frac{1}{8}
  \left\{
   \begin{aligned}
    1, &\quad |i-j|=1, \\
    6, &\quad i=j, \\
    0, &\quad \text {else,}
  \end{aligned}\quad \right.
  \phi_{j}^{\prime\prime}\left(\xi_{i}^{x}\right)=\frac{1}{{\Delta x}^2}
  \left\{
   \begin{aligned}
    1, &\quad |i-j|=1, \\
    -2, &\quad i=j, \\
    0, &\quad \text {else.}
  \end{aligned} \right.
   \end{aligned}
\end{equation*}
$\textit{(}III\textit{)}.$ For the basis function $\phi_{M_{x}+1}(x)$, we have
\begin{equation*}
 \begin{aligned}
  \phi_{M_{x}+1}\left(\xi_{i}^{x}\right)=\frac{1}{8}
   \left\{
    \begin{aligned}
     4, & \quad i=M_{x}+1, \\
     1, & \quad i=M_{x}, \\
     0, & \quad \text{else,}
    \end{aligned}
  \right.
  \quad \phi_{M_{x}+1}^{\prime\prime}\left(\xi_{i}^{x}\right)=\frac{1}{{\Delta x}^2}
  \left\{
   \begin{aligned}
     1, & \quad i=M_{x}, \\
     0, & \quad i=0,1,\ldots,M_{x}-1.
    \end{aligned}
  \right.\\
 \end{aligned}
\end{equation*}
\end{lemma}

The properties in Lemma \ref{lemma:basisf} also holds for the basis functions $\left\{\phi_{j}(y), \ j=0, \cdots ,M_{y}+1 \right\}$.
Taking the collocation points into (\ref{qsc}), together with Lemma \ref{lemma:basisf} , we can get for $k=0, \cdots ,M_{x}+1$, $l=0, \cdots ,M_{y}+1$ that
\begin{equation*}
 \begin{aligned}
  u_{h}^{n}\left(\xi_{k}^{x}, \xi_{l}^{y}\right)
   =\sum_{i=\max \{k-1,0\}}^{\min \left\{k+1, M_{x}+1\right\}} \sum_{j=\max \{l-1,0\}}^{\min \left\{l+1, M_{y}+1\right\}}
    c_{i j} \phi_{i}\left(\xi_{k}^{x}\right) \phi_{j}\left(\xi_{l}^{y}\right)
  :=\boldsymbol{\theta}_{x}\boldsymbol{\theta}_{y} c_{kl}^{n},\\
 \end{aligned}
\end{equation*}
\begin{equation*}
 \begin{aligned}
  \frac{\partial^{2}}{\partial x^{2}} u_{h}^{n}\left(\xi_{k}^{x}, \xi_{l}^{y} \right)
   =\sum_{i=\max \{k-1,0\}}^{\min \left\{k+1, M_{x}+1\right\}} \sum_{j=\max \{l-1,0\}}^{\min \left\{l+1, M_{y}+1\right\}}
     c_{i j} \phi_{i}^{\prime\prime}\left(\xi_{k}^{x}\right) \phi_{j}\left(\xi_{l}^{y}\right)
   := \boldsymbol{\eta}_{x} \boldsymbol{\theta}_{y} c_{kl}^{n},
 \end{aligned}
\end{equation*}
\begin{equation*}
 \begin{aligned}
  \frac{\partial^{2}}{\partial y^{2}} u_{h}^{n}\left(\xi_{k}^{x}, \xi_{l}^{y}\right)
   =\sum_{i=\max \{k-1,0\}}^{\min \left\{k+1, M_{x}+1\right\}} \sum_{j=\max \{l-1,0\}}^{\min \left\{l+1, M_{y}+1\right\}}
    c_{i j} \phi_{i}\left(\xi_{k}^{x}\right) \phi_{j}^{\prime\prime}\left(\xi_{l}^{y}\right)
   := \boldsymbol{\eta}_{y} \boldsymbol{\theta}_{x}c_{kl}^{n},
 \end{aligned}
\end{equation*}
where the operators $\boldsymbol{\theta}_{x}$ and $\boldsymbol{\eta}_{x}$ are defined as
\begin{equation} \label{theta x}
  \boldsymbol{\theta}_{x} c_{k, l}^{n}=\frac{1}{8}
  \left\{
   \begin{aligned}
    &4c_{0, l}^{n}+4c_{1,l}^{n}, & k=0, \\
    &c_{k-1, l}^{n}+6 c_{k, l}^{n}+c_{k+1, l}^{n}, &k=1,2,\cdots, M_{x},\\
    &4c_{M_{x}, l}^{n}+4c_{M_{x}+1, l}^{n}, & k=M_{x}+1,
   \end{aligned}
  \right.
\end{equation}
\begin{equation}\label{delta x}
 \boldsymbol{\eta}_{x} c_{k, l}^{n}=\frac{1}{{\Delta x}^2}
  \left\{
  \begin{aligned}
  &0, & k=0, M_{x}+1, \\
  &(c_{k-1,l}^{n}-2c_{k,l}^{n}+c_{k+1,l}^{n}),\quad &k=1,2,\cdots, M_{x}.\\
  \end{aligned}
  \right.
\end{equation}
Moreover, we define
\begin{equation} \label{1stD}
 \begin{aligned}
  \boldsymbol{\vartheta}_{x} c_{k, l}^{n}=\frac{1}{{\Delta x}}\left( c_{k,l}^{n}-c_{k-1,l}^{n}\right),\quad &k=1,2,\cdots, M_{x}+1,
 \end{aligned}
\end{equation}
then we have
\begin{equation*}
 \begin{aligned}
  \boldsymbol{\eta}_{x} c_{k, l}^{n} = \frac{1}{{\Delta x}} \left( \boldsymbol{\vartheta}_{x} c_{k+1,l}^{n}-\boldsymbol{\vartheta}_{x} c_{k,l}^{n}\right),\quad &k=1,2,\cdots, M_{x}.
 \end{aligned}
\end{equation*}

In addition, the operators $\boldsymbol{\theta}_{y}$, $\boldsymbol{\eta}_{y}$ and $\boldsymbol{\vartheta}_{y} $ are defined along the $y$ direction,
and have the similar expressions as $\boldsymbol{\theta}_{x}$, $\boldsymbol{\eta}_{x}$ and $\boldsymbol{\vartheta}_{x} $, just with $\Delta  x$ replaced by $\Delta  y$.
Next, we will consider the full discretization scheme based on the above approximations.
\subsection{The QSC-L1$^+$ scheme}
We average model (\ref{equation1}) over the time subinterval
$[t_{n-1},t_{n}]$ to get
\begin{equation} \label{L1+}
 \begin{aligned}
   \frac{1}{\tau}\int_{t_{n-1}}^{t_{n}}u_{t}(x,y,t)dt+\frac{1}{\tau} \int_{t_{n-1}}^{t_{n}} {}_{0}^{C}D_{t}^{\alpha(t)}u(x,y,t)dt
   =\frac{1}{\tau}\int_{t_{n-1}}^{t_{n}}\kappa \mathcal{L} u(x,y,t)dt+\frac{1}{\tau} \int_{t_{n-1}}^{t_{n}}f(x,y,t)dt.
 \end{aligned}
\end{equation}
It can be verified for the first term on the left hand side of (\ref{L1+}) that
\begin{equation} \label{Left-1}
 \begin{aligned}
  \frac{1}{\tau}\int_{t_{n-1}}^{t_{n}}u_{t}(x,y,t)dt=\frac{u^n(x,y)-u^{n-1}(x,y)}{\tau}=\boldsymbol{\delta}_{t} u^{n-\frac{1}{2}}(x,y).
 \end{aligned}
\end{equation}
For the first term on right hand side of (\ref{L1+}), we have
\begin{equation} \label{Right-1}
 \begin{aligned}
  \frac{1}{\tau} \int_{t_{n-1}}^{t_{n}} \kappa \mathcal{L} u(x, y, t) dt
  =\kappa \mathcal{L} u^{n-\frac{1}{2}}(x, y)
   +r_{3,n},
 \end{aligned}
\end{equation}
where
\begin{equation*}
 \begin{aligned}
  r_{3,n}
  &=\frac{1}{\tau} \int_{t_{n-1}}^{t_{n}} \frac{\kappa}{2} \mathcal{L}\big[u\left(x, y, t \right)-u^{n-\frac{1}{2}}(x, y)\big] \\
  &=\frac{1}{\tau} \int_{t_{n-1}}^{t_{n}} \frac{\kappa}{2} \mathcal{L}\big[u\left(x, y, t \right)-\Pi u(x, y, t)\big]
  =\frac{1}{\tau} \int_{t_{n-1}}^{t_{n}} \frac{\kappa}{2} \mathcal{L}\theta u\left(x, y, t \right)  ,
 \end{aligned}
\end{equation*}
and satisfies
\begin{equation*}
 \begin{aligned}
  |r_{3,n}|
  \leq
  \frac{\kappa }{2\tau} \int_{t_{n-1}}^{t_{n}} |\mathcal{L}\theta u\left(x, y, t\right)|  dt
  \leq 
  C_{3} \tau\left(t_{n}^{1-\alpha(0)}-t_{n-1}^{1-\alpha(0)}\right)
 \end{aligned}
\end{equation*}
Similarly, for the second term on the right hand side of (\ref{L1+}), we have
\begin{equation} \label{Right-2}
 \begin{aligned}
  \frac{1}{\tau} \int_{t_{n-1}}^{t_{n}} f(x,y,t) dt
  =f^{n-\frac{1}{2}}(x, y)
   +r_{4,n},
 \end{aligned}
\end{equation}
where
\begin{equation*}
 \begin{aligned}
  r_{4,n}
  =\frac{1}{\tau} \int_{t_{n-1}}^{t_{n}} \frac{1}{2} f_{t t}\left(x, y, \rho_{1}\right)\left(t-t_{n}\right)\left(t-t_{n-1}\right) dt,\quad \rho_{1} \in (t_{n-1},t_{n}),
 \end{aligned}
\end{equation*}
and satisfies
\begin{equation*}
 \begin{aligned}
  |r_{4,n}|
  \leq
  \frac{C_{4}}{2\tau} \int_{t_{n-1}}^{t_{n}}\left(t-t_{n}\right)\left(t-t_{n-1}\right) dt
  =\mo\left(\tau^{2}\right),
 \end{aligned}
\end{equation*}
with $C_{4}$ is the bound of $f_{tt}(x,y,t)$.

Based on equations (\ref{Left-1})-(\ref{Right-2}), together with the discretization (\ref{L1+operaer-2}), equation (\ref{L1+}) can be rewritten as
\begin{equation} \label{L1+3}
 \begin{aligned}
  \boldsymbol{\delta}_{t} u^{n-\frac{1}{2}}(x, y)+ \bar{\boldsymbol{\delta}_{t}}^{\tilde{\alpha}_{n}}u^{n-\frac{1}{2}}(x, y)
  =\kappa \mathcal{L} u^{n-\frac{1}{2}}(x, y)+f^{n-\frac{1}{2}}(x, y)+R^{n},
 \end{aligned}
\end{equation}
where
\begin{equation} \label{LTE}
 \begin{aligned}
  R^{n}=r_{1,n}+r_{2,n}+r_{3,n}+r_{4,n}=\mo\left(\tau^{2}t_{n}^{-\tilde{\alpha}_n-\alpha(0)}+\tau\left(t_{n}^{1-\alpha(0)}-t_{n-1}^{1-\alpha(0)}\right)+\tau^{2}\right).
 \end{aligned}
\end{equation}

In order to find numerical solution of model (\ref{equation1}) in the quadratic splines space, we take $u_{h}^{n}(x, y)$ with the form (\ref{qsc}) into (\ref{L1+3})
and drop truncation errors,
which lead to
\begin{equation} \label{L1+2}
 \begin{aligned}
   &\sum_{i=0}^{M_{x}+1} \sum_{j=0}^{M_{y}+1} \boldsymbol{\delta}_{t}c_{ij}^{n-\frac{1}{2}} \phi_{i}(x) \phi_{j}(y)
   +\sum_{i=0}^{M_{x}+1} \sum_{j=0}^{M_{y}+1} \bar{\boldsymbol{\delta}_{t}}^{\tilde{\alpha}_n} c_{ij}^{n-\frac{1}{2}} \phi_{i}(x) \phi_{j}(y)\\
  &=\kappa \sum_{i=0}^{M_{x}+1} \sum_{j=0}^{M_{y}+1} c_{ij}^{n-\frac{1}{2}}
   \left[\phi_{i}^{\prime \prime}(x) \phi_{j}(y)+\phi_{i}(x) \phi_{j}^{\prime \prime}(y)\right] +f^{n-\frac{1}{2}}(x,y),\quad (x,y) \in \Omega,\ 1\leq n \leq N,
 \end{aligned}
\end{equation}
with the initial condition
\begin{equation} \label{ini2}
 \begin{aligned}
  \sum_{i=0}^{M_{x}+1} \sum_{j=0}^{M_{y}+1} c_{ij}^{0} \phi_{i}(x) \phi_{j}(y)=u^{0}(x,y), \quad (x,y) \in \bar{\Omega},
 \end{aligned}
\end{equation}
and the boundary condition
\begin{equation} \label{bound2}
 \begin{aligned}
  \sum_{i=0}^{M_{x}+1} \sum_{j=0}^{M_{y}+1} c_{ij}^{n} \phi_{i}(x) \phi_{j}(y)=0, \quad (x,y) \in \partial\Omega, \ 1\leq n \leq N.
 \end{aligned}
\end{equation}
Now taking the collocation point $(\xi_{i}^{x}, \xi_{j}^{y})$ for $(i,j)\in \bar{\Lambda}$ into (\ref{L1+2})-(\ref{bound2}), respectively,
we directly get the QSC-$L1^+$ scheme,
\begin{equation} \label{QSC-L1+F}
 \begin{aligned}
   \boldsymbol{\delta}_{t} \boldsymbol{\theta}_{x} \boldsymbol{\theta}_{y} c_{ij}^{n-\frac{1}{2}}+\bar{\boldsymbol{\delta}_{t}}^{\tilde{\alpha}_n} \boldsymbol{\theta}_{x} \boldsymbol{\theta}_{y} c_{ij}^{n}
     =\kappa(\boldsymbol{\eta}_{x} \boldsymbol{\theta}_{y} +\boldsymbol{\eta}_{y} \boldsymbol{\theta}_{x})c_{ij}^{n-\frac{1}{2}}+f_{ij}^{n-\frac{1}{2}},
      \quad (i,j)\in \bar{\Lambda},\quad 1\leq n \leq N,
 \end{aligned}
\end{equation}
with the initial condition
\begin{equation} \label{ini2F}
 \begin{aligned}
    \boldsymbol{\theta}_{x} \boldsymbol{\theta}_{y} c_{ij}^{0}= u_{ij}^{0}, \quad (i,j)\in \bar{\Lambda},
 \end{aligned}
\end{equation}
and the boundary condition
\begin{equation} \label{bound2F}
 \begin{aligned}
    \boldsymbol{\theta}_{x} \boldsymbol{\theta}_{y} c_{ij}^{n}=0, \quad (i,j)\in \partial{\Lambda}, \quad 1\leq n \leq N.
 \end{aligned}
\end{equation}

Next, we will analyze the stability and convergence of the scheme.
\section{Numerical analysis of the QSC-$L1^+$ scheme}
\label{sect:stability}
Before the numerical analysis, we need some definitions of the inner products and norms.
Define $\mathcal{M}_h = \big\{u,\ u=\{u_{i,j},\ (i,j)\in \bar{\Lambda}\} \big\}$ as the spatial grid function space with respect to the partition $\triangle$,
and $\mathring{\mathcal{M}}_{h}  = \{ u\in\mathcal{M}_h,\ u_{ij}=0 \ \text{for}\  (i,j)\in \partial{\Lambda} \}$.
For any $u,v \in \mathring{\mathcal{M}}_{h}$, we define the discrete inner product
\begin{equation*}
 \begin{aligned}
  &(u, v):=\Delta x \Delta y \sum_{i=0}^{M_{x}+1 } \sum_{j=0}^{M_{y}+1 } u_{ij} v_{ij},
  \quad
  \left(\boldsymbol{\vartheta}_{x} u, \boldsymbol{\vartheta}_{x} v\right):=\Delta x \Delta y \sum_{i=1}^{M_{x} +1}
       \sum_{j=0}^{M_{y}+1 } \left(\boldsymbol{\vartheta}_{x} u_{ij}\right) \left(\boldsymbol{\vartheta}_{x} v_{ij}\right),\\
  &\left(\boldsymbol{\vartheta}_{y} u, \boldsymbol{\vartheta}_{y} v\right):=\Delta x \Delta y \sum_{i=0}^{M_{x}+1 } \sum_{j=1}^{M_{y}+1} \left(\boldsymbol{\vartheta}_{y} u_{ij} \right) \left(\boldsymbol{\vartheta}_{y} v_{ij}\right).
 \end{aligned}
\end{equation*}
Thus, the corresponding discrete norms can be obtained as
\begin{equation*}
 \begin{aligned}
  \|u\|:=\sqrt{\left(u,u\right)},
  \quad\left|u\right|_{1x}:=\sqrt{\left(\boldsymbol{\vartheta}_{x} u, \boldsymbol{\vartheta}_{x} u\right)}\ ,
  \quad\left|u\right|_{1y}:=\sqrt{\left(\boldsymbol{\vartheta}_{y} u, \boldsymbol{\vartheta}_{y} u\right)}\ .
 \end{aligned}
\end{equation*}

Recalling that
\begin{equation*}
 \begin{aligned}
  \bar{\boldsymbol{\delta}_{t}}^{\tilde{\alpha}_n} v^{n-\f{1}{2}}=\sum_{k=1}^{n} a_{n-k+1}^{(n)}\left(v^{k}-v^{k-1}\right),
 \end{aligned}
\end{equation*}
we can reformulate scheme (\ref{QSC-L1+F}) as
\begin{equation} \label{QSC L1+2}
 \begin{aligned}
   (1+\tau a_{1}^{(n)})\left(\boldsymbol{\theta}_{x} \boldsymbol{\theta}_{y} c_{ij}^{n}
   -\boldsymbol{\theta}_{x} \boldsymbol{\theta}_{y} c_{ij}^{n-1}\right)
   +\tau \sum_{k=1}^{n-1} a_{n-k+1}^{(n)}\left(\boldsymbol{\theta}_{x} \boldsymbol{\theta}_{y} c_{ij}^{k}
   -\boldsymbol{\theta}_{x} \boldsymbol{\theta}_{y} c_{ij}^{k-1}\right)
   =\tau \kappa(\boldsymbol{\eta}_{x} \boldsymbol{\theta}_{y} +\boldsymbol{\eta}_{y} \boldsymbol{\theta}_{x})c_{ij}^{n-\frac{1}{2}}+\tau f_{ij}^{n-\frac{1}{2}},
 \end{aligned}
\end{equation}
for $(i,j)\in \bar{\Lambda}$ and $1\leq n \leq N$. For convenience, we define the coefficients of (\ref{QSC L1+2}) uniformly as
\begin{equation*}
 \left\{
 \begin{aligned}
  &b_{1}^{(n)}=1+\tau a_{1}^{(n)}, \\
  &b_{n-k+1}^{(n)}=\tau a_{n-k+1}^{(n)}, \quad k=1,\ 2,\ \cdots,\ n-1,
 \end{aligned}
 \right.
\end{equation*}
then the QSC-$L1^+$ scheme (\ref{QSC-L1+F}) can be further rewritten as
\begin{equation} \label{new QSC L1+}
 \begin{aligned}
  \sum_{k=1}^{n} b_{n-k+1}^{(n)}\left(\boldsymbol{\theta}_{x} \boldsymbol{\theta}_{y} c_{ij}^{k}-\boldsymbol{\theta}_{x} \boldsymbol{\theta}_{y} c_{ij}^{k-1}\right)
  =\tau \kappa \left(\boldsymbol{\eta}_{x} \boldsymbol{\theta}_{y}+\boldsymbol{\eta}_{y} \boldsymbol{\theta}_{x}\right) c_{ij}^{n-\frac{1}{2}}
   +\tau f_{ij}^{n-\frac{1}{2}},\quad  (i,j)\in \bar{\Lambda},\  1\leq n \leq N,
 \end{aligned}
\end{equation}
Based on Lemma \ref{lemma:L1+CO}, the new coefficients $\left\{b_{n-k+1}^{(n)},k=1,2,\cdots,n\right\}$ in scheme (\ref{new QSC L1+}) satisfy the following lemma.
\begin{lemma}\label{lemma:bnk}
At time instant $t = t_{n}$, if $\tau\leq1$, the coefficients $\left\{b_{n-k+1}^{(n)},k=1,2,\cdots,n\right\}$ in scheme \eqref{new QSC L1+} satisfy
\begin{equation*}
 \begin{aligned}
  b_{1}^{(n)}>b_{2}^{(n)}>\cdots>b_{n}^{(n)}>0.
 \end{aligned}
\end{equation*}
\end{lemma}
{\bf Proof.}
According to the definition of $a_{n-k+1}^{(n)}$, we have
\begin{equation*}
 \begin{aligned}
  a_{1}^{(n)}
  &=\frac{1}{\tau^{2}} \int_{t_{n-1}}^{t_{n}} \int_{t_{n-1}}^{t} \omega_{1-\tilde{\alpha}_n}(t-s) dsdt
=\frac{1}{\tau^{2}} \omega_{3-\tilde{\alpha}_n}(\tau)
=\frac{\tau^{-\tilde{\alpha}_n}}{\Gamma\left(3-\tilde{\alpha}_n\right)},
 \end{aligned}
\end{equation*}
and
\begin{equation*}
 \begin{aligned}
  a_{2}^{(n)}
=\frac{1}{\tau^{2}} \int_{t_{n-1}}^{t_{n}} \int_{t_{n-2}}^{t_{n-1}} \omega_{1-\tilde{\alpha}_n}(t-s) d s d t
=\frac{\tau^{-\tilde{\alpha}_n}}{\Gamma\left(3-\tilde{\alpha}_n\right)}\left(2^{2-\tilde{\alpha}_n}-2\right).
 \end{aligned}
\end{equation*}
Thus, we have
\begin{equation*}
 \begin{aligned}
  b_{1}^{(n)}-b_{2}^{(n)}
 =1+\tau (a_{1}^{(n)}-a_{2}^{(n)})
 =1+\frac{\tau^{1-\tilde{\alpha}_n}}{\Gamma\left(3-\tilde{\alpha}_n\right)}\left(3-2^{2-\tilde{\alpha}_n}\right).
 \end{aligned}
\end{equation*}
When $\tilde{\alpha}_n \in\left[2-\frac{ln3}{ln2},1\right)$, we can get $3-2^{2-\tilde{\alpha}_n} \geq 0$, which verifies that $b_{1}^{(n)}>b_{2}^{(n)}$.
When $\tilde{\alpha}_n \in \left(0,2-\frac{ln3}{ln2}\right)$, we have $-1<3-2^{2-\tilde{\alpha}_n}<0 $. Since $\Gamma\left(3-\tilde{\alpha}_n\right)\geq1$, we
have $\frac{\tau^{1-\tilde{\alpha}_n}}{\Gamma\left(3-\tilde{\alpha}_n\right)}\leq1$, with $\tau\leq1$, which leads to $b_{1}^{(n)}>b_{2}^{(n)}$.
The monotonicity of series $\big\{ b_{k}^{(n)},\ k=2,3,\cdots,n \big\}$ can be found in Lemma \ref{lemma:L1+CO}.\hfill
$\blacksquare$

Lemma \ref{lemma:bnk} means that all the coefficients of the QSC-$L1^+$ scheme (\ref{new QSC L1+}) are monotonic, and this property plays an important role in the following numerical analysis.
\subsection{Auxiliary lemmas}
To proceed with the analysis of stability, we need some auxiliary lemmas. We first investigate some properties of the coefficients in the QSC-$L1^+$ scheme (\ref{new QSC L1+}),
which are exhibited in the following lemmas.
\begin{lemma}\label{lemma:bn}
Suppose that $\alpha^{\prime}(t) \leq 0$, and $\alpha^{\prime}(t)$ is uniformly bounded for $0\leq t\leq T$, then for any fixed $n$ with $ 2\leq n\leq N$, we have
\begin{equation*}
 b_{n-k}^{(n)} \leq\left(1+C_{5} \tau\right) b_{n-k}^{(n-1)} ,\quad k=1, \cdots , n-1,
\end{equation*}
\end{lemma}
where $C_{5}$ is a positive constant.\\
{\bf Proof.}
The proof is generally divided into two parts. In the first part, we consider the case $k=1, \cdots , n-2$, and in the second part, we consider the case $k=n-1$.

(I) For $k=1, \cdots , n-2$, we have $t_{k+1}\leq t_{n-1}$, and
\begin{equation} \label{equationB0}
 \begin{aligned}
  b_{n-k}^{(n)}&=\tau a_{n-k}^{(n)}=\frac{1}{\tau} \int_{t_{n-1}}^{t_{n}} \int_{t_{k}}^{t_{k+1}} \omega_{1-\tilde{\alpha}_n}(t-s)dsdt \\
               &=\frac{1}{\tau}\left[\omega_{3-\tilde{\alpha}_n}\left(t_{n}-t_{k}\right)-\omega_{3-\tilde{\alpha}_n}\left(t_{n-1}-t_{k}\right)
                 -\omega_{3-\tilde{\alpha}_n}\left(t_{n}-t_{k+1}\right)+\omega_{3-\tilde{\alpha}_n}\left(t_{n-1}-t_{k+1}\right)\right]\\
               &=\frac{\tau^{2-\tilde{\alpha}_n}}{\tau\Gamma (3-\tilde{\alpha}_n)}\left[(n-k)^{2-\tilde{\alpha}_n}-2(n-k-1)^{2-\tilde{\alpha}_n}+(n-k-2)^{2-\tilde{\alpha}_n}\right].
 \end{aligned}
\end{equation}
Then the quotient of $b_{n-k}^{(n)}$ and $b_{n-k}^{(n-1)}$ can be simplified as
\begin{equation*}
 \begin{aligned}
  \frac{b_{n-k}^{(n)}}{b_{n-k}^{(n-1)}}
   =\cdot \frac{\Gamma \left(3-\tilde{\alpha}_{n-1}\right)}{\Gamma\left(3-\tilde{\alpha}_n\right)}
    \cdot \frac{t_{n-k}^{2-\tilde{\alpha}_n}-2t_{n-k-1}^{2-\tilde{\alpha}_n}+t_{n-k-2}^{2-\tilde{\alpha}_n}}
               {t_{n-k}^{2-\tilde{\alpha}_{n-1}}-2t_{n-k-1}^{2-\tilde{\alpha}_{n-1}}+t_{n-k-2}^{2-\tilde{\alpha}_{n-1}}}
   :=A^{(n)}\cdot B^{(n)}.
 \end{aligned}
\end{equation*}

We first investigate the quantity $A^{(n)}$. Notice the fact that $0< \alpha(t)\leq  \alpha^{*} <1$, and $\Gamma(x)$ is a increasing and
differentiable function on the interval $[3- \alpha^{*},3]$. Denote $\Gamma_{*} = \min_{3- \alpha^{*} \leq s\leq 3}\left|\Gamma(s)\right|$
and $\Gamma_{*}^{\prime} = \max_{3- \alpha^{*} \leq s\leq 3}|\Gamma^{\prime} (s)|$. Using Taylor expansion, we have
\begin{equation*}
 \begin{aligned}
  A^{(n)} =\frac{\Gamma\left(3-\tilde{\alpha}_n\right)
      +\Gamma^{\prime}\left(\xi_{n}\right)\left(\tilde{\alpha}_n-\tilde{\alpha}_{n-1}\right)}{\Gamma\left(3-\tilde{\alpha}_n\right)}
     = 1+\frac{\Gamma^{\prime}\left(\xi_{n}\right)}{\Gamma\left(3-\tilde{\alpha}_n\right)} \cdot \alpha^{\prime}(\eta_{n}) \cdot \tau
    \leq 1+\tau \left\|\alpha'\right\|_{\infty} \cdot \frac{\Gamma_{*}^{\prime}}{\Gamma_{*}}
    \leq 1+C_{6} \tau,
 \end{aligned}
\end{equation*}
where $\xi_{n} \in (3-\tilde{\alpha}_{n-1},3-\tilde{\alpha}_{n})$ and $\eta_{n} \in (t_{n-1}, t_n)$.
Next, for the quantity $B^{(n)}$, it can be verified that
\begin{equation} \label{equationB2}
 \begin{aligned}
  B^{(n)}
   &=\frac{\Big[t_{n-k}^{2-\tilde{\alpha}_{n}}-t_{n-k-1}^{2-\tilde{\alpha}_{n}}\Big]-\left[t_{n-k-1}^{2-\tilde{\alpha}_{n}}-t_{n-k-2}^{2-\tilde{\alpha}_{n}}\right]}
          {\Big[t_{n-k}^{2-\tilde{\alpha}_{n-1}}-t_{n-k-1}^{2-\tilde{\alpha}_{n-1}}\Big]-\Big[t_{n-k-1}^{2-\tilde{\alpha}_{n-1}}-t_{n-k-2}^{2-\tilde{\alpha}_{n-1}}\Big]}\\
   &=\frac{{\left(2-\tilde{\alpha}_{n}\right)\left(1-\tilde{\alpha}_{n}\right)}}
          {{\left(2-\tilde{\alpha}_{n-1}\right)\left(1-\tilde{\alpha}_{n-1}\right)}}
     \cdot
     \frac{\int_{t_{n-k-1}}^{t_{n-k}} \int_{x-\tau}^{x} s^{-\tilde{\alpha}_{n}}dsdx}
          {\int_{t_{n-k-1}}^{t_{n-k}} \int_{x-\tau}^{x} s^{-\tilde{\alpha}_{n-1}}dsdx}
   :=B_{1}^{(n)} \cdot B_{2}^{(n)}.
 \end{aligned}
\end{equation}

For the term $B_{1}^{(n)}$ in (\ref{equationB2}), since $\alpha^{\prime}(t)$ is bounded
and $0< \alpha(t)\leq  \alpha^{*} <1$,
we can obtain the following estimate
\begin{equation} \label{equationB6}
 \begin{aligned}
  B_{1}^{(n)}&=\frac{{\left(2-\tilde{\alpha}_{n}\right)\left(1-\tilde{\alpha}_{n}\right)}}
          {{\left(2-\tilde{\alpha}_{n-1}\right)\left(1-\tilde{\alpha}_{n-1}\right)}}
  =\left(1+\frac{\tilde{\alpha}_{n-1}-\tilde{\alpha}_{n}}{2-\tilde{\alpha}_{n-1}}\right)
    \cdot
    \left(1+\frac{\tilde{\alpha}_{n-1}-\tilde{\alpha}_{n}}{1-\tilde{\alpha}_{n-1}}\right)\\
  &\leq
    \left(1+\frac{\tau \left\|\alpha'\right\|_{\infty} }{2-\tilde{\alpha}_{n-1}}\right)
    \cdot
    \left(1+\frac{\tau \left\|\alpha'\right\|_{\infty} }{1-\tilde{\alpha}_{n-1}}\right)
  \leq 1+C_{7}\tau.
 \end{aligned}
\end{equation}

For the term $B_{2}^{(n)}$ in (\ref{equationB2}), we define a continuous auxiliary function
\begin{equation*}
 \begin{aligned}
  h_{1}(z)=s^{-z},\quad s>0.
 \end{aligned}
\end{equation*}
Next, we discuss $B_{2}^{(n)}$ separately according to the value of $s$.
If $s \leq 1$, since the function $h_{1}(z)$ is increasing and $\alpha(t)$ is decreasing, we can get $h_{1}(\tilde{\alpha}_{n}) \leq h_{1}(\tilde{\alpha}_{n-1})$.
Thus, we have
\begin{equation} \label{equationB5}
 \begin{aligned}
  B_{2}^{(n)}
  =\frac{\int_{t_{n-k-1}}^{t_{n-k}} \int_{x-\tau}^{x} s^{-\tilde{\alpha}_{n}}dsdx}
             {\int_{t_{n-k-1}}^{t_{n-k}} \int_{x-\tau}^{x} s^{-\tilde{\alpha}_{n-1}}dsdx}
  \leq
    \frac{\int_{t_{n-k-1}}^{t_{n-k}} \int_{x-\tau}^{x} s^{-\tilde{\alpha}_{n-1}}dsdx}
         {\int_{t_{n-k-1}}^{t_{n-k}} \int_{x-\tau}^{x} s^{-\tilde{\alpha}_{n-1}}dsdx}
  = 1.
 \end{aligned}
\end{equation}
If $s>1$, we take the derivative of $h_{1}(z)$,
\begin{equation*}
 \begin{aligned}
  &h_{1}^{\prime}(z)=-\ln s \cdot s^{-z}.
 \end{aligned}
\end{equation*}
Since $h_{1}^{\prime}(z)$ is increasing when $z>0$, we can get
\begin{equation} \label{equationh}
 \begin{aligned}
  h_{1}\left(\tilde{\alpha}_{n}\right)-h_{1}\left(\tilde{\alpha}_{n-1}\right)
  =h_{1}^{\prime}(\gamma_{n})\left(\tilde{\alpha}_{n}-\tilde{\alpha}_{n-1}\right)
  =\ln s\cdot\left(\tilde{\alpha}_{n-1}-\tilde{\alpha}_{n}\right) s^{-\gamma_{n}}
   \leq \ln s \cdot\left\|\alpha^{\prime}\right\|_{\infty} \cdot \tau \cdot s^{-\tilde{\alpha}_{n}},
 \end{aligned}
\end{equation}
where $\gamma_{n} \in (\tilde{\alpha}_{n}, \tilde{\alpha}_{n-1})$.
Based on (\ref{equationh}), we have
\begin{equation*}
 \begin{aligned}
   B_{2}^{(n)}
   &=\frac{\int_{t_{n-k-1}}^{t_{n-k}} \int_{x-\tau}^{x} s^{-\tilde{\alpha}_{n}}-s^{-\tilde{\alpha}_{n-1}}dsdx
          +\int_{t_{n-k-1}}^{t_{n-k}} \int_{x-\tau}^{x} s^{-\tilde{\alpha}_{n-1}}dsdx}
         {\int_{t_{n-k-1}}^{t_{n-k}} \int_{x-\tau}^{x} s^{-\tilde{\alpha}_{n-1}}dsdx}\\
   &\leq
   1+\frac{\ln T \cdot\left\|\alpha^{\prime}\right\|_{\infty} \cdot \tau \cdot \int_{t_{n-k-1}}^{t_{n-k}} \int_{x-\tau}^{x} s^{-\tilde{\alpha}_{n}}dsdx}
          {\int_{t_{n-k-1}}^{t_{n-k}} \int_{x-\tau}^{x} s^{-\tilde{\alpha}_{n-1}}dsdx}
   = 1+C_{8} \tau B_{2}^{(n)},
  \end{aligned}
\end{equation*}
which leads to
\begin{equation} \label{equationB4}
 \begin{aligned}
  B_{2}^{(n)}
    \leq \frac{1}{1-C_{8} \tau}
    \leq 1+C_{9}\tau.
 \end{aligned}
\end{equation}

Combining with (\ref{equationB6}), (\ref{equationB5}) and(\ref{equationB4}), we can obtain
\begin{equation*}
 \begin{aligned}
  B^{(n)} \leq 1+C_{10}\tau.
 \end{aligned}
\end{equation*}

(II) For $k=n-1$, we aim to prove $a_{1}^{(n)} \leq\left(1+C_{5} \tau\right) a_{1}^{(n-1)}$, that is
\begin{equation*}
 \begin{aligned}
  \frac{a_{1}^{(n)}}{a_{1}^{(n-1)}}
  =\tau^{\tilde{\alpha}_{n-1}-\tilde{\alpha}_{n}} \cdot \frac{\Gamma\left(3-\tilde{\alpha}_{n-1}\right)}{\Gamma\left(3- \tilde{\alpha}_{n}\right)}
  \leq A^{(n)}
   \leq 1+C_{6} \tau,
 \end{aligned}
\end{equation*}
which leads to $b_{1}^{(n)} \leq\left(1+C_{5} \tau\right) b_{1}^{(n-1)}$ by the definition of $b_{1}^{(n)} $, and the proof is completed.

\hfill
$\blacksquare$

\begin{lemma}\label{lemma:andown}
For $n\geq 2$, there exists a positive constant $C_{11}$, such that
\begin{equation*}
 a_{n}^{(n)} \geq \frac{T^{-\tilde{\alpha}_n}}{\Gamma(1-\tilde{\alpha}_n)}\geq C_{11},
\end{equation*}
where $C_{11}$ is a constant.
\end{lemma}
{\bf Proof.}
For $n\geq 2$, we have $t_{n-1} \geq t_{1}$. By the definition of $a_{n}^{(n)}$,
\begin{equation*}
 \begin{aligned}
  a_{n}^{(n)}
  =\frac{1}{\tau^{2}} \int_{t_{n-1}}^{t_{n}} \int_{t_{0}}^{t_{1}} \omega_{1-\tilde{\alpha}_n}(t-s) dsdt.
 \end{aligned}
\end{equation*}
Then, by the monotonicity of $\omega_{1-\tilde{\alpha}_n}(t)$, we have
\begin{equation*}
 \begin{aligned}
  a_{n}^{(n)}
  \geq \frac{1}{\tau^{2}} \int_{t_{n-1}}^{t_{n}} \omega_{1-\tilde{\alpha}_n}(t) \tau dt
  =\frac{1}{\tau} \int_{t_{n-1}}^{t_{n}} \frac{t^{-\tilde{\alpha}_n}}{\Gamma\left(1-\tilde{\alpha}_n\right)} dt
  \geq\frac{T^{-\tilde{\alpha}_n}}{\Gamma\left(1-\tilde{\alpha}_n\right)}.
 \end{aligned}
\end{equation*}
If $T\leq1$, we have $T^{-\tilde{\alpha}_n}\geq T^{-\alpha_*}$. Conversely, if $T>1$, we can get $T^{-\tilde{\alpha}_n}\geq T^{-\alpha^*}$.
Thus, we can see that
\begin{equation*}
 \begin{aligned}
    a_{n}^{(n)}
  \geq
  \frac{\min\{T^{-{\alpha}_*},T^{-{\alpha}^*}\}}{\Gamma\left(1-{\alpha}^*\right)}.
 \end{aligned}
\end{equation*}
The proof is completed.\hfill
$\blacksquare$

\begin{lemma} \label{lemma:sumbk}
For $0<\tilde{\alpha}_k \leq \alpha^*<1$, we have
\begin{equation*}
  \sum_{k=1}^{n}b_{k}^{(k)} \leq C_{12},
\end{equation*}
where $C_{12}$ is a positive constant.
\end{lemma}
{\bf Proof.}
We will prove the lemma in two steps. In the first step, we estimate $b_{k}^{(k)}$ individually for $k = 1, 2, \cdots, n$. In the second step,
we consider the summation of $b_{k}^{(k)}$ for $k$ from $1$ to $n$.\\
\textit{Step} 1.
According to the definition of $b_{k}^{(k)}$, we have
\begin{equation*}
 \begin{aligned}
  b_{1}^{(1)}
  = 1+ \tau a_{1}^{(1)}
  = 1+ \frac{\tau^{1- \tilde{\alpha}_1}}{\Gamma\left(3-\tilde{\alpha}_1\right)},
 \end{aligned}
\end{equation*}
and
\begin{equation*}
 \begin{aligned}
  b_{k}^{(k)}
  = \frac{1}{\tau \Gamma\left(3-\tilde{\alpha}_k\right)}
    \left[\left(t_{k}^{2-\tilde{\alpha}_k}-t_{k-1}^{2-\tilde{\alpha}_k}\right)-\left(t_{k-1}^{2-\tilde{\alpha}_k}-t_{k-2}^{2-\tilde{\alpha}_k}\right)\right]
  \leq
   \frac{t_{k-1}^{1-\tilde{\alpha}_k}-t_{k-2}^{1-\tilde{\alpha}_{k}}}{\Gamma\left(2-\tilde{\alpha}_k\right)}, \quad 2 \leq k \leq n.
 \end{aligned}
\end{equation*}
Specially, when $k=2$, we can directly get the estimate  $b_{2}^{(2)} \leq \frac{t_{1}^{1-\tilde{\alpha}_2}}{\Gamma\left(2-\tilde{\alpha}_2\right)}$.
When $k \geq 3$, we have
\begin{equation*}
 \begin{aligned}
  b_{k}^{(k)}
  \leq
  \frac{\tau}
        {\Gamma\left(1-\tilde{\alpha}_k\right)}
            \cdot t_{k-2}^{-\tilde{\alpha}_k}, \quad 3\leq k \leq n.
 \end{aligned}
\end{equation*}
\textit{Step} 2.
Based on the fact $ \Gamma\left(1-\tilde{\alpha}_k\right) >1$, we summate $b_{k}^{(k)}$ for $k$ from $1$ to $n$,
\begin{equation} \label{sum1}
 \begin{aligned}
   \sum_{k=1}^{n} b_{k}^{(k)}
   &\leq 1+ \frac{\tau^{1- \tilde{\alpha}_1}}{\Gamma\left(3-\tilde{\alpha}_1\right)}
              +\frac{t_{1}^{1-\tilde{\alpha}_2}}
                    {\Gamma\left(2-\tilde{\alpha}_2\right)}
             +\tau \sum_{k=3}^{n} \frac{t_{k-2}^{-\tilde{\alpha}_{k}}}
             {\Gamma\left(1-\tilde{\alpha}_k\right)}\\
 &\leq 1+\frac{\tau^{1-\tilde{\alpha}_1}}{\Gamma\left(3-\tilde{\alpha}_1\right)}
    +\frac{\tau^{1-\tilde{\alpha}_2}}
                    {\Gamma\left(2-\tilde{\alpha}_2\right)}
    +\tau \sum_{k=1}^{n-2} t_{k}^{-\tilde{\alpha}_{k+2}}.
  \end{aligned}
\end{equation}

Next, we discuss the summation depending on the value of $t_{n-2}$.

(I) If $t_{n-2}\leq1$, then $t_{k}^{-\tilde{\alpha}_{k+2}} \leq t_{k}^{-\alpha^{*}}$.
The last term of (\ref{sum1}) can be estimated as
\begin{equation} \label{equationBk1}
 \begin{aligned}
  \tau \sum_{k=1}^{n-2} t_{k}^{-\tilde{\alpha}_{k+2}}
\leq
  \tau \sum_{k=1}^{n-2} t_{k}^{-\alpha^{*}}
=\tau^{1-\alpha^{*}} \sum_{k=1}^{n-2} k^{-\alpha^{*}}
\leq
  \tau^{1-\alpha^{*}} \int_{0}^{n-2} s^{-\alpha^{*}}ds
=\frac{t_{n-2}^{1-\alpha^{*}}}{1-\alpha^{*}}
\leq C_{12}.
 \end{aligned}
\end{equation}
The other terms in (\ref{sum1}) are also bounded.

(II) If  $t_{n-2} > 1$, then there exists an integer $k^{*}$  such that $t_{k} \leq 1$ for
$1 \leq k \leq  k^{*}$, and $t_{k} > 1$ for $k^{*} + 1 \leq k \leq  n$.
The summation of $b_{k}^{(k)}$ for $k$ from $1$ to $k^*$ is similar to (\ref{equationBk1}), that is,
$\tau \sum_{k=1}^{k^*} t_{k}^{-\tilde{\alpha}_{k+2}} \leq \frac{t_{k^*}^{1-\alpha^{*}}}{1-\alpha^{*}} $.
Then we have
\begin{equation*}
 \begin{aligned}
   \sum_{k=1}^{n} b_{k}^{(k)}
  &\leq
    1+\frac{\tau^{1-\tilde{\alpha}_1}}
           {\Gamma\left(3-\tilde{\alpha}_1\right)}
    +\frac{\tau^{1-\tilde{\alpha}_2}}
          {\Gamma\left(2-\tilde{\alpha}_2\right)}
    +\tau \left(\sum_{k=1}^{k^{*}} t_{k}^{-\tilde{\alpha}_{k+2}}
    +\sum_{k=k^{*}+1}^{n-2} t_{k}^{-\tilde{\alpha}_{k+2}}\right)\\
  &\leq
     1+\frac{\tau^{1-\tilde{\alpha}_1}}
            {\Gamma\left(3-\tilde{\alpha}_1\right)}
     +\frac{\tau^{1-\tilde{\alpha}_2}}
           {\Gamma\left(2-\tilde{\alpha}_2\right)}
     +\frac{t_{k^*}^{1-\alpha^{*}}}{1-\alpha^{*}}
     +t_{n-k^{*}-2}
  \leq C_{12}.
 \end{aligned}
\end{equation*}
The proof of Lemma \ref{lemma:sumbk} is completed.\hfill
$\blacksquare$

In addition, the following several lemmas on the properties of the operators defined above are necessary in the stability analysis.

\begin{lemma}\label{etax}
For the operator $\boldsymbol{\theta}_{x}$ defined in \eqref{theta x}, there exists an operator $\boldsymbol{\zeta}_{x}$ satisfying $\boldsymbol{\theta}_{x} = \boldsymbol{\zeta}_{x}^2$.
Similarly, there exists an operator $\boldsymbol{\zeta}_{y}$ satisfying $\boldsymbol{\theta}_{y} = \boldsymbol{\zeta}_{y}^2$.
\end{lemma}
{\bf Proof.}
We only prove the result for $\boldsymbol{\theta}_{x}$, and the result for $\boldsymbol{\theta}_{y}$ can be obtained similarly.
According to Lemma \ref{lemma:basisf}, the matrix representation of the operator $\boldsymbol{\theta}_{x}$ for one-dimensional case is
\begin{equation} \label{Matrix:e1}
  \boldsymbol{Q}=\f{1}{8}
    \begin{pmatrix}
           4     &    4   &        &        &    \tbf{0} \\
           1     &    6   &    1   &        &            \\
                 & \ddots & \ddots & \ddots &            \\
                 &        &    1   &    6   &     1      \\
        \tbf{0}  &        &        &    4   &     4      \\
    \end{pmatrix}_{(M_x+2)}.
\end{equation}
We let $\boldsymbol{S}=diag\{2,1,\cdots,1,2\}$ with sizes $M_x+2$,
and we define the matrix $\boldsymbol{A}$ as
\begin{equation*}
  \boldsymbol{A}:=\boldsymbol{S}^{-1}\boldsymbol{Q}\boldsymbol{S}=\frac{1}{8}
    \begin{pmatrix}
           4     &    2   &        &        &            &        &   \tbf{0}  \\
           2     &    6   &    1   &        &            &        &            \\
                 &    1   &    6   &    1   &            &        &            \\
                 &        & \ddots & \ddots &   \ddots   &        &            \\
                 &        &        &    1   &     6      &    1   &            \\
                 &        &        &        &     1      &    6   &      2     \\
       \tbf{0}   &        &        &        &            &    2   &      4     \\
    \end{pmatrix}_{(M_x+2)},
\end{equation*}
which is a symmetric and positive definite matrix. There exists a unique symmetric positive definite matrix $\boldsymbol {B}$
such that $\boldsymbol {A} = \boldsymbol {B}^{2}$. Thus, we have
$\boldsymbol{Q}=\boldsymbol{S}\boldsymbol{B}^{2}\boldsymbol{S}^{-1}=\left(\boldsymbol{S}\boldsymbol {B}\boldsymbol{S}^{-1}\right)^{2}$.
Accordingly, there exist an operator $\boldsymbol{\zeta}_{x}$ satisfying $\boldsymbol{\theta}_{x} = \boldsymbol{\zeta}_{x}^2$.
Similarly, there is an operator $\boldsymbol{\zeta}_{y}$ satisfying $\boldsymbol{\theta}_{y} = \boldsymbol{\zeta}_{y}^2$.\hfill
$\blacksquare$

\begin{lemma}\label{etaxnorm}
For any $v \in \mathring{\mathcal{M}}_{h}$, we have
\begin{equation*}
 \begin{aligned}
  &\frac{3}{16}\|v\|^{2}\leq  \left\|\boldsymbol\zeta_{x} v\right\|^{2}=\left(\boldsymbol\zeta_{x} v, \boldsymbol\zeta_{x} v\right)=\left(\boldsymbol\theta_{x} v, v\right) \leq \|v\|^{2},\\
  &\frac{3}{16}\|v\|^{2}\leq  \left\|\boldsymbol\zeta_{y} v\right\|^{2}=\left(\boldsymbol\zeta_{y} v, \boldsymbol\zeta_{y} v\right)=\left(\boldsymbol\theta_{y} v, v\right) \leq \|v\|^{2}.
 \end{aligned}
\end{equation*}
\end{lemma}
{\bf Proof.}
We only consider the first estimate due to the similarity of them.
Based on the definition of $\boldsymbol\theta_{x}$ in (\ref{theta x}), we can get
\begin{equation} \label{IP-1}
 \begin{aligned}
\left(\boldsymbol\theta_{x} v, v\right)
&=\Delta x \Delta y \sum_{i=0}^{M_{x}+1} \sum_{j=0}^{M_{y}+1}\left(\boldsymbol\theta_{x} v_{i j}\right)\left(v_{i j}\right) \\
&=\Delta x \Delta y \sum_{j=0}^{M_{y}+1}
  \bigg(\frac{1}{8} v_{0 j} v_{1 j}
       +\frac{1}{4} \sum_{i=1}^{M_{x}-1} v_{i j} v_{i+1, j} + \frac{3}{4} \sum_{i=1}^{M_{x}} v_{i j}^{2}
       + \frac{1}{8} v_{M_{x}-1, j} v_{M_{x}, j}\bigg).
\end{aligned}
\end{equation}
We first use the inequality $2ab \leq a^2+b^2$ in equality (\ref{IP-1}) to obtain
\begin{equation*}
 \begin{aligned}
  (\boldsymbol\theta_{x}v, v)
 &\leq
   \Delta x \Delta y \sum_{j=0}^{M_{y}+1}
   \bigg(\frac{1}{16} v_{0 j}^{2} + \frac{15}{16} v_{1 j}^{2}
   +\sum_{i=2}^{M_{x}-1} v_{i j}^{2}
       + \frac{15}{16} v_{M_{x},j}^{2}+ \frac{1}{16} v_{M_{x}+1,j}^{2} \bigg)\\
 &\leq
   \Delta x \Delta y \sum_{i=0}^{M_{x}+1} \sum_{j=0}^{M_{y}+1} v_{i j}^{2}=\|v\|^{2}.
 \end{aligned}
\end{equation*}
Then using the inequality $2ab \geq -a^2-b^2$ in equality (\ref{IP-1}), together with $\boldsymbol\theta_{x} v_{0 j}=\boldsymbol\theta_{x} v_{M_{x}+1, j}=0$
for $j=0,1,\cdots,M_y+1$, we can get
\begin{equation*}
 \begin{aligned}
  (\boldsymbol\theta_{x} v, v)
 &\geq
   \Delta x \Delta y \sum_{j=0}^{M_{y+1}}\bigg(\big(\boldsymbol\theta_{x} v_{0 j}\big)\big(v_{0 j}\big)
   -\frac{1}{16} v_{0 j}^{2}+\frac{9}{16} v_{1 j}^{2}
   +\frac{1}{2} \sum_{i=2}^{M_{x-1}} v_{i j}^{2}
   +\frac{9}{16} v_{M_{x}, j}^{2}
   -\frac{1}{16} v_{M_{x}+1, j}^{2}\\
   &\quad+\big(\boldsymbol\theta_{x} v_{M_{x}+1, j}\big)\big(v_{M_{x}+1,j}\big)\bigg)\\
 &\geq
   \Delta x \Delta y \sum_{j=0}^{M_{y}+1}\bigg(
   \frac{3}{16} v_{0 j}^{2} + \frac{5}{16} v_{1 j}^{2} + \frac{1}{2} \sum_{i=2}^{M_{x}-1} v_{i j}^{2}
   +\frac{5}{16} v_{M_{x},j}^{2} + \frac{3}{16} v_{M_{x}+1,j}^{2}\bigg)\\
 &\geq
   \frac{3}{16}\Delta x \Delta y \sum_{i=0}^{M_{x}+1} \sum_{j=0}^{M_{y}+1} v_{i j}^{2}
   =\frac{3}{16}\|v\|^{2}.
 \end{aligned}
\end{equation*}
The proof of the first estimate is completed, and the second result can be obtained similarly.\hfill
$\blacksquare$

\begin{lemma}[\cite{L21-1}]\label{lemma:left}
If $b_{1}^{(n)}>b_{2}^{(n)}>\cdots>b_{n}^{(n)}>0$, $n = 1,2, \cdots , N$, then for any quadratic spline solution $u_h\in \mathcal{V}^{0}$, the following estimate holds,
\begin{equation*}
   \sum_{k=1}^{n} b_{n-k+1}^{(n)}\left(\boldsymbol{\theta}_{x} \boldsymbol{\theta}_{y} c^{k}-\boldsymbol{\theta}_{x} \boldsymbol{\theta}_{y} c^{k-1}, \boldsymbol{\theta}_{x} \boldsymbol{\theta}_{y} c^{n}\right)
   \geq \frac{1}{2}\left[\sum_{k=1}^{n}b_{n-k+1}^{(n)}\left(\left\|\boldsymbol{\theta}_{x} \boldsymbol{\theta}_{y} c^{k}\right\|^{2}-\left\|\boldsymbol{\theta}_{x} \boldsymbol{\theta}_{y} c^{k-1}\right\|^{2}\right)\right],
 \end{equation*}
where $c^k=\big\{ c_{ij}^k,\ (i,j)\in \bar{\Lambda} \big\}$, for $k=1,2,\cdots,n $, are the DOFs in the expression \eqref{qsc} of $u_h^{k}$.
\end{lemma}

\begin{lemma}\label{lemma:right}
For the QSC-$L1^+$ scheme \eqref{QSC-L1+F}-\eqref{bound2F} for model \eqref{equation1}-\eqref{equation3}, we have the following estimate,
\begin{equation}\label{rightproof}
 \begin{aligned}
 \left(\big(\boldsymbol{\eta}_{x} \boldsymbol{\theta}_{y}+\boldsymbol{\eta}_{y} \boldsymbol{\theta}_{x}\big) c^{n-\frac{1}{2}}, \boldsymbol{\theta}_{x} \boldsymbol{\theta}_{y} c^{n}\right)
 \leq
 -\frac{1}{4}\left(\left|\boldsymbol{\zeta}_{x} \boldsymbol{\theta}_{y} c^{n}\right|_{1 x}^{2}+\left|\boldsymbol{\zeta}_{y}\boldsymbol{\theta}_{x} c^{n}\right|_{1 y}^{2}\right)
 +\frac{1}{4}\left(\left|\boldsymbol{\zeta}_{x}\boldsymbol{\theta}_{y} c^{n-1}\right|_{1 x}^{2}+\left|\boldsymbol{\zeta}_{y}\boldsymbol{\theta}_{x} c^{n-1}\right|_{1 y}^{2}\right).
\end{aligned}
\end{equation}
\end{lemma}
{\bf Proof.}
Recalling the notation $c^{n-\frac{1}{2}} = \frac{1}{2}\left( c^n+c^{n-1} \right)$, the left hand side of (\ref{rightproof}) can be separated as
\begin{equation} \label{rightEQ1}
 \begin{aligned}
  & \left(\big(\boldsymbol{\eta}_{x} \boldsymbol{\theta}_{y}+\boldsymbol{\eta}_{y} \boldsymbol{\theta}_{x}\big) c^{n-\frac{1}{2}}, \boldsymbol{\theta}_{x} \boldsymbol{\theta}_{y} c^{n}\right)\\
  &=\frac{1}{2}\left(\boldsymbol{\eta}_{x} \boldsymbol{\theta}_{y} c^{n-1}, \boldsymbol{\theta}_{x}\boldsymbol{\theta}_{y} c^{n}\right)
    +\frac{1}{2}\left(\boldsymbol{\eta}_{x} \boldsymbol{\theta}_{y} c^{n}, \boldsymbol{\theta}_{x}\boldsymbol{\theta}_{y} c^{n}\right)
    +\frac{1}{2}\left(\boldsymbol{\eta}_{y} \boldsymbol{\theta}_{x} c^{n-1}, \boldsymbol{\theta}_{y}\boldsymbol{\theta}_{x} c^{n}\right)
    +\frac{1}{2}\left(\boldsymbol{\eta}_{y} \boldsymbol{\theta}_{x} c^{n}, \boldsymbol{\theta}_{y}\boldsymbol{\theta}_{x} c^{n}\right)\\
  &:= \sum_{i}^{4} P_{i}.\\
\end{aligned}
\end{equation}
Since $P_{1}$ and $P_{2}$ are similar as $P_{3}$ and $P_{4}$, respectively, and we just investigate $P_{1}$ and $P_{2}$.
For simplicity, we turn to the one-dimensional case for the term $P_{1}$. Together with $\boldsymbol\theta_{x}c_{0}=\boldsymbol\theta_{x}c_{M_{x}+1}=0$,
we have
\begin{equation*}
 \begin{aligned}
  P_1
=&\frac{1}{2} \Delta x \sum_{i=0}^{M_{x}+1}\left(\boldsymbol\eta_{x} c_{i}^{n-1}\right)\left(\boldsymbol\theta_{x}c_{i}^{n}\right)\\
=&\frac{1}{2} \sum_{i=1}^{M_{x}}\left(\boldsymbol{\vartheta}_{x} c_{i+1}^{n-1}-\boldsymbol{\vartheta}_{x} c_{i}^{n-1}\right)\left(\boldsymbol\theta_{x}c_{i}^{n}\right)\\
=&-\frac{1}{2}\left(\boldsymbol{\vartheta}_{x} c_{1}^{n-1}\right)\left(\boldsymbol\theta_{x} c_{1}^{n}\right)
    -\frac{1}{2} \sum_{i=2}^{M_{x}}\left(\boldsymbol{\vartheta}_{x} c_{i}^{n-1}\right)\left(\boldsymbol\theta_{x} c_{i}^{n}-\boldsymbol\theta_{x} c_{i-1}^{n}\right)
    +\frac{1}{2}\left(\boldsymbol{\vartheta}_{x} c_{M_x +1}^{n-1}\right)\left(\boldsymbol\theta_{x} c_{M_x}^{n}\right)\\
=&-\frac{\Delta x}{2} \sum_{i=1}^{M_{x}+1}\left(\boldsymbol{\vartheta}_{x} c_{i}^{n-1}\right) \left(\boldsymbol{\vartheta}_{x} \boldsymbol\theta_{x} c_{i}^{n}\right)
  =-\frac{1}{2}\left(\boldsymbol{\vartheta}_{x} c^{n-1},\boldsymbol{\vartheta}_{x} \boldsymbol\theta_{x} c^{n}\right).
 \end{aligned}
\end{equation*}
Thus, we can obtain the following result for two-dimensional case with $c \Leftarrow \boldsymbol\theta_{y} c$ and Lemma \ref{etax},
\begin{equation*}
 \begin{aligned}
  P_1
  =-\frac{1}{2}\left(\boldsymbol{\vartheta}_{x} \boldsymbol\theta_{y} c^{n-1}, \boldsymbol{\vartheta}_{x} \boldsymbol\theta_{x} \boldsymbol\theta_{y} c^{n}\right)
  =-\frac{1}{2}\left(\boldsymbol{\vartheta}_{x} \boldsymbol\theta_{y} c^{n-1}, \boldsymbol{\vartheta}_{x} \boldsymbol{\zeta}_{x}^2 \boldsymbol\theta_{y} c^{n}\right)
  =-\frac{1}{2}\left(\boldsymbol{\vartheta}_{x}\boldsymbol{\zeta}_{x} \boldsymbol\theta_{y} c^{n-1}, \boldsymbol{\vartheta}_{x} \boldsymbol{\zeta}_{x} \boldsymbol\theta_{y} c^{n}\right).
 \end{aligned}
\end{equation*}
Using the similar routine for $P_{1}$, we can get
\begin{equation*}
 \begin{aligned}
   P_2
   =-\frac{1}{2}\left(\boldsymbol{\vartheta}_{x}\boldsymbol{\zeta}_{x} \boldsymbol\theta_{y} c^{n}, \boldsymbol{\vartheta}_{x} \boldsymbol{\zeta}_{x} \boldsymbol\theta_{y} c^{n}\right)
   =-\frac{1}{2}\left|\boldsymbol{\zeta}_{x} \boldsymbol{\theta}_{y} c^{n}\right|_{1 x}^{2}.
 \end{aligned}
\end{equation*}
Furthermore, with the equality $2ab =(a+b)^2-a^2 - b^2$, we will obtain
\begin{equation*}
 \begin{aligned}
  P_1+P_2
=&-\frac{1}{2}\left(\boldsymbol{\vartheta}_{x}\boldsymbol{\zeta}_{x} \boldsymbol\theta_{y} c^{n-1}, \boldsymbol{\vartheta}_{x} \boldsymbol{\zeta}_{x} \boldsymbol\theta_{y} c^{n}\right)
   -\frac{1}{2}\left|\boldsymbol{\zeta}_{x} \boldsymbol{\theta}_{y} c^{n}\right|_{1 x}^{2}\\
=&-\frac{1}{2}\left|\boldsymbol{\zeta}_{x} \boldsymbol{\theta}_{y} c^{n}\right|_{1 x}^{2}
    -\frac{1}{4}\left(\left\|\boldsymbol{{\vartheta}}_{x} \boldsymbol{\zeta}_{x} \boldsymbol{\theta}_{y} c^{n}
    +\boldsymbol{{\vartheta}}_{x}\boldsymbol{\zeta}_{x}\boldsymbol{\theta}_{y} c^{n-1}\right\|^2-\left|\boldsymbol{\zeta}_{x} \boldsymbol{\theta}_{y} c^{n}\right|_{1 x} ^{2}
    -\left|\boldsymbol{\zeta}_{x} \boldsymbol{\theta}_{y} c^{n-1}\right|_{1 x}^{2}\right) \\
\leq&
-\frac{1}{4}\left(\left|\boldsymbol{\zeta}_{x} \boldsymbol{\theta}_{y} c^{n}\right|_{1 x}^{2}-\left|\boldsymbol{\zeta}_{x}\boldsymbol{\theta}_{y} c^{n-1}\right|_{1 x}^{2}\right).
 \end{aligned}
\end{equation*}
The terms $P_3$ and $P_4$ in (\ref{rightEQ1}) have similar results. Therefore, we can get the following estimate
\begin{equation*}
 \begin{aligned}
  \sum_{i}^{4} P_{i}
 \leq
 -\frac{1}{4}\left(\left|\boldsymbol{\zeta}_{x} \boldsymbol{\theta}_{y} c^{n}\right|_{1 x}^{2}+\left|\boldsymbol{\zeta}_{y}\boldsymbol{\theta}_{x} c^{n}\right|_{1 y}^{2}\right)
 +\frac{1}{4}\left(\left|\boldsymbol{\zeta}_{x}\boldsymbol{\theta}_{y} c^{n-1}\right|_{1 x}^{2}+\left|\boldsymbol{\zeta}_{y}\boldsymbol{\theta}_{x} c^{n-1}\right|_{1 y}^{2}\right).
 \end{aligned}
\end{equation*}
The proof of Lemma \ref{lemma:right} is completed.\hfill
$\blacksquare$

At last, we need the discrete Gr\"{o}nwall inequality.
\begin{lemma}[\cite{Gron}]\label{lemma:gronwall}
Let $v$, $w$ $\in \mathfrak{T}$ be nonnegative temporal grid functions,  and $C_{13}$ is a nonnegative constant.
If $v^{n} \leq (1 + \tau C_{13})v^{n-1} + \tau w^{n-1} $ for $1 \leq  n \leq  N$, then
\begin{equation*}
 \begin{aligned}
  v^{n} \leq  e^{C_{13}n\tau }\left[ v^{0} + \tau \sum_{l=0}^{n-1} w^{l}\right],\quad n = 1,2, \cdots , N.
 \end{aligned}
\end{equation*}
\end{lemma}

With all the lemmas above, we next consider the stability of the QSC-$L1^+$ scheme.
\subsection{The stability of the QSC-$L1^+$ scheme}
\begin{theorem}\label{QSC-L1+-Pr}
Assume that $\alpha^{\prime} (t) \leq 0$,
and suppose that $c^n=\{c_{ij}^n,\ (i,j)\in\bar{\Lambda},0\leq n\leq N\}$ is the solution of the QSC-$L1^+$ scheme \eqref{new QSC L1+}.
Then we have
\begin{equation} \label{PRtheorem}
 \begin{aligned}
  &\left\|\boldsymbol{\theta}_{x} \boldsymbol{\theta}_{y} c^{n}\right\|^{2}
   +\frac{3\tau\kappa}{32}\left(\left|\boldsymbol{\theta}_{y}c^{n}\right|_{1 x}^{2}+\left| \boldsymbol{\theta}_{x} c^{n}\right|^{2}_{1y}\right)\\
  &\leq
   C_{14} \left[\left\|\boldsymbol{\theta}_{x} \boldsymbol{\theta}_{y} c^{0}\right\|^{2}
   +\frac{\tau\kappa}{2}\left(\left| \boldsymbol{\theta}_{y} c^{0}\right|_{1 x} ^{2}
   +\left| \boldsymbol{\theta}_{x} c^{0}\right|_{1 y}^{2}\right)\right]
   +C_{15}\tau\sum_{k=1}^{n}\left\|f^{k-\frac{1}{2}}\right\|^{2}.
 \end{aligned}
\end{equation}
\end{theorem}
{\bf Proof.}
Multiplying both sides of equation (\ref{new QSC L1+}) by $2\Delta x\Delta y\boldsymbol\theta_{x}\boldsymbol\theta_{y}c_{ij}^{n}$,
and summing up for $i$ from $0$ to $M_{x}+1$ and for $j$ from $0$ to $M_{y}+1$, we can get
\begin{equation} \label{Sum-1}
 \begin{aligned}
  \sum_{k=1}^{n} b_{n-k+1}^{(n)}\left(\boldsymbol{\theta}_{x} \boldsymbol{\theta}_{y}c^{k}-\boldsymbol{\theta}_{x} \boldsymbol{\theta}_{y} c^{k-1}, 2\boldsymbol{\theta}_{x} \boldsymbol{\theta}_{y} c^{n}\right)
   =\tau\kappa \left((\boldsymbol{\eta}_{x} \boldsymbol{\theta}_{y}+\boldsymbol{\eta}_{y} \boldsymbol{\theta}_{x})c^{n-\frac{1}{2}}, 2\boldsymbol{\theta}_{x}\boldsymbol{\theta}_{y} c^{n}\right)
   +\tau\left(f^{n-\frac{1}{2}}, 2\boldsymbol{\theta}_{x} \boldsymbol{\theta}_{y} c^{n}\right).
 \end{aligned}
\end{equation}
For the summation term on the left hand side of (\ref{Sum-1}), we have by Lemma \ref{lemma:left} that
\begin{equation*}
   \sum_{k=1}^{n} b_{n-k+1}^{(n)}\left(\boldsymbol{\theta}_{x} \boldsymbol{\theta}_{y} c^{k}-\boldsymbol{\theta}_{x} \boldsymbol{\theta}_{y} c^{k-1}, 2\boldsymbol{\theta}_{x} \boldsymbol{\theta}_{y} c^{n}\right)
   \geq \sum_{k=1}^{n}b_{n-k+1}^{(n)}\left(\left\|\boldsymbol{\theta}_{x} \boldsymbol{\theta}_{y} c^{k}\right\|^{2}-\left\|\boldsymbol{\theta}_{x} \boldsymbol{\theta}_{y} c^{k-1}\right\|^{2}\right).
\end{equation*}
For the first term on the right hand side of (\ref{Sum-1}), we have by Lemma \ref{etax} and Lemma \ref{lemma:right} that
\begin{equation*}
 \begin{aligned}
 &\tau\kappa\left(\big(\boldsymbol{\eta}_{x} \boldsymbol{\theta}_{y}+\boldsymbol{\eta}_{y} \boldsymbol{\theta}_{x}\big) c^{n-\frac{1}{2}}, 2\boldsymbol{\theta}_{x} \boldsymbol{\theta}_{y} c^{n}\right)\\
 &\leq
 -\frac{\tau \kappa}{2}\left(\left|\boldsymbol{\zeta}_{x} \boldsymbol{\theta}_{y} c^{n}\right|_{1x}^{2}+\left|\boldsymbol{\zeta}_y \boldsymbol{\theta}_{x} c^{n}\right|_{1y}^{2}\right)
             +\frac{\tau\kappa}{2}\left(\left|\boldsymbol{\zeta}_{x} \boldsymbol{\theta}_{y} c^{n-1}\right|_{1 x}^{2}
             +\left|\boldsymbol{\zeta}_y \boldsymbol{\theta}_{x} c^{n-1}\right|_{1 y}^{2}\right).
\end{aligned}
\end{equation*}
Then, we can obtain
\begin{equation*}
 \begin{aligned}
  &\sum_{k=1}^{n} b_{n-k+1}^{(n)}\left\|\boldsymbol{\theta}_{x} \boldsymbol{\theta}_{y} c^{k}\right\|^{2}
   +\frac{\tau\kappa}{2}\left(\left|\boldsymbol{\zeta}_{x} \boldsymbol{\theta}_{y} c^{n}\right|_{1 x}^{2}
   +\left|\boldsymbol{\zeta}_y \boldsymbol{\theta}_{x} c^{n}\right|_{1 y}^{2}\right)\\
 & \leq
  \sum_{k=1}^{n-1} b_{n-k}^{(n)}\left\|\boldsymbol{\theta}_{x} \boldsymbol{\theta}_{y} c^{k}\right\|^{2}
   +\frac{\tau\kappa}{2}\left(\left|\boldsymbol{\zeta}_{x} \boldsymbol{\theta}_{y} c^{n-1}\right|_{1 x}^{2}
   +\left|\boldsymbol{\zeta}_y \boldsymbol{\theta}_{x} c^{n-1}\right|_{1 y}^{2}\right)
  +b_{n}^{(n)}\left\|\boldsymbol{\theta}_{x} \boldsymbol{\theta}_{y} c^{0}\right\|^{2}+2\tau\left(f^{n-\frac{1}{2}},\boldsymbol{\theta}_{x} \boldsymbol{\theta}_{y} c^{n}\right).
 \end{aligned}
\end{equation*}
Combining with Lemma \ref{lemma:bn}, we obtain
\begin{equation} \label{Sta-1}
 \begin{aligned}
  &\sum_{k=1}^{n} b_{n-k+1}^{(n)}\left\|\boldsymbol{\theta}_{x} \boldsymbol{\theta}_{y} c^{k}\right\|^{2}
   +\frac{\tau\kappa}{2}\left(\left|\boldsymbol{\zeta}_x \boldsymbol{\theta}_{y} c^{n}\right|_{1 x}^{2}
   +\left|\boldsymbol{\zeta}_y \boldsymbol{\theta}_{x} c^{n}\right|_{1 y}^{2}\right)\\
  &\leq
(1+C_{5}\tau)\sum_{k=1}^{n-1} b_{n-k}^{(n-1)}\left\|\boldsymbol{\theta}_{x} \boldsymbol{\theta}_{y} c^{k}\right\|^{2}
   +\frac{\tau\kappa}{2}\left(\left|\boldsymbol{\zeta}_x \boldsymbol{\theta}_{y} c^{n-1}\right|_{1 x}^{2}
   +\left|\boldsymbol{\zeta}_{y} \boldsymbol{\theta}_{x} c^{n-1}\right|_{1 y}^{2}\right) \\
   &\quad +b_{n}^{(n)}\left\|\boldsymbol{\theta}_{x} \boldsymbol{\theta}_{y} c^{0}\right\|^{2}+2\tau\left(f^{n-\frac{1}{2}},\boldsymbol{\theta}_{x} \boldsymbol{\theta}_{y} c^{n}\right).
 \end{aligned}
\end{equation}
Denote
\begin{equation*}
 \begin{aligned}
  G^{0}
  =\frac{\tau\kappa}{2}\left(\left|\boldsymbol{\zeta}_{x} \boldsymbol{\theta}_{y}c^{0}\right|_{1 x}^{2}
   +\left|\boldsymbol{\zeta}_y \boldsymbol{\theta}_{x} c^{0}\right|_{1y}^{2}\right),
 \end{aligned}
\end{equation*}
and
\begin{equation*}
 \begin{aligned}
  G^{n}
  =\sum\limits_{k=1}^{n} b_{n-k+1}^{(n)}\left\|\boldsymbol{\theta}_{x} \boldsymbol{\theta}_{y} c^{k}\right\|^{2}
   +\frac{\tau\kappa}{2}\left(\left|\boldsymbol{\zeta}_{x} \boldsymbol{\theta}_{y}c^{n}\right|_{1 x}^{2}
   +\left|\boldsymbol{\zeta}_y \boldsymbol{\theta}_{x} c^{n}\right|_{1y}^{2}\right),\quad1\leq n \leq N,
 \end{aligned}
\end{equation*}
 inequality (\ref{Sta-1}) can be simplified as
\begin{equation*}
 \begin{aligned}
  G^{n} \leq\left(1+C_{5} \tau\right) G^{n-1}+b_{n}^{(n)}\left\|\boldsymbol{\theta}_{x} \boldsymbol{\theta}_{y} c^{0}\right\|^{2}+2 \tau\left(f^{n-\frac{1}{2}}, \boldsymbol{\theta}_{x} \boldsymbol{\theta}_{y} c^{n}\right).
 \end{aligned}
\end{equation*}
It is easy to proof the theorem when $n=1$. Then, for $n \geq 2$, applying Lemma \ref{lemma:gronwall} to deduce that
\begin{equation}\label{Gn1}
 \begin{aligned}
  G^{n} \leq e^{C_{5} n \tau}\left[G^{0}+\sum_{k=1}^{n} b_{k}^{(k)}\left\|\boldsymbol{\theta}_{x} \boldsymbol{\theta}_{y} c^{0}\right\|^{2}
   +2 \tau \sum_{k=1}^{n}\left(f^{k-\frac{1}{2}}, \boldsymbol{\theta}_{x} \boldsymbol{\theta}_{y}c^{k}\right)\right],\quad 2\leq n\leq N.
 \end{aligned}
\end{equation}

According to the definition of $\Big\{ b_{n-k+1}^{(n)},\ k = 1,2, \cdots , n \Big\}$, we have $b_{1}^{(n)}=1+\tau a_{1}^{(n)} > 1+C_{11}\tau$, and $b_{n-k+1}^{(n)}>C_{11}\tau$, for $k=1,2,\cdots,n-1$, by Lemma \ref{lemma:andown}.
Then $G^{n}$ has the following lower bound
\begin{equation} \label{Gn2}
 \begin{aligned}
   G^{n} \geq\left\|\boldsymbol{\theta}_{x} \boldsymbol{\theta}_{y} c^{n}\right\|^{2}+\tau C_{11} \sum_{k=1}^{n}\left\|\boldsymbol{\theta}_{x} \boldsymbol{\theta}_{y} c^{k}\right\|^{2}
    +\frac{\tau\kappa}{2}\left(\left|\boldsymbol{\zeta}_x \boldsymbol{\theta}_{y} c^{n}\right|_{1x}^{2}+\left|\boldsymbol{\zeta}_y \boldsymbol{\theta}_{x}c^{n}\right|_{1 y}^{2}\right).
 \end{aligned}
\end{equation}
We combine estimates (\ref{Gn1}) and (\ref{Gn2}) to conclude that for $1\leq n\leq N$
\begin{equation}\label{young1}
 \begin{aligned}
  &\left\|\boldsymbol{\theta}_{x} \boldsymbol{\theta}_{y} c^{n}\right\|^{2}+ C_{11}\tau \sum_{k= 1}^{n} \left\|\boldsymbol{\theta}_{x} \boldsymbol{\theta}_{y} c^{k}\right\|^{2}
   +\frac{\tau\kappa}{2}\left(\left|\boldsymbol{\zeta}_x \boldsymbol{\theta}_{y} c^n\right|_{1 x}^{2}+\left|\boldsymbol{\zeta}_y \boldsymbol{\theta}_{x} c^{n}\right|_{1y}^{2}\right)\\
 & \leq
e^{C_{5} T}\left[\frac{\tau\kappa}{2}\left(\left|\boldsymbol{\zeta}_x \boldsymbol{\theta}_{y} c^0\right|_{1 x}^{2}
            +\left|\boldsymbol{\zeta}_y \boldsymbol{\theta}_{x} c^{0}\right|_{y}^{2}\right)
  +\sum_{k= 1}^{n} b_{k}^{(k)}\left\|\boldsymbol{\theta}_{x} \boldsymbol{\theta}_{y} c^{0}\right\|^{2}\right]
  +2e^{C_{5}T}\tau\sum_{k = 1}^{n}\left(f^{k-\frac{1}{2}},\boldsymbol{\theta}_{x} \boldsymbol{\theta}_{y} c^{k}\right).
 \end{aligned}
\end{equation}
For the third term on the left hand side and the first term on the right hand side of (\ref{young1}), we have from Lemma \ref{etaxnorm} that
\begin{equation} \label{right norm 1}
 \begin{aligned}
  \frac{\tau\kappa}{2}\left(\left|\boldsymbol{\zeta}_x \boldsymbol{\theta}_{y} c^n\right|_{1 x}^{2}
  +\left|\boldsymbol{\zeta}_y \boldsymbol{\theta}_{x} c^{n}\right|_{1y}^{2}\right)
  \geq
  \frac{3\tau\kappa}{32}\left(\left|\boldsymbol{\theta}_{y}c^{n}\right|_{1 x}^{2}+\left| \boldsymbol{\theta}_{x} c^{n}\right|^{2}_{1y}\right),
 \end{aligned}
\end{equation}
and
\begin{equation}\label{right norm 2}
 \begin{aligned}
  \frac{\tau\kappa}{2}\left(\left|\boldsymbol{\zeta}_x \boldsymbol{\theta}_{y} c^0\right|_{1 x}^{2}
            +\left|\boldsymbol{\zeta}_y \boldsymbol{\theta}_{x} c^{0}\right|_{1 y}^{2}\right)
  \leq
  \frac{\tau\kappa}{2}\left(\left| \boldsymbol{\theta}_{y} c^{0}\right|_{1 x} ^{2}
   +\left| \boldsymbol{\theta}_{x} c^{0}\right|_{1 y}^{2}\right).
 \end{aligned}
\end{equation}
For the last term in (\ref{young1}), recalling the inequality $ab\leq\ \varepsilon a^2 +(1/4\varepsilon) b^2 $, we can obtain
\begin{equation} \label{right young}
 \begin{aligned}
  2 e^{C_{5} T} \tau \sum_{k=1}^{n}\left(f^{k-\frac{1}{2}}, \boldsymbol{\theta}_{x} \boldsymbol{\theta}_{y} c^{k}\right)
   \leq  C_{11} \tau \sum_{k=1}^{n}\left\|\boldsymbol{\theta}_{x} \boldsymbol{\theta}_{y} c^{k}\right\|^{2}
   +\frac{\tau e^{2 C_{5} T}}{C_{11}} \sum_{k=1}^{n}\left\|f^{k-\frac{1}{2}}\right\|^{2}.
 \end{aligned}
\end{equation}
Then we can deduce from (\ref{young1}) -- (\ref{right young}) and Lemma \ref{lemma:sumbk} that
\begin{equation*}
 \begin{aligned}
  &\left\|\boldsymbol{\theta}_{x} \boldsymbol{\theta}_{y} c^{n}\right\|^{2}
   +\frac{3\tau\kappa}{32}\left(\left|\boldsymbol{\theta}_{y} c^{n}\right|_{1 x}^{2}+\left| \boldsymbol{\theta}_{x} c^{n}\right|^{2}_{1y}\right)\\
  &\leq
  C_{16} e^{C_{5} T}\left[\left\|\boldsymbol{\theta}_{x} \boldsymbol{\theta}_{y} c^{0}\right\|^{2}\right.
  \left.
  +\frac{\tau\kappa}{2}\left(\left| \boldsymbol{\theta}_{y} c^{0}\right|_{1 x} ^{2}
  +\left| \boldsymbol{\theta}_{x} c^{0}\right|_{1 y}^{2}\right)
  \right]
  +\frac{\tau e^{2 C_{5} T}}{C_{11}} \sum_{k= 1}^{n} \left\|f^{k-\frac{1}{2}}\right\|^{2}.
 \end{aligned}
\end{equation*}
The proof of Theorem \ref{QSC-L1+-Pr} is completed.\hfill
$\blacksquare$

\subsection{Convergence of the QSC-$L1^+$ scheme}
Based on the stability, we will investigate the convergence of the QSC-$L1^+$ scheme in this subsection.
For a function $w(x,y)\in C^4\big(\bar{\Omega}\big)$, we let $\mal I w(x,y)\in \mathcal{V}^{0}$
be the quadratic spline interpolation of $w(x,y)$, such that
\begin{equation}\label{QSC:Int}
  \mal I w\left(\xi _{i}^{x}, \xi _{j}^{y}\right)=w\left(\xi _{i}^{x}, \xi _{j}^{y}\right), \quad (i,j)\in \bar{\Lambda},
\end{equation}
where $\big(\xi _{i}^{x}, \xi _{j}^{y}\big)$ for $(i,j)\in \bar{\Lambda}$ are the collocation points described in Section \ref{sect:QSCL1+}.
Let $\|\cdot\|_c$ denote the maximum norm over all the
collocation points, \textit{i.e.},
\begin{equation*}
    \|w\|_c = \max_{(i,j)\in\bar{\Lambda}}|w\big(\xi _{i}^{x}, \xi _{j}^{y}\big)|.
\end{equation*}
Then it follows from the conclusions in \cite{KE,Iu-u} that the interpolation error $(\mal I w- w)$ satisfies
\begin{equation}\label{err:int}
    \begin{aligned}
      \| (\mal I w- w)_{xx} \|_c & \leq \f{\Delta x^{2}}{12} \| w^{(4)} \|_\infty + \mo(\Delta x^{3}),   \\
      \| (\mal I w- w)_{yy} \|_c & \leq \f{\Delta y^{2}}{12} \| w^{(4)} \|_\infty + \mo(\Delta y^{3}).   \\
    \end{aligned}
\end{equation}

We denote by
\begin{equation*}
 \begin{aligned}
  u^n = \left\{ u^n(\xi _{i}^{x}, \xi _{j}^{y}),\ (i,j)\in \bar{\Lambda} \right\}
 \end{aligned}
 \quad\text{and}\quad
 \begin{aligned}
  u_{h}^n = \left\{ u_{h}^n(\xi _{i}^{x}, \xi _{j}^{y}),\ (i,j)\in \bar{\Lambda} \right\}
 \end{aligned}
\end{equation*}
the true solution of the problem (\ref{equation1})-(\ref{equation3}) and the quadratic spline collocation solution of the
the QSC-$L1^+$ scheme (\ref{L1+2})-(\ref{ini2}), respectively, at the collocation points, where $u_{h}^n(x,y)$ has the expression (\ref{qsc}).
Then, we have the following conclusion.
\begin{theorem}\label{QSC-L1+-Con}
For $0<\alpha_{*} \leq\alpha(t)\leq  \alpha^{*} <1$, there exists a positive constant $C_{17}$, such that
\begin{equation*}
 \begin{aligned}
  \left\|u^{n}-u_{h}^{n}\right\|
  \leq
  C_{17}\left(\tau^{\min{\{3-\alpha^*-\alpha(0),2\}}}+\Delta x^{2}+\Delta y^{2}\right), \quad 1\leq n\leq N .
 \end{aligned}
\end{equation*}

\end{theorem}
{\bf Proof.}
According to equation (\ref{L1+3}), the interpolation function $\mal I u^{n}(x,y)$, for $n=1,2,\cdots,N$, satisfy the following equation

\begin{equation} \label{Con-1}
 \begin{aligned}
  \boldsymbol{\delta}_{t} \mal Iu^{n-\frac{1}{2}}(x, y)+ \bar{\boldsymbol{\delta}_{t}}^{\tilde{\alpha}_{n}}\mal Iu^{n-\frac{1}{2}}(x, y)
  =\kappa \left[\mal I u_{x x}^{n-\frac{1}{2}}(x, y)
                       +\mal I u_{y y}^{n-\frac{1}{2}}(x, y)\right]+f^{n-\frac{1}{2}}(x, y)+g^{n-\frac{1}{2}}(x, y),
 \end{aligned}
\end{equation}
where
\begin{equation} \label{Con-2}
 \begin{aligned}
  g^{n-\frac{1}{2}}(x, y)
  =&\boldsymbol{\delta}_{t} (\mal Iu-u)^{n-\frac{1}{2}}(x, y)+ \bar{\boldsymbol{\delta}_{t}}^{\tilde{\alpha}_{n}}(\mal Iu-u)^{n-\frac{1}{2}}(x, y)\\
   &- \kappa \left[\left(\mal I u-u\right)_{x x}^{n-\frac{1}{2}}(x, y)
                       +\left(\mal I u-u\right)_{y y}^{n-\frac{1}{2}}(x, y)\right]
                       +R^n.
 \end{aligned}
\end{equation}
We take the the collocation point $(\xi_{i}^{x}, \xi_{j}^{y})$ for $(i,j)\in \bar{\Lambda}$ into (\ref{Con-1}) and (\ref{Con-2}), and they can be rewritten as
\begin{equation}\label{equationC1}
 \begin{aligned}
  &\sum_{k=1}^{n} b_{n-k+1}^{(n)}\left[\mal I u^{k}\left(\xi _{i}^{x}, \xi _{j}^{y}\right)-\mal I u^{k-1}\left(\xi _{i}^{x}, \xi _{j}^{y}\right)\right]\\
  &=\tau\kappa
  \left[\mal I u_{x x}^{n-\frac{1}{2}}\left(\xi _{i}^{x}, \xi _{j}^{y}\right)
   +\mal I u_{y y}^{n-\frac{1}{2}}\left(\xi _{i}^{x}, \xi _{j}^{y}\right)\right]
  +\tau f^{n-\frac{1}{2}}\left(\xi _{i}^{x}, \xi _{j}^{y}\right)
  +\tau g^{n-\frac{1}{2}}\left(\xi _{i}^{x}, \xi _{j}^{y}\right),
 \end{aligned}
\end{equation}
where
\begin{equation*}
 \begin{aligned}
  g^{n-\frac{1}{2}}\left(\xi _{i}^{x}, \xi _{j}^{y}\right)
  =- \kappa \left[\left(\mal I u-u\right)_{x x}^{n-\frac{1}{2}}\left(\xi _{i}^{x}, \xi _{j}^{y}\right)
                       +\left(\mal I u-u\right)_{y y}^{n-\frac{1}{2}}\left(\xi _{i}^{x}, \xi _{j}^{y}\right)\right]
                       +R^n, \quad (i,j)\in \bar{\Lambda} ,
 \end{aligned}
\end{equation*}
$R^n$ is defined in (\ref{LTE}), and can be bounded by $\mo\left(\tau^{2}t_{n}^{-\tilde{\alpha}_n-\alpha(0)}+\tau\left(t_{n}^{1-\alpha(0)}-t_{n-1}^{1-\alpha(0)}\right)+\tau^{2}\right)$.

Since $\mal I u^n(x,y)\in \mathcal{V}^{0}$, it is reasonable to suppose that $\mal I u^n(x,y)$ can be expressed in the form
\begin{equation*}
 \begin{aligned}
  \mal I u^{n}(x,y)=\sum_{i=0}^{M_{x}+1} \sum_{j=0}^{M_{y}+1} d_{ij}^{n} \phi_{i}(x) \phi_{j}(y),
 \end{aligned}
\end{equation*}
where $d_{ij}^{n}$ are DOFs corresponding to $\mal I u^{n}(x,y)$.
Then equation (\ref{equationC1}) can be  rewritten as
\begin{equation}\label{equationC2}
 \begin{aligned}
  \sum_{k=1}^{n} b_{n-k+1}^{(n)}\left(\boldsymbol\theta_{x} \boldsymbol\theta_{y} d_{ij}^{k}
    -\boldsymbol\theta_{x} \boldsymbol\theta_{y} d_{ij}^{k-1}\right)
  =\tau \kappa \left(\boldsymbol\eta_{x} \boldsymbol\theta_{y}+\boldsymbol\eta_{y} \boldsymbol\theta_{x}\right) d_{ij}^{n-\frac{1}{2}}
   +\tau f_{ij}^{n-\frac{1}{2}}
   +\tau g_{ij}^{n-\frac{1}{2}}, \quad (i,j)\in \bar{\Lambda},
 \end{aligned}
\end{equation}
where $g_{ij}^{n-\frac{1}{2}} = g^{n-\frac{1}{2}}\left(\xi _{i}^{x}, \xi _{j}^{y}\right)$.
Denote $e^n = d^n - c^n$, we substitute (\ref{new QSC L1+}) from (\ref{equationC2}) to obtain
\begin{equation*}
 \begin{aligned}
  \sum_{k=1}^{n} b_{n-k+1}^{(n)}\left(\boldsymbol\theta_{x} \boldsymbol\theta_{y} e_{ij}^{k}
    -\boldsymbol\theta_{x} \boldsymbol\theta_{y} e_{ij}^{k-1}\right)
  =\tau \kappa \left(\boldsymbol\eta_{x} \boldsymbol\theta_{y}+\boldsymbol\eta_{y} \boldsymbol\theta_{x}\right) e_{ij}^{n-\frac{1}{2}}
   +\tau g_{ij}^{n-\frac{1}{2}},\quad (i,j)\in \bar{\Lambda},\ 1\leq n\leq N.
 \end{aligned}
\end{equation*}
Applying Theorem
\ref{QSC-L1+-Pr}, together with $e^0=0$, we can get
\begin{equation} \label{LTE-1}
 \begin{aligned}
  &\left\|\boldsymbol{\theta}_{x} \boldsymbol{\theta}_{y} e^{n}\right\|^{2}
   +\frac{3\tau\kappa}{32}\left(\left|\boldsymbol{\theta}_{y}e^{n}\right|_{1 x}^{2}+\left| \boldsymbol{\theta}_{x} e^{n}\right|^{2}_{1y}\right)\\
&\leq
  C_{14}
   \left[\left\|\boldsymbol{\theta}_{x} \boldsymbol{\theta}_{y} e^{0}\right\|^{2}
   +\frac{\tau\kappa}{2}\left(\left| \boldsymbol{\theta}_{y} e^{0}\right|_{1 x} ^{2}
   +\left| \boldsymbol{\theta}_{x} e^{0}\right|_{1 y}^{2}\right)\right]
   +C_{14}\tau \sum_{k=1}^{n} \left\|g^{k-\frac{1}{2}}\right\|^{2}
\leq
   C_{18}\tau \sum_{k=1}^{n}\left\|g^{k-\frac{1}{2}}\right\|_{c}^2,
 \end{aligned}
\end{equation}
where $C_{18} = (x_R-x_L)(y_R-y_L)C_{14}$.
Based on (\ref{LTE}) and (\ref{err:int}), we see
\begin{equation} \label{LTE-2}
 \begin{aligned}
  C_{18}\tau \sum_{k=1}^{n}\left\|g^{k-\frac{1}{2}}\right\|_{c}^2
  \leq
  &C_{18}\tau \sum_{k=1}^{n}\left[\kappa \left(\Delta x^{2}+\Delta y^{2}\right)+\mo\left(\tau^{2}t_{k}^{-\tilde{\alpha}_k-\alpha(0)}+\tau\left(t_{k}^{1-\alpha(0)}-t_{k-1}^{1-\alpha(0)}\right)+\tau^{2}\right)\right]^2\\
  \leq
  &C_{19}\left[\kappa T \left(\Delta x^{2}+\Delta y^{2}\right) + \tau\sum_{k=1}^{n} \left(\tau^{2}t_{k}^{-\tilde{\alpha}_k-\alpha(0)}+\tau\left(t_{k}^{1-\alpha(0)}-t_{k-1}^{1-\alpha(0)}\right)+\tau^{2}\right) \right]^2 \\
  \leq
  &C_{19}\left[\kappa T \left(\Delta x^{2}+\Delta y^{2}\right) +\tau^{2}t_{n}^{1-\alpha(0)}+ \tau\sum_{k=1}^{n} \left(\tau^{2}t_{k}^{-\tilde{\alpha}_k-\alpha(0)}+\tau^{2}\right) \right]^2.
 \end{aligned}
\end{equation}
Next, we give further discussions according to the value of $t_{n}$.

(I) If $t_{n} \leq 1$, we have $t_{k}^{-\tilde{\alpha}_k-\alpha(0)} \leq t_{k}^{-\alpha^{*}-\alpha(0)}$ for $ k\leq n$. Therefore, we have
\begin{equation} \label{LTE-3}
 \begin{aligned}
  \tau\sum_{k=1}^{n} \left(\tau^{2}t_k^{-\tilde{\alpha}_k-\alpha(0)}\right)
  &\leq \tau^{2}\left(\tau \sum_{k = 1}^{n} t_{k}^{-\alpha^{*}-\alpha(0)}\right)
  \leq \tau^{3-\alpha^{*}-\alpha(0)} + \tau^{2} \int_{t_{1}}^{t_{n}} t^{-\alpha^{*}-\alpha(0)} dt \\
  &= \frac{t_{n}^{1-\alpha^{*}-\alpha(0)}}{1-\alpha^{*}-\alpha(0)} \tau^{2} - \frac{\alpha^{*}+\alpha(0)}{1-\alpha^{*}-\alpha(0)} \tau^{3-\alpha^{*}-\alpha(0)}.
 \end{aligned}
\end{equation}

(II) If $t_{n} > 1$, then there exists an integer $k^{*}$, such that $t_{k} \leq 1$ for
$ k \leq  k^{*}$ and $t_{k} > 1$ for $k^{*} + 1 \leq k \leq  n$. With $t_{k}^{-\tilde{\alpha}_k-\alpha(0)} \leq t_{k}^{-\alpha_{*}-\alpha(0)}$ for $k^{*} + 1 \leq k \leq  n$, we can obtain
\begin{equation}\label{LTE-4}
 \begin{aligned}
  \tau \sum_{k=1}^{n}\left(\tau^{2} t_{k}^{-\alpha_{k}}\right)
& \leq \tau^{2}\left(\tau \sum_{k=1}^{k^{*}} t_{k}^{-\alpha^{*}-\alpha(0)}\right)+\tau^{2}\left(\sum_{k=k^{*}+1}^{n} t_{k}^{-\alpha_{*}-\alpha(0)}\right) \\
&=\frac{t_{k^{*}}^{1-\alpha^{*}-\alpha(0)}}{1-\alpha^{*}-\alpha(0)} \tau^{2}
- \frac{\alpha^{*}+\alpha(0)}{1-\alpha^{*}-\alpha(0)} \tau^{3-\alpha^{*}-\alpha(0)}
+\frac{t_{n}^{1-\alpha_*-\alpha(0)}-t_{k^{*}}^{1-\alpha_{*}-\alpha(0)}}{1-\alpha_{*}-\alpha(0)}\tau^{2}.
 \end{aligned}
\end{equation}
Thus, we can  substitute (\ref{LTE-3}) and (\ref{LTE-4}) into (\ref{LTE-2}) to obtain
\begin{equation} \label{LTE-5}
 \begin{aligned}
  C_{18}\tau \sum_{k=1}^{n}\left\|g^{k-\frac{1}{2}}\right\|_{c}^2
  \leq
  &C_{20}\left( \tau^{\min{\{3-\alpha^*-\alpha(0),2\}}}  + \Delta x^{2}+\Delta y^{2}\right)^2.
 \end{aligned}
\end{equation}
Based on (\ref{LTE-1}) and (\ref{LTE-5}), we can have
\begin{equation} \label{equationC3}
 \begin{aligned}
  \left\|\boldsymbol{\theta}_{x} \boldsymbol{\theta}_{y} e^{n}\right\|^{2}
   +\frac{3\tau\kappa}{32}\left(\left|\boldsymbol{\theta}_{y}e^{n}\right|_{1 x}^{2}+\left| \boldsymbol{\theta}_{x} e^{n}\right|^{2}_{1y}\right)
  \leq
  C_{20}\left( \tau^{\min{\{3-\alpha^*-\alpha(0),2\}}}  + \Delta x^{2}+\Delta y^{2} \right)^2.
 \end{aligned}
\end{equation}
Since $\big(\mal I u - u \big)^n\left(\xi _{i}^{x}, \xi _{j}^{y}\right)= 0$ for $(i,j)\in \bar{\Lambda}$, we can get
\begin{equation*}
    \left\|u^{n}-u_{h}^{n}\right\|
    =\| \mal I u^n - u_{h}^n \|
    = \| \boldsymbol{\theta}_{x} \boldsymbol{\theta}_{y} e^n \|.
\end{equation*}
This together with the estimation (\ref{equationC3}) complete the proof.
\hfill
$\blacksquare$

It can be seen that, if $\alpha(t)$ satisfies the condition $\alpha^{*}+ \alpha(0) < 1$, we can get $3-\alpha^{*}- \alpha(0) > 2$. Thus, the following 
corollary can be obtained directly from Theorem \ref{QSC-L1+-Con}.
\begin{corollary} \label{corollary-1}
   If the fractional order $\alpha(t)$ satisfies $0<\alpha_{*} \leq\alpha(t)\leq  \alpha^{*} <1$ and $ \alpha^{*}+ \alpha(0) < 1 $, there exists a positive constant $C_{21}$, such that
   \begin{equation*}
   	\begin{aligned}
   		\left\|u^{n}-u_{h}^{n}\right\|
   		\leq
   		C_{21}\left(\tau^{2}+\Delta x^{2}+\Delta y^{2}\right), \quad 1\leq n\leq N .
   	\end{aligned}
   \end{equation*}
\end{corollary}

\begin{remark}
If the variable fractional order $\alpha(t)$ satisfies the condition $\alpha^{*}+ \alpha(0) = 1$, 
the estimations (\ref{LTE-3})--(\ref{LTE-4}) in the proof of Theorem \ref{QSC-L1+-Con} need to be modified slightly, and the resulting convergence order is $\mo({\tau^{2} |\ln{\tau} |} + \Delta x^{2}+ \Delta y^{2})$, which is consistent with the result in Ref. \cite{WangL1+}.
\end{remark}
 
\begin{remark}
We can see from the proof of Theorem \ref{QSC-L1+-Con} that, the truncation error $\mo\big(\tau^{3-\alpha^{*}-\alpha(0)}\big)$ is from the error estimation near the initial time. We will investigate the behavior of the numerical solution near the initial time point and the final time point, respectively, in Example 6.1 in Section \ref{sect:numerical}. It will be seen that, the 
convergence order near the initial time behaves indeed as we have estimated, and the convergence order at the final time point behaves well.

\end{remark}

\begin{remark} \label{remark-2}
	If the solution of model (\ref{equation1})--(\ref{equation3}) has better regularity, such as $u \in C^2[0,T]$, the QSC-$L1^+$ scheme can achieve 
	the second temporal convergence order, without the restriction $\alpha^{*}+ \alpha(0) < 1$, which is consistent with the example confirmation in Ref. \cite{L1+2}.
\end{remark}

Next we will take the QSC-$L1^+$ scheme into the ADI framework for model (\ref{equation1})-(\ref{equation3}), defined in the two-dimensional space domain.
\section{The ADI-QSC-L1$^+$ scheme}
\label{sect:ADIQSCL1+}
It is known that the computational cost for multi-dimensional FPDEs is usually expensive.
The ADI method is able to change the solution of the multi-dimensional problem to the solutions of a series of one-dimensional subproblems,
and the computational cost can be efficiently reduced.
In this section, we will investigate the QSC-$L1^+$ scheme in the ADI framework.

We reformulate the QSC-$L1^+$ scheme (\ref{QSC-L1+F}) as
\begin{equation}\label{ADI3-1}
 \begin{aligned}
    &\boldsymbol\theta_{x} \boldsymbol\theta_{y} c_{ij}^{n}-\gamma_{n}\Big(\boldsymbol\eta_{x}  \boldsymbol\theta_{y}
    +\boldsymbol\eta_{y} \boldsymbol\theta_{x}\Big) c_{ij}^{n}\\
 &=\boldsymbol\theta_{x} \boldsymbol\theta_{y} c_{ij}^{n-1}
    +\gamma_{n}\Big(\boldsymbol\eta_{x} \boldsymbol\theta_{y}+\boldsymbol\eta_{y} \boldsymbol\theta_{x}\Big)c_{ij}^{n-1}
    -\frac{2 \gamma_{n}}{\kappa}\sum_{k=1}^{n-1} a_{n-k+1}^{(n)}\Big(\boldsymbol\theta_{x} \boldsymbol\theta_{y} c_{ij}^{k}
     -\boldsymbol\theta_{x} \boldsymbol\theta_{y} c_{ij}^{k-1}\Big)
   +\frac{2 \gamma_{n}}{\kappa}f_{ij}^{n-\frac{1}{2}},
 \end{aligned}
\end{equation}
for $(i,j)\in \bar{\Lambda}$ and $1\leq n\leq N$, where $\gamma_{n}=\frac{\tau \kappa }{2(1+\tau a_{1}^{(n)})}=\mo(\tau)$.
In order to implement the ADI method, we need to add a small perturbation, which depends on time levels.

(I) For $n = 1$, we define
\begin{equation*}
 \begin{aligned}
  Q_{ij}^{1}
   =\gamma_{1}^2  \boldsymbol\eta_{x} \boldsymbol\eta_{y} \left( c_{ij}^{1}-c_{ij}^{0}\right)
   =\frac{\tau^{3} \kappa^{2}}{4\left(1+\tau a_{1}^{(1)}\right)^2}\boldsymbol\eta_{x} \boldsymbol\eta_{y} \boldsymbol\delta_{t}c_{ij}^{\frac{1}{2}},\quad (i,j)\in \bar{\Lambda}.
 \end{aligned}
\end{equation*}
Adding the term ${Q} _{ij}^{1}$ on the left side of (\ref{ADI3-1}), we can get
\begin{equation} \label{ADI5}
 \begin{aligned}
  \big(\boldsymbol\theta_{x}-\gamma_{1} \boldsymbol\eta_{x} \big)\big(\boldsymbol\theta_{y}-\gamma_{1} \boldsymbol\eta_{y} \big) \tilde{c}_{ij}^{1}
  =\big(\boldsymbol\theta_{x}+\gamma_{1} \boldsymbol\eta_{x} \big)\big(\boldsymbol\theta_{y}+\gamma_{1} \boldsymbol\eta_{y} \big) {c}_{ij}^{0}
   +\frac{2 \gamma_{1}}{\kappa}f_{ij}^{\frac{1}{2}},\quad (i,j)\in \bar{\Lambda}.
 \end{aligned}
\end{equation}
Therefore, the ADI-QSC-$L1^+$ scheme for $n=1$ is implemented as following two steps:\\
\textit{Step} 1.
For each
$0\leq j\leq M_y+1$, solve the following one-dimensional linear systems in the $x$ direction for $\tilde{c}_{ij}^{1,*}$
\begin{equation} \label{eqADI1}
 \begin{aligned}
  \big(\boldsymbol\theta_{x}-\gamma_{1} \boldsymbol\eta_{x} \big)\tilde{c}_{ij}^{1,*}
  =\big(\boldsymbol\theta_{x}+\gamma_{1} \boldsymbol\eta_{x} \big) \big(\boldsymbol\theta_{y}+\gamma_{1} \boldsymbol\eta_{y} \big) {c}_{ij}^{0}
   +\frac{2 \gamma_{1}}{\kappa}f_{ij}^{\frac{1}{2}}, \quad 0\leq i\leq M_x+1.
 \end{aligned}
\end{equation}
\textit{Step} 2.
For each $0\leq i\leq M_x+1$, solve the following one-dimensional linear systems in the $y$ direction for $\tilde{c}_{ij}^{1}$
\begin{equation} \label{eqADI2}
 \begin{aligned}
  \big(\boldsymbol\theta_{y}-\gamma_{1} \boldsymbol\eta_{y} \big) \tilde{c}_{ij}^{1}=\tilde{c}_{ij}^{1,*}, \quad 0\leq j\leq M_y+1 .
 \end{aligned}
\end{equation}

(II) For the case $n \geq 2$, if we add the similar perturbation
$Q_{ij}^{n}=\gamma_{n}^2  \boldsymbol\eta_{x} \boldsymbol\eta_{y}\big( c_{ij}^{n}-c_{ij}^{n-1}\big)=\mo\big(\tau^3\big)$ as $Q_{ij}^{1}$ in the case (I),
which is equivalent to add $\mo\big(\tau^2\big)$ term to the truncation error of the QSC-$L1^+$ scheme.
Although adding such a small perturbation will not change the convergence order, some fundamental numerical tests
show that the observation error increase obviously. Thus, in this part, we aim to seek for a special small perturbation which has higher order.

We consider the term
\begin{equation*}
 \begin{aligned}
  Q_{ij}^n
   =\gamma_{n}^2 \boldsymbol\eta_{x} \boldsymbol\eta_{y}
     \left(c_{ij}^{n}-2c_{ij}^{n-1}+c_{ij}^{n-2}\right)
   =\frac{\tau^{2} \kappa^{2}}{4\left(1+\tau a_{1}^{(n)}\right)^2}\boldsymbol\eta_{x} \boldsymbol\eta_{y} \left(c_{ij}^{n}-2c_{ij}^{n-1}+c_{ij}^{n-2}\right),
    \quad (i,j)\in \bar{\Lambda}.
 \end{aligned}
\end{equation*}
It can be verified that
\begin{equation} \label{Q1}
 \begin{aligned}
  &\gamma_{n}^{2} \frac{\partial^{4}}{\partial x^{2} \partial y^{2}}\left[u(x,y,t_{n})-2 u(x,y,t_{n-1})+u(x,y,t_{n-2})\right]\\
&=\frac{\tau^{2}\gamma_{n}^{2} \partial^{6}}{\partial x^{2} \partial y^{2} \partial t^{2}} u(x,y,t_{n-\frac{1}{2}})
  -\frac{ \tau^{3}\gamma_{n}^{2} }{2} \cdot   \frac{\partial^{7}}{\partial x^{2} \partial y^2 \partial t^{3}} u(x,y,\rho_{2}),\quad \rho_{2} \in (t_{n-1},t_{n})
 \end{aligned}
\end{equation}
Recalling the regularity assumption of the solution, we have $\big\|\frac{\partial^{3}}{\partial t^{3}} u(\rho_{2} )\big\|_{{\mathcal{X}}} \leq C_{22} t^{-\alpha(0)-1}$ for $t_2\leq t\leq t_N$.
We can see that both the two terms on the right hand side of (\ref{Q1}) is equivalent to $\mo\big(\tau^{4-\alpha(0)}\big)$ when $t$ is
near the initial time point.

Adding the term ${Q}_{ij}^n$ on the left side of (\ref{ADI3-1}), we can get for $n\geq 2$ that
\begin{equation} \label{ADI7}
 \begin{aligned}
  &\left(\boldsymbol\theta_{x}-\gamma_{n} \boldsymbol\eta_{x} \right)\left(\boldsymbol\theta_{y}-\gamma_{n} \boldsymbol\eta_{y} \right) \tilde{c}_{ij}^{n}\\
&=\boldsymbol\theta_{x} \boldsymbol\theta_{y} \tilde{c}_{ij}^{n-1}
    +\gamma_{n}\left(\boldsymbol\eta_{x} \boldsymbol\theta_{y}+\boldsymbol\eta_{y} \boldsymbol\theta_{x}\right) \tilde{c}_{ij}^{n-1}
    +2\gamma_{n}^2 \boldsymbol\eta_{x}  \boldsymbol\eta_{y}  \tilde{c}_{ij}^{n-1}
    -\gamma_{n}^2 \boldsymbol\eta_{x}  \boldsymbol\eta_{y}  \tilde{c}_{ij}^{n-2}\\
   &\quad -\frac{2 \gamma_{n}}{\kappa}\sum_{k=1}^{n-1} a_{n-k+1}^{(n)}\left(\boldsymbol\theta_{x}\boldsymbol\theta_{y} \tilde{c}_{ij}^{k}
     -\boldsymbol\theta_{x} \boldsymbol\theta_{y}\tilde{c}_{ij}^{k-1}\right)
   +\frac{2 \gamma_{n}}{\kappa}f_{ij}^{n-\frac{1}{2}}, \quad (i,j)\in \bar{\Lambda}.
 \end{aligned}
\end{equation}
Compared with (\ref{ADI3-1}), scheme (\ref{ADI7}) has two extra terms, which is of the order $\mo\big(\tau^{3-\alpha(0)}\big)$ near the initial time.
The implementation of the ADI-QSC-$L1^+$ scheme for $n\geq 2$ is similar to the case (I).
Some numerical experiments show that
the ADI-QSC-$L1^+$ scheme (\ref{ADI7}) preserves almost the same observation error as the QSC-$L1^+$ scheme,
but the CPU time is effectively reduced in the ADI framework.

Next, we present the stability and convergence of the ADI-QSC-$L1^+$ scheme (\ref{ADI5}) and (\ref{ADI7}).
\begin{theorem}\label{ADI-Sta}
  The ADI-QSC-$L1^+$ scheme \eqref{ADI5} and \eqref{ADI7} are unconditionally stable.
Moreover, we denote by
$u^n = \left\{ u^n(\xi _{i}^{x}, \xi _{j}^{y}),\ (i,j)\in \bar{\Lambda} \right\}$
the true solution of the problem \eqref{equation1}-\eqref{equation3} and
$u_{h}^n = \left\{ u_{h}^n(\xi _{i}^{x}, \xi _{j}^{y}),\ (i,j)\in \bar{\Lambda} \right\}$
the numerical solution
by the ADI-QSC-$L1^+$ scheme \eqref{ADI5} and \eqref{ADI7} at the collocation points. Then, there exist a constant $C_{23}$ such that
\begin{equation*}
 \begin{aligned}
  \left\|u^{n}-u_{h}^{n}\right\|
  \leq
  C_{23}\left( \tau^{\min{\{3-\alpha^*-\alpha(0),2\}}} +\Delta x^{2}+\Delta y^{2}\right), \quad n=1,2, \ldots, N.
 \end{aligned}
\end{equation*}
\end{theorem}
{\bf Proof.}
The stability of the ADI-QSC-$L1^+$ scheme can be proved
by a similar routine in the proof of Theorem \ref{QSC-L1+-Pr}.
In fact, the ADI-QSC-$L1^+$ scheme can be expressed in the equivalent form
\begin{equation} \label{eqADI6}
 \begin{aligned}
  b_{1}^{(1)}\Big(\boldsymbol{\theta}_{x} \boldsymbol{\theta}_{y} c_{ij}^{1}-\boldsymbol{\theta}_{x} \boldsymbol{\theta}_{y} c_{ij}^{0}\Big)
   +\frac{\tau^3 \kappa ^{2}}{4 b_{1}^{(1)}} \boldsymbol\eta_{x}  \boldsymbol\eta_{y} \boldsymbol\delta_{t} c_{ij}^{\frac{1}{2}}
  =\tau \kappa \Big(\boldsymbol{\eta}_{x} \boldsymbol{\theta}_{y}+\boldsymbol{\eta}_{y} \boldsymbol{\theta}_{x}\Big) c_{ij}^{\frac{1}{2}}
   +\tau f_{ij}^{\frac{1}{2}}
 \end{aligned}
\end{equation}
for $n=1$, and
\begin{equation} \label{eqADI5}
 \begin{aligned}
  \sum_{k=1}^{n} b_{n-k+1}^{(n)}\Big(\boldsymbol{\theta}_{x} \boldsymbol{\theta}_{y} c_{ij}^{k}-\boldsymbol{\theta}_{x} \boldsymbol{\theta}_{y} c_{ij}^{k-1}\Big)
   +\frac{\tau^{2} \kappa ^{2}}{4 b_{1}^{(n)}} \boldsymbol\eta_{x}  \boldsymbol\eta_{y} \Big(c_{ij}^{n}-2c_{ij}^{n-1}+c_{ij}^{n-2}\Big)
  =\tau \kappa \Big(\boldsymbol{\eta}_{x} \boldsymbol{\theta}_{y}+\boldsymbol{\eta}_{y} \boldsymbol{\theta}_{x}\Big) c_{ij}^{n-\frac{1}{2}}
   +\tau f_{ij}^{n-\frac{1}{2}}
 \end{aligned}
\end{equation}
for $n\geq 2$. Next, we give further discussion based on the different value of $n$.

(I) For $n=1$,
equation (\ref{eqADI6}) has the same form as equation (\ref{new QSC L1+}) with $n=1$, excepting the second term on the left hand
side of (\ref{eqADI6}). We use the technique in Lemma \ref{lemma:right} to deal with this term as
\begin{equation*}
 \begin{aligned}
   \Big(\boldsymbol\eta_{x} \boldsymbol\eta_{y} \boldsymbol\delta_{t} c^{\frac{1}{2}}, \boldsymbol\theta_{x} \boldsymbol\theta_{y} c^{1}\Big)
  =&\Big(\boldsymbol\delta_{t} \boldsymbol{\vartheta}_{x} \boldsymbol{\vartheta}_{y} c^{\frac{1}{2}},
           \boldsymbol\theta_{x} \boldsymbol\theta_{y} \boldsymbol{\vartheta}_{x} \boldsymbol{\vartheta}_{y} c^{1}\Big)\\
  =&\frac{1}{\tau}
         \Big(\boldsymbol\zeta_{x} \boldsymbol\zeta_{y} \boldsymbol{\vartheta}_{x} \boldsymbol{\vartheta}_{y} c^{1} ,
               \boldsymbol\zeta_{x} \boldsymbol\zeta_{y} \boldsymbol{\vartheta}_{x} \boldsymbol{\vartheta}_{y} c^{1}\Big)
    -\frac{1}{\tau}
         \Big(\boldsymbol\zeta_{x} \boldsymbol\zeta_{y} \boldsymbol{\vartheta}_{x} \boldsymbol{\vartheta}_{y} c^{0},
               \boldsymbol\zeta_{x} \boldsymbol\zeta_{y} \boldsymbol{\vartheta}_{x} \boldsymbol{\vartheta}_{y} c^{1}\Big),
 \end{aligned}
\end{equation*}
where $\boldsymbol\zeta_{x}^2=\boldsymbol\theta_{x}$ and $\boldsymbol\zeta_{y}^2=\boldsymbol\theta_{y}$.
Then using the relation $2ab = a^2 + b^2-(a-b)^2$, we can get the following estimate
\begin{equation*}
 \begin{aligned}
  &\Big(\boldsymbol\eta_{x} \boldsymbol\eta_{y} \boldsymbol\delta_{t} c^{ \frac{1}{2}}, \boldsymbol\theta_{x} \boldsymbol\theta_{y} c^{1}\Big) \\
  &=
\frac{1}{\tau}\left\|\boldsymbol\zeta_{x} \boldsymbol\zeta_{y} \boldsymbol{\vartheta}_{x} \boldsymbol{\vartheta}_{y} c^{1}\right\|^{2}
   +\frac{1}{2\tau}
    \left[\left\|\boldsymbol\zeta_{x} \boldsymbol\zeta_{y} \boldsymbol{\vartheta}_{x} \boldsymbol{\vartheta}_{y} c^{1}
                -\boldsymbol\zeta_{x} \boldsymbol\zeta_{y} \boldsymbol{\vartheta}_{x} \boldsymbol{\vartheta}_{y} c^{0}\right\|^{2}
    -\left\|\boldsymbol\zeta_{x} \boldsymbol\zeta_{y} \boldsymbol{\vartheta}_{x} \boldsymbol{\vartheta}_{y} c^{1}\right\|^{2}
          -\left\|\boldsymbol\zeta_{x} \boldsymbol\zeta_{y} \boldsymbol{\vartheta}_{x} \boldsymbol{\vartheta}_{y} c^{0}\right\|^{2}\right]\\
  &\geq
  \frac{1}{2\tau}\left\|\boldsymbol\zeta_{x} \boldsymbol\zeta_{y} \boldsymbol{\vartheta}_{x} \boldsymbol{\vartheta}_{y} c^{1}\right\|^{2}
   -\frac{1}{2\tau}\left\|\boldsymbol\zeta_{x} \boldsymbol\zeta_{y} \boldsymbol{\vartheta}_{x} \boldsymbol{\vartheta}_{y} c^{0}\right\|^{2}.
 \end{aligned}
\end{equation*}
Then, we use the similar technique in Theorem \ref{QSC-L1+-Pr} and Theorem \ref{QSC-L1+-Con} to get
\begin{equation*}
 \begin{aligned}
  \left\|u^{1}-u_h^{1}\right\|^2
  \leq
  C_{24}\left( \tau^{\min{\{3-\alpha^*-\alpha(0),2\}}} + \Delta x^{2}+\Delta y^{2} \right)^2.
 \end{aligned}
\end{equation*}

(II) For $n\geq 2$, we can rewrite (\ref{eqADI5}) as
\begin{equation} \label{eqADI7}
 \begin{aligned}
  \sum_{k=1}^{n} b_{n-k+1}^{(n)}\left(\boldsymbol{\theta}_{x} \boldsymbol{\theta}_{y} c_{ij}^{k}-\boldsymbol{\theta}_{x} \boldsymbol{\theta}_{y} c_{ij}^{k-1}\right)
  =\tau \kappa \left(\boldsymbol{\eta}_{x} \boldsymbol{\theta}_{y}+\boldsymbol{\eta}_{y} \boldsymbol{\theta}_{x}\right) c_{ij}^{n-\frac{1}{2}}
   +\tau s_{ij}^{n-\frac{1}{2}},
 \end{aligned}
\end{equation}
where
\begin{equation}\label{eqADI8}
 \begin{aligned}
  s_{ij}^{n-\frac{1}{2}}
  =f_{ij}^{n-\frac{1}{2}}
   -\frac{\tau \kappa ^{2}}{4 b_{1}^{(n)}} \boldsymbol\eta_{x}  \boldsymbol\eta_{y}\left(c_{ij}^{n}-2c_{ij}^{n-1}+c_{ij}^{n-2}\right).
 \end{aligned}
\end{equation}
Recalling that the second term on the right hand of (\ref{eqADI8}) is of the order $\mo\big(\tau^{3-\alpha(0)}\big)$,
we use the similar routine in Theorem \ref{QSC-L1+-Pr} and Theorem \ref{QSC-L1+-Con} to get
\begin{equation*}
 \begin{aligned}
  \left\|u^{n}-u_h^{n}\right\|^2
  \leq
  C_{25}\tau \sum_{k=1}^{n}\left\|s^{k-\frac{1}{2}}\right\|_{c}^2
  \leq
  C_{26}\left( \tau^{\min{\{3-\alpha^*-\alpha(0),2\}}} + \Delta x^{2}+\Delta y^{2} \right)^2.
 \end{aligned}
\end{equation*}
Combining the case (I) and (II), the proof of theorem is completed.
\hfill
$\blacksquare$

\section{Acceleration techniques}
\label{sect:speed}
Computational efficiency of the numerical schemes for FPDEs has always been concerned.
In this section, we consider two kinds of techniques to accelerate the implementation of numerical schemes.
One is the fast computation based on the ESA technique along the time direction, the other is the optimal QSC method from the view of space direction.
\subsection{Fast computation in time direction}

We can see that the time discretization of Caputo fractional differential operator involves the numerical solutions at all previous time levels, the computation is extremely
expensive for long-time simulations.
In order to reduce the computational cost, we employ ESA technique to accelerate the evaluation of the
$L1^+$ scheme for variable-order FPDEs. The main purpose is to approximate the singular kernel $t^{-\beta}$ of the Caputo fractional differential operator
on the interval $[\tau, T]$ efficiently.
\begin{lemma}[\cite{Fast-ESA-1}]\label{Fast2}
At any time instant $t_n$, for $\tilde{\alpha}_{n} \in[\alpha_*, \alpha^*] \subset(0,1)$,
$t \in\left[t_{n-1},t_{n}\right]$, $s \in\left[0, t_{n-2}\right]$ and the expected accuracy $0<\epsilon \leq 1 / e$,
if we choose constants $h$, $\overline{N}$ and $\underline{N}$ as
\begin{equation*}
  \begin{aligned}
   h             &=\frac{2 \pi}{\log 3+\alpha^* \log (\cos 1)^{-1}+\log \epsilon^{-1}} ,\quad
   \underline{N} =\left \lceil \frac{1}{h} \frac{1}{\alpha_*}(\log \epsilon+\log \Gamma(1+\alpha^*))\right \rceil , \\
   \overline{N}  &=\left \lfloor\frac{1}{h}\left(\log \frac{T}{\Delta t}+\log \log \epsilon^{-1}+\log \alpha_*+2^{-1}\right) \right \rfloor,
  \end{aligned}
\end{equation*}
then the quantity $\left(\frac{t-s}{T}\right)^{-\tilde{\alpha}_{n }}$ can be approximated by
\begin{equation*}
 \begin{aligned}
  \left|\left(\frac{t-s}{T}\right)^{-\tilde{\alpha}_{n }}-\sum_{r=\underline{N}+1}^{\overline{N}} \varpi^{(n,r)} e^{\frac{-\lambda^{(r)}\left({t-s}\right)}{T}}\right|
  \leq \left(\frac{t-s}{T}\right)^{-\tilde{\alpha}_{n }} \epsilon,
 \end{aligned}
\end{equation*}
where the quadrature exponents and weights are given by
\begin{equation*}
 \begin{aligned}
  \lambda^{(r)}=e^{r h}, \quad \varpi^{(n,r)}=\frac{h e^{\tilde{\alpha}_{n } r h}}{\Gamma\left(\tilde{\alpha}_{n }\right)}.
 \end{aligned}
\end{equation*}
\end{lemma}

 Now for any $v \in \mathfrak{T}$, the non-local term $\bar{\delta}_{t}^{\tilde{\alpha}_{n}} v(t_{n-\frac{1}{2}})$ with $n\geq3$ defined in (\ref{L1+operaer-2}) for the QSC-$L1^+$ scheme
can be decomposed as
\begin{equation} \label{FL1+-1}
 \begin{aligned}
  \begin{aligned}
   \bar{\delta}_{t}^{\tilde{\alpha}_{n}} v\left(t_{n-\frac{1}{2}}\right)
   &=\frac{1}{\tau} \int_{t_{n-1}}^{t_{n}} \int_{0}^{t_{n-2}} \partial_{t} \Pi v(s)\ \omega_{1-\tilde{\alpha}_{n }}(t-s)dsdt
    +\frac{1}{\tau} \int_{t_{n-1}}^{t_{n}} \int_{t_{n-2}}^{t} \partial_{t} \Pi v(s)\ \omega_{1-\tilde{\alpha}_{n }}(t-s) dsdt \\
   &:=I_{\tau}^{t_{0}, t_{n-2}}\left(t_{n-\frac{1}{2}}\right)+I_{\tau}^{t_{n-2}, t}\left(t_{n-\frac{1}{2}}\right).
\end{aligned}
 \end{aligned}
\end{equation}

First, for the local term in (\ref{FL1+-1}), it can be computed directly as
\begin{equation} \label{I0-n-1}
 \begin{aligned}
  I_{\tau}^{t_{n-2}, t}\left(t_{n-\frac{1}{2}}\right)=\sum_{k=n-1}^{n} a_{n-k+1}^{(n)}\left(v^{k}-v^{k-1}\right),
 \end{aligned}
\end{equation}
where the coefficients $a_{1}^{(n)}$ and $a_{2}^{(n)}$ can be found in (\ref{L1+-cof}).
Second, based on the definition of $\omega_{1-\beta}(t)$ for the singular kernel, we get the non-local term as
\begin{equation*}
 \begin{aligned}
  I_{\tau}^{t_{0}, t_{n-2}}\left(t_{n-\frac{1}{2}}\right)
  &=\frac{T^{-\tilde{\alpha}_{n }}}{\tau \Gamma(1-\tilde{\alpha}_{n })} \int_{t_{n-1}}^{t_{n}} \int_{0}^{t_{n-1}} \partial_{t} \Pi v(s)\left(\frac{t-s}{T}\right)^{-\tilde{\alpha}_{n }} dsdt.
 \end{aligned}
\end{equation*}
According to Lemma \ref{Fast2}, the term $(\frac{t-s}{T})^{-\tilde{\alpha}_{n }}$ in the integral can be approximated. Then we have
\begin{equation}\label{In-1-t}
 \begin{aligned}
  I_{\tau}^{t_{0}, t_{n-2}}\left(t_{n-\frac{1}{2}}\right)
  &\approx\frac{T^{-\tilde{\alpha}_{n }}}{\tau \Gamma(1-\tilde{\alpha}_{n })} \int_{t_{n-1}}^{t_{n}} \int_{0}^{t_{n-2}} \partial_{t} \Pi v(s)
    \sum_{r=\underline{N}+1}^{\overline{N}} \varpi^{(n,r)} e^{-\lambda^{(r)} \frac{t-s}{T}} dsdt\\
  &=\frac{T^{-\tilde{\alpha}_{n }}}{\tau \Gamma(1-\tilde{\alpha}_{n })} \int_{t_{n-1}}^{t_{n}} \int_{0}^{t_{n-2}} \partial_{t} \Pi v(s)
    \sum_{r=\underline{N}+1}^{\overline{N}} \varpi^{(n,r)} e^{-\lambda^{(r)} \frac{t-t_{n-2}}{T}} e^{-\lambda^{(r)} \frac{t_{n-2}-s}{T}} dsdt\\
  &:=\frac{T^{-\tilde{\alpha}_{n }}}{\tau \Gamma(1-\tilde{\alpha}_{n })}
    \sum_{r=\underline{N}+1}^{\overline{N}} \varpi^{(n,r)} b^{(n,r)} V^{(n,r)},
 \end{aligned}
\end{equation}
where
\begin{equation*}
 \begin{aligned}
  b^{(n,r)}= \int_{t_{n-1}}^{t_{n}} e^{-\lambda^{(r)} \frac{t-t_{n-2}}{T}}dt,
\quad
  V^{(n,r)}=\int_{0}^{t_{n-2}} \partial_t \Pi v(s) e^{-\lambda^{(r)} \frac{t_{n-2}-s}{T}} ds,\quad r=\underline{N}+1,\underline{N}+2,\cdots,\overline{N}.
 \end{aligned}
\end{equation*}

We note that $V^{(2,r)}=0$ and $V^{(n,r)}$ can be got recursively by
\begin{equation} \label{Vni}
 \begin{aligned}
  V^{(n,r)}
  &=\int_{0}^{t_{n-3}} \partial_{t} \Pi v(s) e^{-\lambda^{(r)} \frac{t_{n-2}-s}{T}} ds+\int_{t_{n-3}}^{t_{n-2}} \partial_{t} \Pi v(s) e^{-\lambda^{(r)} \frac{t_{n-2}-s}{T}} ds \\
  &=e^{-\lambda^{(r)} \frac{\tau}{T}} V^{(n-1,r)}+\frac{T}{\lambda^{(r)} \tau}\left(1-e^{-\lambda^{(r)} \frac{\tau}{T}}\right)\left(v^{n-2}-v^{n-3}\right).
 \end{aligned}
\end{equation}
Taking expressions (\ref{I0-n-1}) and (\ref{In-1-t}) into (\ref{FL1+-1}), we can obtain the following fast computational version of $L1^+$ formula
\begin{equation} \label{Fast scheme}
 \begin{aligned}
  \bar{\delta}_{t}^{\tilde{\alpha}_n} v\left(t_{n-\frac{1}{2}}\right)
  =\sum_{k=n-1}^{n} a_{n-k+1}^{(n)}\left(v^{k}-v^{k-1}\right)
   +\frac{T^{-\tilde{\alpha}_{n }}}{\tau \Gamma(1-\tilde{\alpha}_{n})}
    \sum_{r=\underline{N}+1}^{\overline{N}} \varpi^{(n,r)} b^{(n,r)} V^{(n,r)}.
 \end{aligned}
\end{equation}

We use the fast evaluation (\ref{Fast scheme}) for the variable-order fractional operator to take the place of
the $L1^+$ formula in the ADI-QSC-$L1^+$ scheme (\ref{ADI7}),
which results in an improved numerical scheme, say the ADI-QSC-F$L1^+$ scheme. In fact, when $n=1$ and $2$, we
still employ $(\ref{ADI5})$ and $(\ref{ADI7})$ to simulate model (\ref{equation1}).
When $n \geq 3$, the ADI-QSC-F$L1^+$ scheme can be implemented as follows:\\
\textit{Step} 1.
For each $0\leq j\leq M_y +1$, we solve the following one-dimensional linear systems in the $x$ direction for $\tilde{c}_{ij}^{n,*}$
\begin{equation} \label{FQL1}
 \begin{aligned}
  &\left(\boldsymbol\theta_{x}-\gamma_{n} \boldsymbol\eta_{x} \right) \tilde{c}_{ij}^{n,*}\\
  &=\left(1-\frac{b_{2}^{(n)}}{b_{1}^{(n)}} \right)\boldsymbol\theta_{x} \boldsymbol\theta_{y} \tilde{c}_{ij}^{n-1}
    +\frac{b_{2}^{(n)}}{b_{1}^{(n)}}\boldsymbol\theta_{x} \boldsymbol\theta_{y} \tilde{c}_{ij}^{n-2}
    +\gamma_{n}\big(\boldsymbol\eta_{x} \boldsymbol\theta_{y}+\boldsymbol\eta_{y} \boldsymbol\theta_{x}\big) \tilde{c}_{ij}^{n-1}
    +2\gamma_{n}^2 \boldsymbol\eta_{x}  \boldsymbol\eta_{y}  \tilde{c}_{ij}^{n-1}
    -\gamma_{n}^2 \boldsymbol\eta_{x}  \boldsymbol\eta_{y}  \tilde{c}_{ij}^{n-2}\\
    &\quad-\frac{2 \gamma_{n}T^{-\tilde{\alpha}_{n }}}{\tau \kappa \Gamma(1-\tilde{\alpha}_{n})}
    \sum_{r=\underline{N}+1}^{\overline{N}} \varpi^{(n,r)} b^{(n,r)} \boldsymbol\theta_{x} \boldsymbol\theta_{y}\tilde{V}_{ij}^{(n,r)}
   +\frac{2 \gamma_{n}}{\kappa}f_{ij}^{n-\frac{1}{2}},\quad 0\leq i\leq M_x +1,
 \end{aligned}
\end{equation}
where
\begin{equation*}
 \begin{aligned}
  \tilde{V}^{(n,r)}
  =e^{-\lambda^{(r)} \frac{\tau}{T}} \tilde{V}^{(n-1,r)}+\frac{T}{\lambda^{(r)} \tau}\left(1-e^{-\lambda^{(r)} \frac{\tau}{T}}\right)\left(c^{n-2}-c^{n-3}\right),
 \end{aligned}
\end{equation*}
which can be computed similarly with (\ref{Vni}).\\
\textit{Step} 2.
For each $0\leq i\leq M_x +1$, solving the following one-dimensional linear systems in the $y$ direction for $\tilde{c}_{ij}^{n}$
\begin{equation} \label{FQL2}
 \begin{aligned}
  \big(\boldsymbol\theta_{y}-\gamma_{n} \boldsymbol\eta_{y} \big) \tilde{c}_{ij}^{n}= \tilde{c}_{ij}^{n,*},\quad 0\leq j\leq M_y +1 .
 \end{aligned}
\end{equation}

With the expected accuracy $\epsilon \leq \mathcal{O}\big( \tau^{\min{\{3-\alpha^*-\alpha(0),2\}}} \big)$, the fast scheme requires $\mathcal{O}(n\log^{2} n )$
computational cost to approximate the variable-order Caputo fractional derivative.
Therefore, for the ADI-QSC-F$L1^+$ scheme (\ref{FQL1})-(\ref{FQL2}) for model (\ref{equation1}),
the computational cost is reduced from $\mathcal{O}(M_{x}M_{y}N^2)$ to $\mathcal{O}(M_{x}M_{y}N\log^{2} N)$.
Moreover, with the fast evaluation scheme (\ref{Fast scheme}), the storage requirement is reduced from for $\mathcal{O}(M_{x}M_{y}N)$
of the ADI-QSC-$L1^+$ scheme to $\mathcal{O}(M_{x}M_{y}\log^{2} N)$ for the ADI-QSC-F$L1^+$ scheme.

\subsection{Optimal QSC method}
In this subsection, we consider the acceleration in space domain. In fact,
the standard QSC method can be improved by introducing high order perturbations,
which leads to the optimal QSC method with fourth-order spatial convergence order. Therefore, we assume the solution of model (\ref{equation1})-(\ref{equation3})
satisfies $u(x,y,\cdot)\in C^6\big(\bar{\Omega}\big)$.
Therefore, we introduce the perturbation $\mathcal{P}_{\mathcal{S}x}$ as
\begin{equation*}
 \begin{aligned}
 \mathcal{P}_{\mathcal{S}x} \tilde{c}_{ij}
  =\frac{1}{24\Delta x^2}
  \left\{
   \begin{aligned}
    &0,&i= 0,\\
    & -11\tilde{c}_{1j}      + 16  \tilde{c}_{2j}        - 14 \tilde{c}_{3j}         +  6  \tilde{c}_{4j}         - \tilde{c}_{5j}          , &i= 1,\\
    & -5 \tilde{c}_{1j}      + 6   \tilde{c}_{2j}        - 4  \tilde{c}_{3j}         +     \tilde{c}_{4j}                                   , &i= 2,\\
    &    \tilde{c}_{i-2,j}   - 4   \tilde{c}_{i-1,j}     + 6  \tilde{c}_{ij}         -  4  \tilde{c}_{i+1,j}      + \tilde{c}_{i+2,j}       , &i= 3, \ldots, M_{x}-2,\\
    & -5 \tilde{c}_{M_{x},j} + 6   \tilde{c}_{M_{x}-1,j} - 4  \tilde{c}_{M_{x}-2,j}  +     \tilde{c}_{M_{x}-3,j}                            , &i= M_{x}-1,\\
    & -11\tilde{c}_{M_{x}j}  + 16  \tilde{c}_{M_{x}-1,j} - 14 \tilde{c}_{M_{x}-2,j}  +  6  \tilde{c}_{M_{x}-3,j}  - \tilde{c}_{M_{x}-4,j}   , &i= M_{x},\\
    &0,&i= M_{x}+1,\\
   \end{aligned}
  \right.
 \end{aligned}
\end{equation*}
and the perturbation $\mathcal{P}_{\mathcal{S}y}$ can be defined similarly. Specially, the derivation of the expressions of
$\mathcal{P}_{\mathcal{S}x}$ and $\mathcal{P}_{\mathcal{S}y}$ can be found in \cite{PSoper2,PSoper1} for detail.
We only need to take the perturbations
$\mathcal{P}_{\mathcal{S}x}$ and $\mathcal{P}_{\mathcal{S}y}$ together with $\boldsymbol\eta_{x}$ and $\boldsymbol\eta_{y}$, respectively, for improving the
accuracy of spatial approximation. Taking the optimal QSC framework into the ADI-QSC-$L1^+$ scheme, we can get the optimal ADI-QSC-$L1^+$ scheme.
Without loss of generality, we consider the optimal ADI-QSC-F$L1^+$ scheme with fast computation directly as follows.

(I) For $n=1$, 
\begin{equation} \label{FADI1}
	\begin{aligned}
		&\Big[\boldsymbol\theta_{x}-\gamma_{1} \big(\boldsymbol\eta_{x}+\mathcal{P}_{\mathcal{S}x} \big)\Big]\Big[\boldsymbol\theta_{y}-\gamma_{1} \big(\boldsymbol\eta_{y}+\mathcal{P}_{\mathcal{S}y} \big) \Big] \tilde{c}_{ij}^{1} \\
		&=\Big[\boldsymbol\theta_{x}+\gamma_{1}\big(\boldsymbol\eta_{x}+\mathcal{P}_{\mathcal{S}x} \big)\Big] \Big[\boldsymbol\theta_{y}
		+\gamma_{1} \big(\boldsymbol\eta_{y}+\mathcal{P}_{\mathcal{S}y} \big) \Big] {c}_{ij}^{0}
		+\frac{2 \gamma_{1}}{\kappa}f_{ij}^{\frac{1}{2}},\quad (i,j)\in \bar{\Lambda}.
	\end{aligned}
\end{equation}

(II) For $n=2$, 
\begin{equation} \label{FADI3}
	\begin{aligned}
		&\Big[\boldsymbol\theta_{x}-\gamma_{2} \big(\boldsymbol\eta_{x}+\mathcal{P}_{\mathcal{S}x} \big)\Big] \Big[\boldsymbol\theta_{y}-\gamma_{2} \big(\boldsymbol\eta_{y}+\mathcal{P}_{\mathcal{S}y} \big) \Big] \tilde{c}_{ij}^{2}\\
		&=\boldsymbol\theta_{x} \boldsymbol\theta_{y} \tilde{c}_{ij}^{1}
		+\gamma_{2}\Big[\big(\boldsymbol\eta_{x} +\mathcal{P}_{\mathcal{S}x}\big)\boldsymbol\theta_{y}
		+\big(\boldsymbol\eta_{y} +\mathcal{P}_{\mathcal{S}y}\big)\boldsymbol\theta_{x}\Big] \tilde{c}_{ij}^{1}
		+2\gamma_{2}^2 \big(\boldsymbol\eta_{x} +\mathcal{P}_{\mathcal{S}x} \big)\big(\boldsymbol\eta_{y} +\mathcal{P}_{\mathcal{S}y} \big)  \tilde{c}_{ij}^{1}\\
		&\quad -\gamma_{2}^2 \big(\boldsymbol\eta_{x} +\mathcal{P}_{\mathcal{S}x} \big)\big(\boldsymbol\eta_{y} +\mathcal{P}_{\mathcal{S}y} \big) {c}_{ij}^{0}
		-\frac{2 \gamma_{2}}{\kappa} a_{2}^{(2)}\left(\boldsymbol\theta_{x}\boldsymbol\theta_{y} \tilde{c}_{ij}^{1}
		-\boldsymbol\theta_{x} \boldsymbol\theta_{y} {c}_{ij}^{0}\right)
		+\frac{2 \gamma_{2}}{\kappa}f_{ij}^{\frac{3}{2}}, \  (i,j)\in \bar{\Lambda}.
	\end{aligned}
\end{equation}

(III) For $n\geq 3$, 
\begin{equation} \label{eqOPT1}
	\begin{aligned}
		&\Big[\boldsymbol\theta_{x}-\gamma_{n} \big(\boldsymbol\eta_{x} +\mathcal{P}_{\mathcal{S}x}\big)\Big]  \Big[\boldsymbol\theta_{y}-\gamma_{n} \big(\boldsymbol\eta_{y} +\mathcal{P}_{\mathcal{S}y}\big)\Big] \tilde{c}_{ij}^{n}\\
		&=\boldsymbol\theta_{x} \boldsymbol\theta_{y} \tilde{c}_{ij}^{n-1}
		+\gamma_{n}\Big[\big(\boldsymbol\eta_{x} +\mathcal{P}_{\mathcal{S}x}\big)\boldsymbol\theta_{y}
		+\big(\boldsymbol\eta_{y} +\mathcal{P}_{\mathcal{S}y}\big)\boldsymbol\theta_{x}\Big]\tilde{c}_{ij}^{n-1}\\
		&\quad +2\gamma_{n}^2 \big(\boldsymbol\eta_{x} +\mathcal{P}_{\mathcal{S}x} \big)\big(\boldsymbol\eta_{y} +\mathcal{P}_{\mathcal{S}y} \big)\tilde{c}_{ij}^{n-1}
		-\gamma_{n}^2 \big(\boldsymbol\eta_{x} +\mathcal{P}_{\mathcal{S}x} \big)\big(\boldsymbol\eta_{y} +\mathcal{P}_{\mathcal{S}y} \big) \tilde{c}_{ij}^{n-2}\\
		&\quad -\frac{2 \gamma_{n}T^{-\tilde{\alpha}_{n }}}{\tau \kappa \Gamma(1-\tilde{\alpha}_{n})}
		\sum_{r=\underline{N}+1}^{\overline{N}} \varpi^{(n,r)} b^{(n,r)} \boldsymbol\theta_{x} \boldsymbol\theta_{y}\tilde{V}_{ij}^{(n,r)}
		+\frac{2 \gamma_{n}}{\kappa}f_{ij}^{n-\frac{1}{2}},\  (i,j)\in \bar{\Lambda}.
	\end{aligned}
\end{equation}

The optimal ADI-QSC-F$L1^+$ scheme can achieve fourth-order accuracy in space, which means that we can get a desired accuracy with much less mesh grids.

\section{Numerical Experiments}
\label{sect:numerical}
 In this section, we consider numerical experiments to numerically support the accuracy and efficiency of schemes developed in this paper.
All schemes are programmed in Matlab R2018b, and implemented on a Windows server with Intel(R) Xeon(R) E5-2650 CPU @ 2.30 GHz.
We consider model (\ref{equation1})-(\ref{equation3}) in the space domain $\Omega =(0,1)\times (0,1)$ and time interval $[0,1]$,
and choose four different variable time fractional order
\begin{equation*}
 \begin{aligned}
  &\alpha_{0}(t)=0.45-0.3t,
   \quad \alpha_{1}(t)=0.4+0.5(1-t)-\frac{1}{4\pi}\left[sin\big(2\pi(1-t)\big)\right],\\
  &\alpha_{2}(t)=0.8-0.5(1-t),
   \quad \alpha_{3}(t)=\left|3(t-0.5)^2-0.2\right|+0.3.
 \end{aligned}
\end{equation*}

\textbf{Example 6.1.} \label{ex-1}
  We choose the diffusivity coefficient $\kappa = 1$,
  the initial date and the source function as
\begin{equation*}
  	\begin{aligned}
  		&u^{0}(x, y)=\sin x \sin y, \quad f\left(x,y,t \right)=0.
  	\end{aligned}
\end{equation*}
 The true solution of model (\ref{equation1})-(\ref{equation3}) is unknown. 
 
 We first fix the values of $M_x = M_y$ big enough to investigate the temporal errors and convergence orders. We compute numerical solutions $u_h^n(x,y)$ on a coarse time mesh with size $\tau$, then we refine the time mesh with sizes $\tau/2$. The resulting errors in the discrete $L_{2}$ -norm at time $t_n$ can be calculated on the coarse mesh as
\begin{equation*}
	\begin{aligned}
	Err^2:= \Delta x \Delta y \sum_{i=0}^{M_{x}+1 } \sum_{j=0}^{M_{y}+1 }
	\left|u_{h}^{n}\big(\xi_{i}^{x},\xi_{j}^{y}\big)-u_{h}^{2n}\big(\xi_{i}^{x},\xi_{j}^{y}\big)\right|^2 .
\end{aligned}
\end{equation*}
We show the observing errors and the temporal convergence orders at the time instance near the initial time and at time instance $t=T$, respectively, in Table \ref{Ex2-Table0} and Table \ref{Ex2-Table1}. We can see that, near the initial time, the fractional index $\alpha_{0}(t)$ satisfies the assumptions in Corollary \ref{corollary-1} and the local temporal convergence order is preserved. But the fractional indices $\alpha_{1}(t)$ and $\alpha_{3}(t)$ do not satisfy the assumptions, and the local temporal convergence orders are about $\mo\big(\tau^{3-\alpha^{*}-\alpha(0)}\big)$. On the other hand, the singularity near the initial time almost do not affect the convergence order at the final time point, which reaches a satisfied second-order.
 
Next, we choose the values of $N$ big enough. Using the similar routine by the coarse and fine space mesh, we can get computing errors and convergence orders in Table \ref{Ex2-Table2}, which fit the theoretical spatial convergence orders well.

 \begin{table}[htbp!]
 	\centering
 	\caption{\small{Errors and temporal convergence orders of QSC-$L1^+$ near the initial time point, with  $M_x=M_y=2^{12}$ }}\label{Ex2-Table0}
 	\begin{tabular}{ccccccc}
 		\hline
 		\multicolumn{1}{c}{}&\multicolumn{2}{c}{$\alpha_{0}(t)$}   & \multicolumn{2}{c}{$\alpha_{1}(t)$}   & \multicolumn{2}{c}{$\alpha_{3}(t)$}                     \\ 
 		\hline
 		$\tau$        & $Err$  & $Order$               & $Err$  & $Order$                                & $Err$  & $Order$                                  \\
 		\hline\hline
 		$2^{-9}$     & 5.01e-06     & ---              &     1.10e-05       & ---     &            1.62e-05             & ---                                         \\
 		$2^{-10}$    & 8.71e-07     & 2.52             &     4.17e-06       & 1.39    &            5.81e-06             & 1.47                                        \\
 		$2^{-11}$    & 1.57e-07     & 2.47             &     1.72e-06       & 1.27    &            2.26e-06             & 1.36                                        \\
 		\hline\hline
 		&              &$\approx$ 2.10   &            &             $\approx$ 1.20   &            &                   $\approx$ 1.30                             \\
 		\hline
 	\end{tabular}
 \end{table}
 \begin{table}[htbp!]
 	\centering
 	\caption{\small{Errors and temporal convergence orders of QSC-$L1^+$ at the final time point, with  $M_x=M_y=2^{11}$ }}\label{Ex2-Table1}
 	\begin{tabular}{ccccccc}
 		\hline
 		\multicolumn{1}{c}{}&\multicolumn{2}{c}{$\alpha_{0}(t)$}   & \multicolumn{2}{c}{$\alpha_{1}(t)$}   & \multicolumn{2}{c}{$\alpha_{3}(t)$}                     \\ 
 		\hline
 		$\tau$        & $Err$  & $Order$               & $Err$  & $Order$                                & $Err$  & $Order$                                  \\
 		\hline\hline
 		$2^{-4}$    & 6.78e-06     & ---              &     4.23e-06       & ---     &            4.27e-05             & ---                                         \\
 		$2^{-5}$    & 1.71e-06     & 1.98             &     9.70e-07       & 2.12    &            1.03e-05             & 2.05                                        \\
 		$2^{-6}$    & 4.27e-07     & 2.00             &     2.37e-07       & 2.03    &            3.08e-06             & 1.74                                        \\
 		$2^{-7}$    & 1.07e-07     & 1.99             &     5.89e-08       & 2.00    &            5.95e-07             & 2.37                                        \\
 		\hline\hline
 		&              &$\approx$ 2.00   &            &             $\approx$ 2.00   &            &                   $\approx$ 2.00                              \\
 		\hline
 	\end{tabular}
 \end{table}
 \begin{table}[htbp!]
	\centering
	\caption{\small{Errors and spatial convergence orders of QSC-$L1^+$ for Example 6.1, with $N=2^{11}$ }}\label{Ex2-Table2}
	\begin{tabular}{ccccccc}
		\hline
		\multicolumn{1}{c}{}&\multicolumn{2}{c}{$\alpha_{0}(t)$}  & \multicolumn{2}{c}{$\alpha_{1}(t)$}   & \multicolumn{2}{c}{$\alpha_{3}(t)$}     \\ 
		\hline
		$\Delta{x}=\Delta{y}$        & $Err$  & $Order$          & $Err$  & $Order$               & $Err$  & $Order$                                           \\
		\hline\hline
		$2^{-4}$    & 1.28e-05     & ---                &        1.08e-05       & ---       &     8.87e-06             & ---                                         \\
		$2^{-5}$    & 3.20e-06     & 2.00               &        2.70e-06       & 2.00      &     2.23e-06             & 1.99                                        \\
		$2^{-6}$    & 8.01e-07     & 1.99               &        6.75e-07       & 2.00      &     5.57e-07             & 2.00                                        \\
		$2^{-7}$    & 2.00e-07     & 2.00               &        1.69e-07       & 1.99      &     1.39e-07             & 2.00                  \\
		\hline\hline
		&              &$\approx$ 2.00   &                           & $\approx$ 2.00      &                           & $\approx$ 2.00                        \\
		\hline
	\end{tabular}
\end{table}

\textbf{Example 6.2.} 
We choose the diffusivity coefficient $\kappa = 1$,
the initial date and the source function as
\begin{equation*}
 \begin{aligned}
  &u^{0}(x, y)=\sin x \sin y,\\
  &f\left(x,y,t \right)=\left[3 t^{2}+\frac{\Gamma(4)}{\Gamma\big(4-\alpha(t)\big)} t^{3-\alpha(t)}+2 \pi^{2}\left(1+t^{3}\right)\right]\sin(\pi x)\sin(\pi y),
 \end{aligned}
\end{equation*}
such that the true solution of model (\ref{equation1})-(\ref{equation3}) is
\begin{equation*}
 \begin{aligned}
  u(x, y, t)=\Big(1+t^{3}\Big) \sin (\pi x) \sin (\pi y).
 \end{aligned}
\end{equation*}
The error is measured in the discrete $L_{2}$ -norm as
\begin{equation*}
	\begin{aligned}
		Err^2:= \Delta x \Delta y \sum_{i=0}^{M_{x}+1 } \sum_{j=0}^{M_{y}+1 }
		\left|u_{h}^{n}\big(\xi_{i}^{x},\xi_{j}^{y}\big)-u^{n}\big(\xi_{i}^{x},\xi_{j}^{y}\big)\right|^2 .
	\end{aligned}
\end{equation*}

We first fix $N=2^{11}$, and the observation errors and spatial convergence orders of the QSC-$L1^+$ scheme, the ADI-QSC-$L1^+$ scheme and the ADI-QSC-F$L1^+$ scheme are shown in Table \ref{Ex1-Table1}. Then, we fix $M_{x}=M_{y}=2^{11}$, and the temporal convergence orders are shown in Table \ref{Ex1-Table2}. It can be seen that, all the three schemes have second-order convergence orders in both space and time, which conform the theoretical convergence orders in Remark \ref{remark-2}. Due to the high convergence order, the optimal ADI-QSC-$L1^+$ scheme is considered separately.
we fix $N=2^{17}$ to observe the spatial convergence orders in Table \ref{Ex1-Table3}, and fix
$M_{x}=M_{y}=2^{5}$ to get the temporal convergence orders in Table \ref{Ex1-Table4}.
The results show that the optimal ADI-QSC-$L1^+$ scheme has fourth-order convergence orders in space, which is consistent with
the theoretical results.

\begin{table}[htbp!]
  \centering
  \caption{\small{Errors and spatial convergence orders of proposed schemes for Example 6.2, with $N=2^{11}$}}      \label{Ex1-Table1}
  \begin{tabular}{cccccccc}
  \hline
              \multicolumn{2}{c}{}&\multicolumn{2}{c}{QSC-$L1^+$}                 &\multicolumn{2}{c}{ADI-QSC-$L1^+$}           & \multicolumn{2}{c}{ADI-QSC-F$L1^+$} \\
  \hline
               &$\Delta{x}=\Delta{y}$        & $Err$  & $Order$   &               $Err$  & $Order$                  &$Err$  & $Order$     \\
   \hline\hline
               &$2^{-4}$       & 1.42e-03     & ---             &                 1.42e-03      & ---            &      1.42e-03                        & ---                  \\
$\alpha_{1}(t)$&$2^{-5}$       & 3.55e-04     & 2.00            &                 3.55e-04      & 2.00           &      3.55e-04                        & 2.00                  \\
               &$2^{-6}$       & 8.89e-05     & 1.99            &                 8.89e-05      & 1.99           &      8.88e-05                        & 1.99                  \\
               &$2^{-7}$       & 2.23e-05     & 1.99            &                 2.23e-05      & 1.99           &      2.22e-05                        & 2.00                  \\
  \hline
               &$2^{-4}$       & 1.41e-03     & ---             &                 1.41e-03      & ---            &      1.41e-03                        & ---                  \\
$\alpha_{2}(t)$&$2^{-5}$       & 3.53e-04     & 1.99            &                 3.53e-04      & 1.99           &      3.53e-04                        & 1.99                  \\
               &$2^{-6}$       & 8.84e-05     & 1.99            &                 8.84e-05      & 1.99           &      8.84e-05                        & 1.99                  \\
               &$2^{-7}$       & 2.22e-05     & 1.99            &                 2.22e-05      & 1.99           &      2.22e-05                        & 1.99                  \\
  \hline
               &$2^{-4}$       & 1.41e-03     & ---             &                 1.41e-03      & ---            &      1.41e-03                        & ---                  \\
$\alpha_{3}(t)$&$2^{-5}$       & 3.53e-04     & 1.99            &                 3.53e-04      & 1.99           &      3.53e-04                        & 1.99                  \\
               &$2^{-6}$       & 8.84e-05     & 1.99            &                 8.84e-05      & 1.99           &      8.84e-05                        & 1.99                  \\
               &$2^{-7}$       & 2.22e-05     & 1.99            &                 2.22e-05      & 1.99           &      2.22e-05                        & 1.99                  \\
  \hline\hline
               &               &              &$\approx$ 2.00   &                               & $\approx$ 2.00 &                              & $\approx$ 2.00  \\
  \hline
\end{tabular}
\end{table}

\begin{table}[htbp!]
  \centering
  \caption{\small{Errors and temporal convergence orders of proposed schemes for Example 6.2, with  $M_x=M_y=2^{11}$}}      \label{Ex1-Table2}
  \begin{tabular}{cccccccc}
  \hline
  \multicolumn{2}{c}{}       &\multicolumn{2}{c}{QSC-$L1^+$}                 &\multicolumn{2}{c}{ADI-QSC-$L1^+$}           & \multicolumn{2}{c}{ADI-QSC-F$L1^+$}         \\
  \hline
               &$\tau$        & $Err$  & $Order$   &               $Err$  & $Order$                  &$Err$  & $Order$                   \\
   \hline\hline
               &$2^{-4}$       & 1.42e-03     & ---             &                 2.11e-03      & ---            &                2.02e-03              & ---                \\
$\alpha_{1}(t)$&$2^{-5}$       & 3.54e-04     & 2.00            &                 4.46e-04      & 2.24           &                4.26e-04              & 2.24                 \\
               &$2^{-6}$       & 8.86e-05     & 1.99            &                 1.00e-04      & 2.15           &                9.65e-05              & 2.14                 \\
               &$2^{-7}$       & 2.22e-05     & 1.99            &                 2.37e-05      & 2.07           &                2.28e-05              & 2.08                 \\
  \hline
               &$2^{-4}$       & 1.40e-03     & ---             &                 1.95e-03      & ---            &                1.94e-03              & ---                 \\
$\alpha_{2}(t)$&$2^{-5}$       & 3.51e-04     & 1.99            &                 4.24e-04      & 2.20           &                4.22e-04              & 2.20                  \\
               &$2^{-6}$       & 8.80e-05     & 1.99            &                 9.76e-05      & 2.11           &                9.73e-05              & 2.11                  \\
               &$2^{-7}$       & 2.21e-05     & 1.99            &                 2.33e-05      & 2.06           &                2.33e-05              & 2.06                  \\
  \hline

               &$2^{-4}$       & 2.81e-03     & ---             &                 1.99e-03      & ---            &                1.98e-03              & ---                  \\
$\alpha_{3}(t)$&$2^{-5}$       & 7.04e-04     & 1.99            &                 4.29e-04      & 2.21           &                4.28e-04              & 2.20                  \\
               &$2^{-6}$       & 1.76e-04     & 2.00            &                 9.82e-05      & 2.12           &                9.81e-05              & 2.12                  \\
               &$2^{-7}$       & 4.42e-05     & 1.99            &                 2.34e-05      & 2.06           &                2.34e-05              & 2.06                  \\
  \hline\hline
               &               &              &$\approx$ 2.00   &                               & $\approx$ 2.00 &                              & $\approx$ 2.00  \\
  \hline
\end{tabular}
\end{table}

\begin{table}[htbp!]
  \centering
  \caption{\small{Errors and spatial convergence orders of optimal ADI-QSC-$L1^+$ for Example 6.2, with $N=2^{17}$}}\label{Ex1-Table3}
  \begin{tabular}{ccccccccc}
  \hline
   \multicolumn{1}{c}{} &\multicolumn{2}{c}{$\alpha_{1}(t)$}     & \multicolumn{2}{c}{$\alpha_{2}(t)$}    & \multicolumn{2}{c}{$\alpha_{3}(t)$} \\ \hline
   $\Delta x=\Delta y$                  & $Err$       & $Order$       & $Err$       & $Order$     &  $Err$  & $Order$     \\
   \hline\hline
  $2^{-4}$              & 1.25e-06    &  ---          & 1.26e-06    & ---         & 1.25e-06    & ---                     \\
  $2^{-5}$              & 5.98e-08    & 4.38          & 5.99e-08    & 4.39        & 5.95e-08    & 4.39                    \\
  $2^{-6}$              & 3.57e-09    & 4.06          & 3.58e-09    & 4.06        & 3.55e-09    & 4.07                    \\
  $2^{-7}$              & 2.02e-10    & 4.14          & 2.03e-10    & 4.14        & 2.01e-10    & 4.14                    \\

  \hline\hline
                        &              &$\approx$ 4.00   &          & $\approx$ 4.00 &  & $\approx$ 4.00  \\
  \hline
\end{tabular}
\end{table}

\begin{table}[htbp!]
  \centering
  \caption{\small{Errors and temporal convergence orders of optimal ADI-QSC-$L1^+$ for Example 6.2, $M_x=M_y=2^{5}$}}\label{Ex1-Table4}
  \begin{tabular}{cccccccccc}
  \hline
   \multicolumn{1}{c}{}&\multicolumn{2}{c}{$\alpha_{1}(t)$}     & \multicolumn{2}{c}{$\alpha_{2}(t)$}    & \multicolumn{2}{c}{$\alpha_{3}(t)$} \\ \hline
   $\tau$        & $Err$  & $Order$               & $Err$  & $Order$               &$Err$  & $Order$     \\
   \hline\hline
  $2^{-4}$    & 2.11e-03     & ---                  &            1.95e-03             & ---                       &  1.99e-03             & ---                  \\
  $2^{-5}$    & 4.45e-04     & 2.24                 &            4.24e-04             & 2.20                      &  4.29e-04             & 2.21                  \\
  $2^{-6}$    & 1.00e-04     & 2.15                 &            9.74e-05             & 2.12                      &  9.81e-05             & 2.12                  \\
  $2^{-7}$    & 2.36e-05     & 2.08                 &            2.32e-05             & 2.06                      &  2.33e-05             & 2.07                  \\
  \hline\hline
              &              &$\approx$ 2.00   &                           & $\approx$ 2.00 &                          & $\approx$ 2.00  \\
  \hline
\end{tabular}
\end{table}

\textbf{Example 6.3.} 
We choose the same equations in Example 6.2 to show the effectiveness of acceleration techniques.
We present the comparisons among the QSC-$L1^+$ scheme, the ADI-QSC-$L1^+$ scheme,
the ADI-QSC-F$L1^+$ scheme, the optimal ADI-QSC-$L1^+$ scheme,
and the optimal ADI-QSC-F$L1^+$ scheme.
For the first three schemes, the observation errors and CPU time can be found in Table \ref{Ex1-Table5}.
We can see that, the ADI method and the ESA technique can improve the efficiency greatly. For example, for $M_{x}=M_{y}=N=2^{11}$ in the case of $\alpha_{2}(t)$,
ADI method can reduce the CPU time from $100197$ seconds to $58852$ seconds, and the fast evaluation can further reduce the CPU time
to $26782$ seconds, while almost preserving the same accuracy.
For the last two optimal QSC based schemes, the observation errors and CPU time can be found in Table \ref{Ex1-Table6}.
It can be seen that, the fourth spatial convergence order by the optimal QSC method allows us to employ much more sparse meshes which will
save enormous computational cost, though higher smooth assumption for the solution is required.

\begin{table}[htbp!]
  \centering
  \tabcolsep=0.15cm
  \caption{\small{Errors and CPU time of proposed schemes for Example 6.3}}\label{Ex1-Table5}
  \begin{tabular}{cccccccc}
  \hline
  \multicolumn{2}{c}{}&\multicolumn{2}{c}{QSC-$L1^+$} &\multicolumn{2}{c}{ADI-QSC-$L1^+$ } &\multicolumn{2}{c}{ADI-QSC-F$L1^+$}      \\
  \hline
               &$\tau=\Delta x=\Delta y$        & $Err$        &     $Time(Sec.)$       & $Err$        &   $Time(Sec.)$    & $Err$        &$Time(Sec.)$               \\
  \hline\hline
               &$2^{-8} $       & 2.22e-05     &       72.4       & 2.26e-05     &      9.3     & 2.20e-05     &     13.2           \\
$\alpha_{1}(t)$&$2^{-9} $       & 5.54e-06     &      719.1       & 5.59e-06     &    215.8     & 5.47e-06     &    202.8           \\
               &$2^{-10}$       & 1.39e-06     &     7449.6       & 1.39e-06     &   3145.4     & 1.30e-06     &   1930.1           \\
               &$2^{-11}$       & 3.46e-07     &    98969         & 3.47e-07     &  59383       & 2.29e-07     &  22452             \\
  \hline
               &$2^{-8} $       & 2.20e-05     &       73.8       & 2.24e-05     &     18.7     & 2.24e-05     &     15.4           \\
$\alpha_{2}(t)$&$2^{-9} $       & 5.51e-06     &      735.7       & 5.55e-06     &    212.4     & 5.55e-06     &    237.6           \\
               &$2^{-10}$       & 1.38e-06     &     7579.7       & 1.38e-06     &   2962.3     & 1.38e-06     &   2298.2           \\
               &$2^{-11}$       & 3.45e-07     &   100197         & 3.46e-07     &  58852       & 3.45e-07     &  26782             \\
  \hline
               &$2^{-8} $       & 2.21e-05     &       74.8       & 2.24e-05     &      9.4     & 2.24e-05     &     16.1           \\
$\alpha_{3}(t)$&$2^{-9} $       & 5.52e-06     &      730.5       & 5.56e-06     &    209.5     & 5.53e-06     &    237.3           \\
               &$2^{-10}$       & 1.38e-06     &     7544.0       & 1.39e-06     &   3036.6     & 1.37e-06     &   2322.8           \\
               &$2^{-11}$       & 3.45e-07     &   100328         & 3.46e-07     &  59004       & 3.24e-07     &  27161             \\
  \hline\hline
\end{tabular}
\end{table}

\begin{table}[htbp!]
  \centering
  \caption{\small{Errors and CPU time of optimal QSC based schemes for Example 6.3}}\label{Ex1-Table6}
  \begin{tabular}{ccccccc}
  \hline
  \multicolumn{3}{c}{} &\multicolumn{2}{c}{optimal ADI-QSC-$L1^+$}&\multicolumn{2}{c}{optimal ADI-QSC-F$L1^+$}\\
  \hline
               &$\tau$          &$\Delta x=\Delta y$          & $Err$  &       $Time(Sec.)$   & $Err$  &       $Time(Sec.)$     \\
  \hline\hline
               &$2^{-8} $       &$2^{-4} $                    &4.49e-06    &  0.2             &4.24e-06    &  0.2    \\
$\alpha_{1}(t)$&$2^{-9} $       &$2^{-5} $                    &1.35e-06    &  1.3             &1.29e-06    &  1.1    \\
               &$2^{-10}$       &$2^{-5} $                    &2.89e-07    &  3.7             &2.19e-07    &  2.3    \\
               &$2^{-11}$       &$2^{-6} $                    &8.33e-08    & 34.0             &2.34e-08    & 16.6    \\
  \hline
               &$2^{-8} $       &$2^{-4} $                    &4.44e-06    &  0.2             &4.44e-06    &  0.2    \\
$\alpha_{2}(t)$&$2^{-9} $       &$2^{-5} $                    &1.34e-06    &  1.3             &1.34e-06    &  1.2    \\
               &$2^{-10}$       &$2^{-5} $                    &2.88e-07    &  3.6             &2.88e-07    &  2.5    \\
               &$2^{-11}$       &$2^{-6} $                    &8.29e-08    & 33.8             &8.29e-08    & 18.3    \\
  \hline
               &$2^{-8} $       &$2^{-4} $    &4.45e-06   &  0.2      &4.44e-06    &  0.2    \\
$\alpha_{3}(t)$&$2^{-9} $       &$2^{-5} $    &1.34e-06   &  1.3      &1.32e-06    &  1.2    \\
               &$2^{-10}$       &$2^{-5} $    &2.88e-07   &  3.5      &2.65e-07    &  2.5    \\
               &$2^{-11}$       &$2^{-6} $    &8.30e-08   & 33.3      &4.97e-08    & 18.4    \\
  \hline\hline
\end{tabular}
\end{table}

\section{Conclusions}
\label{sect:conclusions}

In this paper, we develop the QSC-$L1^+$ scheme for the variable-order TF-MID equation in two dimensional space domain.
The scheme is proved to be unconditionally stable and convergent with accuracy $\mathcal{O}\big( \tau^{\min{\{3-\alpha^*-\alpha(0),2\}}} +\Delta{x}^2+\Delta{y}^2\big)$,
for proper assumptions on $\alpha(t)$.
Based on the QSC-$L1^+$ scheme, we design a novel ADI framework to achieve the ADI-QSC-$L1^+$ scheme, and we also analyze
its unconditional stability and convergence.
The numerical experiments show that the results fit well with the theoretical analysis, even $\alpha(t)$ do not satisfy the restrictions.
Then we employ the fast evaluation based on ESA technique to obtain the ADI-QSC-F$L1^+$ scheme, which can reduce the computation cost greatly.
Furthermore, the optimal QSC method is also applied to get the optimal ADI-QSC-F$L1^+$ scheme, the numerical results
show that the higher order schemes lead to much better computational efficiency.

\vskip +5pt
\small
\noindent \tbf{Funding} The work of J. Liu was supported in part by the Shandong Provincial Natural Science Foundation (Nos. ZR2021MA020, ZR2020MA039), the Fundamental Research Funds for the Central Universities (Nos. 22CX03016A, 20CX05011A), and the Major Scientific and Technological Projects of CNPC under Grant (No. ZD2019-184-001). The work of H. Fu was supported in part by the National Natural Science Foundation of China (Nos. 11971482, 12131014), the Fundamental Research Funds for the Central Universities (No. 202264006), and by the OUC Scientific Research Program for Young Talented Professionals.

\noindent \tbf{Data availability}  Enquiries about data availability should be directed to the authors.

\section*{Declarations} \small
\noindent  \tbf{Conflict of interest}  The authors declare that they have no conflict of interest.

\appendix
	\section{Estimate of $r_{1,n}$} \label{r1n}

Based on the definition of $r_{1, n}$, we have
\begin{equation*}
\begin{aligned}
	r_{1, n} 
	& =\frac{1}{\tau} \int_{t_{n-1}}^{t_{n}} {}_0^C D_{t}^{\alpha(t)} v(t) d t-\frac{1}{\tau} \int_{t_{n-1}}^{t_{n}}{}_0^C D_{t}^{\tilde{\alpha}_{n}} v(t)  dt \\
	& =\frac{1}{\tau} \int_{t_{n-1}}^{t_{n}} \int_{0}^{t} \left[\omega_{1-\alpha(t)}(t-s)-\omega _{1-\tilde{\alpha}_{n}} (t-s)\right]\partial_{s} v(s) dsdt.
\end{aligned}
\end{equation*}
We can obtain from Lemma \ref{ZhengZe} that
\begin{equation*}
	\begin{aligned}
\left|r_{1, n}\right| 
&\leq 
Q_0\left|\frac{1}{\tau} \int_{t_{n-1}}^{t_{n}}\int_{0}^{t} \big[\omega_{1-\alpha(t)}(t-s)-\omega_{1-\tilde{\alpha}_{n}} (t-s)\big] d s d t\right| \\
&=
Q_0\left|\frac{1}{\tau} \int_{t_{n-1}}^{t_{n}} \left[\frac{t^{1-\alpha(t)}}{\Gamma\left(2-\alpha(t)\right)}-\frac{t^{1-\tilde{\alpha}_{n}}}{\Gamma\left(2-\tilde{\alpha}_{n}\right)} \right] d t\right| \\
&=
Q_0\left|\frac{1}{\tau} \int_{t_{n-1}}^{t_{n}} \frac{\Gamma\left(2-\tilde{\alpha}_{n}\right) t^{1-\alpha(t)}-\Gamma(2-\alpha(t)) t^{1-\tilde{\alpha}_{n}}}{\Gamma(2-\alpha(t)) \Gamma\left(2-\tilde{\alpha}_{n}\right)} d t\right| \\
&\leq 
Q_0Q_1\left|\frac{1}{\tau} \int_{t_{n-1}}^{t_{n}} \Gamma\left(2-\tilde{\alpha}_{n}\right)\left(t^{1-\alpha(t)}
-t^{1-\tilde{\alpha}_{n}}\right) d t\right|
+Q_0 Q_{1} \left| \frac{1}{\tau} \int_{t_{n-1}}^{t_{n}} t^{1-\tilde{\alpha}_{n}}\big[\Gamma\left(2-\tilde{\alpha}_{n}\right)-\Gamma(2-\alpha(t))\big]d t \right|,
	\end{aligned}
\end{equation*}
where $\Gamma(x)$ is bounded when $x\in(1,2)$.
By Taylor's expansion, we can get 
\begin{equation*}
	\begin{aligned}
\left.t^{1-\alpha(\xi)}\right|_{\xi=t}
=&
t^{1-\tilde{\alpha}_{n}}-t^{1-\tilde{\alpha}_{n}} (\ln t) \alpha^{\prime}(t_{n-\frac{1}{2}})\left(t-t_{n-\frac{1}{2}}\right)  \\
&+\frac{1}{2}\left[t^{1-\alpha(\eta_1)}(\ln t)^{2}\left(\alpha^{\prime}\left(\eta_{1}\right)\right)^{2}-t^{1-\alpha(\eta_1)} (\ln t) \alpha^{\prime \prime}(\eta_1)\right]\left(t-t_{n-\frac{1}{2}}\right)^{2},\\
	\end{aligned}
\end{equation*}
and
\begin{equation*}
	\begin{aligned}
\Gamma(2-\alpha(t))
=&
\Gamma \left(2-\tilde {\alpha}_{n}\right)-\Gamma^{\prime}\left(2-\tilde{\alpha}_{n}\right) {\alpha}^{\prime}(t_{n-\frac{1}{2}})\left(t-t_{n-\frac{1}{2}}\right)\\
&+
\frac{1}{2}\left[-\Gamma^{\prime \prime}\left(2-\alpha\left(\eta_{2}\right)\right)\left(\alpha^{\prime}\left(\eta_{2}\right)\right)^{2}
+
\Gamma^{\prime}\left(2-\alpha\left(\eta_{2}\right)\right) \alpha^{\prime \prime}\left(\eta_{2}\right)\right]\left(t-t_{n-\frac{1}{2}}\right)^{2},
	\end{aligned}
\end{equation*}
where $\eta_{1}$ and $\eta_{2}$ are both between $t$ and $t_{n-\frac{1}{2}}$. 
Thus, we have
\begin{equation*}
\begin{aligned}
\left|r_{1, n}\right| 
&\leq 
Q_{0}Q_{1}Q_{2}\left|\frac{1}{\tau} \int_{t_{n-1}}^{t_{n}} \left[ t^{1-\alpha(t)}-t^{1-\tilde{\alpha}_{n}}\right] d t\right|
+
Q_{0}Q_{1}Q_{3}\left|\frac{1}{\tau} \int_{t_{n-1}}^{t_{n}} \big[\Gamma(2-\alpha(t))-\Gamma\left(2-\tilde{\alpha}_{n}\right) \big]d t\right|		\\
&\leq 
Q_{0}Q_{1}Q_{2}\left|\frac{1}{\tau} \int_{t_{n-1}}^{t_{n}} \left[Q_{4}\left(t-t_{n-\frac{1}{2}}\right)+Q_{5}\left(t-t_{n-\frac{1}{2}}\right)^{2} \right]d t\right| \\
&\quad +
Q_{0}Q_{1}Q_{3}\left|\frac{1}{\tau} \int_{t_{n-1}}^{t_{n}} \left[Q_{6}\left(t-t_{n-\frac{1}{2}}\right)+Q_{7}\left(t-t_{n-\frac{1}{2}} \right)^{2} \right] dt\right|,
\end{aligned}
\end{equation*}
Since 
     $\int_{t_{n-1}}^{t_{n}} \left(t-t_{n-\frac{1}{2}}\right)dt=0$ and
     $\int_{t_{n-1}}^{t_{n}} \left(t-t_{n-\frac{1}{2}} \right)^{2} dt = \tau^{3}$,
we can get $\left|r_{1, n}\right|  = \mo\left(\tau^{2}\right)$.
\section{Estimate of $r_{2,n}$} \label{r2n}
Similar to the routine in Refs. \cite{L1+1, WangL1+}, we can get the following proof.
By exchanging the order of integration and using integration by parts on each subinterval,
we can get
\begin{equation} \label{B-1}
	\begin{aligned}
			r_{2,n} 
			& =\frac{1}{\tau}\left[\int_0^{t_n} \int_0^t \omega_{1-\tilde {\alpha}_{n}}(t-s) \partial_s \theta v(s) dsdt-\int_0^{t_{n-1}} \int_0^t \omega_{1-\tilde {\alpha}_{n}}(t-s) \partial_s \theta v(s) dsdt\right] \\
			& =\frac{1}{\tau} \int_{t_0}^{t_n} \omega_{1-\tilde {\alpha}_{n}}\left(t_n-s\right) \theta v\left(s\right) ds-\frac{1}{\tau} \int_{t_0}^{t_{n-1}} \omega_{1-\tilde {\alpha}_{n}}\left(t_{n-1}-s\right) \theta v\left(s\right)ds.
	\end{aligned}
\end{equation}

(I) For $1\leq n\leq 3$, we have from (\ref{theta-v}) and (\ref{B-1}) that
\begin{equation} \label{B-2}
	\begin{aligned}
	\left|r_{2,n}\right| 
	& \leq 
	\frac{1}{\tau} \sum_{k=1}^{n-1} \int_{t_{k-1}}^{t_k}
	\left[\omega_{1-\tilde {\alpha}_{n}}\left(t_{n-1}-s\right)-\omega_{1-\tilde {\alpha}_{n}}\left(t_n-s\right)\right]
	\left|\theta v\left(s\right)\right| ds
	+\frac{1}{\tau} \int_{t_{n-1}}^{t_n} \omega_{1-\tilde {\alpha}_{n}}\left(t_n-s\right)\left|\theta v\left(s\right)\right| ds \\
	& 
	\leq 
	Q_8 \sum_{k=1}^{n-1}\left(t_k^{1-\alpha(0)}-t_{k-1}^{1-\alpha(0)}\right) \tau^{1-\tilde {\alpha}_{n}}
	+Q_9\left(t_n^{1-\alpha(0)}-t_{n-1}^{1-\alpha(0)}\right) \tau^{1-\tilde {\alpha}_{n}}
	=Q_{10} t_{n}^{1-\alpha(0)} \tau^{1-\tilde {\alpha}_{n}} 
	\leq Q_{11} t_n^{-\tilde {\alpha}_{n}-\alpha(0)} \tau^2 .
	\end{aligned}
\end{equation}

(II) For $n\geq 4$, we set $n_0=\left\lceil\frac{n}{2}\right\rceil$ so that $\frac{n}{2} \leq n_0 \leq \frac{n}{2}+1$ and $n \geq n_0+2$.
According to (\ref{B-1}), we now split $r_{2,n}= r_{2,n}^1+r_{2,n}^2+r_{2,n}^3 $, where 
\begin{equation} \label{B-3}
	\begin{aligned}
	r_{2,n}^1 & =\frac{1}{\tau} \sum_{k=1}^{n_0} \int_{t_{k-1}}^{t_k}\left[\omega_{1-\tilde {\alpha}_{n}}\left(t_n-s\right)
	              -\omega_{1-\tilde {\alpha}_{n}}\left(t_{n-1}-s\right)\right] \theta v\left(s\right) ds, \\
	r_{2,n}^2 & =\frac{1}{\tau} \int_{t_{n_{0}}}^{t_{n_{0}+1}} \omega_{1-\tilde {\alpha}_{n}}\left(t_n-s\right) \theta v\left(s\right)ds, \\
	r_{2,n}^3 & =\frac{1}{\tau} \sum_{k=n_0+1}^{n-1}\left[\int_{t_k}^{t_{k+1}} \omega_{1-\tilde {\alpha}_{n}}\left(t_n-s\right) \theta v\left(s\right) ds
	             -\int_{t_{k-1}}^{t_k} \omega_{1-\tilde {\alpha}_{n}}\left(t_{n-1}-s\right) \theta v\left(s\right) ds\right].
	\end{aligned}
\end{equation}
Since $0 \leq \omega_{1-\tilde {\alpha}_{n}}\left(t_{n-1}-s\right)-\omega_{1-\tilde {\alpha}_{n}}\left(t_{n}-s\right) 
\leq Q_{12} \tau\left(t_{n-1}-t_{n_{0}}\right)^{-\tilde {\alpha}_{n}-1} \leq Q_{13} \tau t_{n}^{-\tilde {\alpha}_{n}-1}$ for $t_{0} \leq s \leq t_{n_{0}}$,
we can get from (\ref{theta-v}) that
\begin{equation} \label{B-4}
	\begin{aligned}
	\left|r_{2,n}^1\right| 
	\leq 
	Q_{14} \tau^{2} t_{n}^{-\tilde {\alpha}_{n}-1} \sum_{k=1}^{n_{0}}\left(t_{k}^{1-\alpha(0)}-t_{k-1}^{1-\alpha(0)}\right)
	=
	Q_{14} \tau^{2} t_{n}^{-\tilde {\alpha}_{n}-1} t_{n_{0}}^{1-\alpha(0)} 
	\leq 
	Q_{15} t_{n}^{-\tilde {\alpha}_{n}-\alpha(0)} \tau^{2}.
	\end{aligned}
\end{equation}
Similarly, we can obtain
\begin{equation} \label{B-5}
	\begin{aligned}
    \left|r_{2,n}^2\right| 
    \leq 
    Q_{16} \tau\left(t_{n_{0}+1}^{1-\alpha(0)}-t_{n_{0}}^{1-\alpha(0)}\right)\left(t_{n}-t_{n_{0}+1}\right)^{-\tilde {\alpha}_{n}} 
    \leq
     Q_{17} \tau^{2} t_{n_{0}}^{-\alpha(0)} t_{n}^{-\tilde {\alpha}_{n}} 
     \leq 
     Q_{18} t_{n}^{-\tilde {\alpha}_{n}-\alpha(0)} \tau^{2}.
	\end{aligned}
\end{equation}
According to that
\begin{equation*} 
	\begin{aligned}
	\int_{t_{k}}^{t_{k+1}} \omega_{1-\tilde {\alpha}_{n}}\left(t_{n}-s\right)\left(s-t_{k}\right)\left(s-t_{k+1}\right)ds
	=
	\int_{t_{k-1}}^{t_{k}} \omega_{1-\tilde {\alpha}_{n}}\left(t_{n-1}-s\right)\left(s-t_{k-1}\right)\left(s-t_{k}\right)ds,
	\end{aligned}
\end{equation*}
we can rewrite $r_{2,n}^3$ as $r_{2,n}^3=\frac{1}{\tau} \sum_{k=n_{0}+1}^{n-1}\left(\eta_{n}^{k}-\widetilde{\eta}_{n}^{k}\right)$, where
\begin{equation*} 
	\begin{aligned}
	&\eta_{n}^{k}
	=
	\int_{t_{k}}^{t_{k+1}} \omega_{1-\tilde {\alpha}_{n}}\left(t_{n}-s\right)\left[\theta v\left(s\right)-\frac{\partial_{t}^{2} v\left( t_{k}\right)}{2}\left(s-t_{k}\right)\left(s-t_{k+1}\right)\right]ds, \\
	&\tilde{\eta}_{n}^{k}
	=
	\int_{t_{k-1}}^{t_{k}} \omega_{1-\tilde {\alpha}_{n}}\left(t_{n-1}-s\right)\left[\theta v\left(s\right)-\frac{\partial_{t}^{2} v\left( t_{k}\right)}{2}\left(s-t_{k-1}\right)\left(s-t_{k}\right)\right]ds .
	\end{aligned}
\end{equation*}
When $n_{0}+1 \leq k \leq n-1$ and $t_{k} \leq s \leq t_{k+1}$, $ \theta v\left(s\right)=\frac{\partial_{t}^{2} v\left(\hat{\rho}_{k}\right)}{2}\left(s-t_{k}\right)\left(s-t_{k+1}\right)$, 
where $\hat{\rho}_{k} \in\left(t_{k}, t_{k+1}\right)$.
Thus,  there exists $\hat{\zeta }_{k} \in \left(t_{k}, \hat{\rho}_{k}\right)$ such that 
\begin{equation*} 
	\begin{aligned}
	\left|\eta_{n}^{k}\right| 
	\leq 
	\frac{1}{2} \int_{t_{k}}^{t_{k+1}} \omega_{1-\tilde {\alpha}_{n}}\left(t_{n}-s\right)\left|\partial_{t}^{3} v\left( \hat{\zeta }_{k}\right)\right|\left(\hat{\rho}_{k}-t_{k}\right)\left(s-t_{k}\right)\left(t_{k+1}-s\right) ds 
	\leq 
	Q_{19} \tau^{3} t_{k}^{-1-\alpha(0)} \int_{t_{k}}^{t_{k+1}} \omega_{1-\tilde {\alpha}_{n}}\left(t_{n}-s\right)ds.
	\end{aligned}
\end{equation*}
Similarly, we have
\begin{equation*} 
	\begin{aligned}
	\left|\tilde{\eta}_{n}^{k}\right| 
	\leq 
	Q_{20} \tau^{3} t_{k-1}^{-1-\alpha(0)} \int_{t_{k-1}}^{t_{k}} \omega_{1-\tilde {\alpha}_{n}}\left(t_{n-1}-s\right)ds, \quad n_{0}+1 \leq k \leq n-1.
	\end{aligned}
\end{equation*}
Based on that $t_{k} \geq t_{k-1} \geq t_{n_{0}} \geq \frac{1}{2} t_{n}$ when $n_{0}+1 \leq k \leq n-1$, we can get
\begin{equation} \label{B-6}
	\begin{aligned}
		\left|r_{2,n}^3\right| 
		\leq 
		\frac{1}{\tau} \sum_{k=n_{0}+1}^{n-1}\left|\eta_{n}^{k}-\tilde{\eta}_{n}^{k}\right| 
		\leq 
		Q_{21} \tau^{2} t_{n}^{-1-\alpha(0)}\left[\left(t_{n}-t_{n_{0}+1}\right)^{1-\tilde {\alpha}_{n}}+\left(t_{n-1}-t_{n_{0}}\right)^{1-\tilde {\alpha}_{n}}\right] 
		\leq 
		Q_{22} t_{n}^{-\tilde {\alpha}_{n}-\alpha(0)} \tau^{2}.
	\end{aligned}
\end{equation}
By (\ref{B-4}), (\ref{B-5}) and (\ref{B-6}), we obtain $\left|r_{2,n}\right| \leq Q_{23} t_{n}^{-\tilde {\alpha}_{n}-\alpha(0)} \tau^{2}$ for $n \geq 4$. 
The proof is completed.

\end{document}